\documentclass[12pt]{article}
\usepackage{mathtools,amsfonts,amssymb}
\usepackage[a4paper, margin=1.5cm]{geometry}
\usepackage{subfiles}
\newcommand{\onlyinsubfile}[1]{#1}
\newcommand{\notinsubfile}[1]{}

\everymath{\displaystyle}

\usepackage{url}
\usepackage{faktor}
\usepackage{amsthm}
\usepackage{thmtools}
\usepackage{tikz-cd}
\usepackage{stmaryrd}
\usepackage{comment}
\usepackage{hyperref}

\newtheoremstyle{kkkk}
{\topsep}
{0pt}
{}
{0pt}
{\bfseries}
{:}
{\newline}
{}

\newtheoremstyle{llll}
{\topsep}
{0pt}
{}
{0pt}
{\bfseries}
{}
{\newline}
{}

\newtheoremstyle{pf}
{0pt}
{\topsep}
{}
{0pt}
{\bfseries}
{:}
{\newline}
{}

\theoremstyle{kkkk}
\newtheorem{jthm}{Theorem}[section]

\newtheorem{jdef}[jthm]{Definition}

\newtheorem{jlem}[jthm]{Lemma}

\newtheorem{jprop}[jthm]{Proposition}
\newtheorem{jcor}[jthm]{Corollary}
\newtheorem*{jcorr}{Corollary}

\theoremstyle{pf}
\newtheorem*{jpf}{Proof}
\newcommand{\h}{\hangindent=2em}
\newcommand{\hsmall}{\hangindent=1em}

\theoremstyle{llll}

\newtheorem*{jex}{Example}

\declaretheoremstyle[
  spaceabove=0pt,
  spacebelow=0pt,
  headfont=\bfseries,
  headpunct={},
  postheadspace=\newline,
  notefont=\bfseries,
  notebraces={of~}{:},
]{ofstyle}
\declaretheorem[name={Proof},style=ofstyle, unnumbered]{jpfof}

\newcommand{\jhom}{\mathrm{Hom}}
\newcommand{\Tor}{\textrm{Tor}}
\newcommand{\Ker}{\textrm{Ker }}
\newcommand{\Res}{\textrm{Res}}
\newcommand{\Ind}{\textrm{Ind}}

\newcommand{\jref}[1]{[\textbf{\small{#1}}]}
\newcommand{\jcite}[1]{\notinsubfile{\cite{#1}}\onlyinsubfile{\jref{#1}}}

\begin{document}

\renewcommand{\onlyinsubfile}[1]{}
\renewcommand{\notinsubfile}[1]{#1}

\begin{center}
\LARGE \textbf{Coherence of Augmented Iwasawa Algebras}\vspace{10pt}\\
\large James Timmins \vspace{10pt}
\end{center}

\begin{center}
\large{\textbf{Abstract}} \vspace{-4pt}
\end{center}

The augmented Iwasawa algebra of a $p$-adic Lie group is a generalisation of the Iwasawa algebra of a compact $p$-adic Lie group. We prove that a split-semisimple group over a $p$-adic field has a coherent augmented Iwasawa algebra if and only if its root system is of rank one. We deduce that the general linear group of degree $n$ has a coherent augmented Iwasawa algebra precisely when $n$ is at most two. We also characterise when certain solvable $p$-adic Lie groups have a coherent augmented Iwasawa algebra.


\section*{Introduction}

Given a profinite group $G$, and a commutative ring $R$, one can form the completed group ring $RG$. This is a generalisation of the group ring of a finite group. A commonly studied class of profinite groups is the compact $p$-adic Lie groups. If $G$ is a compact $p$-adic Lie group and $R = k$ is a field of characteristic $p$, we call its completed group ring $kG$ the (mod $p$) Iwasawa algebra. The modules over the Iwasawa algebra correspond to continuous representations of $G$ over $k$. This means the representation theory of compact $p$-adic Lie groups is bound up with the study of Iwasawa algebras and their modules. One of the most important and useful facts in this theory is that Iwasawa algebras are always Noetherian rings. \\ 

The representation theory of non-compact $p$-adic Lie groups, such as $GL_n(\mathbb{Q}_p)$, is of considerable interest, particularly in the context of the Langlands programme. Kohlhaase has defined in \jcite{Koh17} a generalisation of the Iwasawa algebra for a non-compact $p$-adic Lie group $G$, which we call the augmented Iwasawa algebra of $G$. Our choice of terminology recognises that an augmented representation of $G$, defined by Emerton in \jcite{Eme10}, is exactly a module over the augmented Iwasawa algebra of $G$. The augmented Iwasawa algebra of $G$ is faithfully flat over the augmented Iwasawa algebra of any closed subgroup, Theorem \ref{non-compact faithfully flat thm}, and flat if $G$ is a locally profinite group, Proposition \ref{locally profinite flatness prop}.\\

The augmented Iwasawa algebra of $G$ is almost never Noetherian if $G$ is not compact. One might hope that augmented Iwasawa algebras are coherent rings, meaning every finitely-generated ideal is finitely-presented. Indeed, Shotton has shown in \jcite{Sho20} that the augmented Iwasawa algebra of $SL_2(F)$ is coherent, for $F$ a finite extension of $\mathbb{Q}_p$. We show that the augmented Iwasawa algebra of a unipotent group is coherent, see Corollary \ref{unipotent coherent cor}. \\

However, not all augmented Iwasawa algebras are coherent, as we show in this article. We give a characterisation of which solvable $p$-adic Lie groups in a certain class have a coherent mod $p$ augmented Iwasawa algebra, Theorem \ref{Main result general coherence thm}. From this we deduce a characterisation for split-solvable algebraic groups over $F$, Corollary \ref{Main result solvable algebraic group non-coherence cor}, for split-semisimple algebraic groups, Theorem \ref{Main result semisimple algebraic group non-coherence thm}, and for general linear groups, Corollary \ref{Main result general linear group coherence cor}. In particular we show that $SL_2(F)$ and $PGL_2(F)$ are the only split-semisimple algebraic groups with a coherent augmented Iwasawa algebra, and that the augmented Iwasawa algebra of $GL_n(F)$ is coherent precisely when $n \leq 2$.\\

Shotton has shown in \jcite{Sho20} that the category of finitely-presented smooth representations of $G$ is an abelian category if the augmented Iwasawa algebra of $G$ is coherent. Our results raise the question of whether this category remains abelian in the absence of the augmented Iwasawa algebra being coherent. Emerton, Gee, and Hellmann have conjectured that this holds for the general linear groups in \jcite{EmGeHe23}. I do not know the answer in this case, or when $G$ is reductive, or semisimple. However, in subsection \ref{The category of finitely-presented smooth mod $p$ representations} we give a module-theoretic criterion, and answer the question in the negative for some solvable $p$-adic Lie groups.\\

\textit{Acknowledgements.} I thank the anonymous referee for their very careful reading of the manuscript, which greatly improved this article. I thank Konstantin Ardakov for many valuable discussions. This research was financially supported by EPSRC.

\section{Main results} \label{Main results}

Let $p$ be a prime, let $k$ be a perfect field of characteristic $p$, and $F$ be a finite extension of $\mathbb{Q}_p$. In this article we prove that augmented Iwasawa algebras of semisimple groups over $F$ are rarely coherent.

\begin{restatable}{jthm}{semisimplealgebraicgroups} \h \label{Main result semisimple algebraic group non-coherence thm}
Let $\mathbb{G}$ be a split-semisimple affine group scheme defined over $F$. Let $G = \mathbb{G}(F)$. Then the augmented Iwasawa algebra $kG$ is coherent if and only if the rank of the root system of $\mathbb{G}$ is $1$. \\
\end{restatable}

This means that any split-semisimple algebraic group $G = \mathbb{G}(F)$ has a non-coherent augmented Iwasawa algebra, unless $\mathbb{G} = \mathbb{SL}_2$ or $\mathbb{PGL}_2$. A further consequence is that the augmented Iwasawa algebras of general linear groups are rarely coherent.

\begin{restatable}{jcor}{generallineargroups} \h \label{Main result general linear group coherence cor}
The augmented Iwasawa algebra of $GL_n(F)$ is coherent if and only if $n \leq 2$.\\
\end{restatable}

We deduce Theorem \ref{Main result semisimple algebraic group non-coherence thm} from the following corollary, which in turn is a special case of Theorem \ref{Main result general coherence thm} below.

\begin{restatable}{jcor}{solvablealgebraicgroups} \h \label{Main result solvable algebraic group non-coherence cor}
Let $\mathbb{G}$ be a finite-dimensional split-solvable affine group scheme defined over $F$, and $G=\mathbb{G}(F)$. If the root system of $\mathbb{G}$ has rank 2 or greater, then $kG$ is not coherent. \\
\end{restatable}

In Theorem \ref{Main result general coherence thm} we give a characterisation of when certain solvable $p$-adic Lie groups $G$ have a coherent augmented Iwasawa algebra.\\

We describe the class of such groups and state the theorem.\\

Let $n \in \mathbb{N}$. Let $\mathbb{U}_n$ be the affine group scheme of upper unitriangular matrices in $\mathbb{GL}_n$ over $F$. Let $\mathbb{U}$ be a closed affine group subscheme of $\mathbb{U}_n$, defined and split over $F$. Let $U = \mathbb{U}(F)$.\\
Let $T$ be a closed subgroup -- in the $p$-adic topology -- of the diagonal elements $\mathbb{D}_n(F) \leq GL_n(F)$, such that $T$ normalises $U$. Let $G = \langle T, U \rangle \cong T \ltimes U$, so $G$ is a solvable $p$-adic Lie group, of upper triangular elements in $GL_n(F)$. \\ We define a set of roots of $U$ with respect to $T$, similarly to section 8.17 of \jcite{Bo91}.\vspace{2pt}\\
Let $\mathfrak{u} = \textrm{Lie}(U)$. Since $T$ acts on $U$ via conjugation, $T$ acts on $\mathfrak{u}$. Let
\[X(T) = \{ \beta: T \rightarrow F^{\times} \mid \beta \textrm{ a continuous group homomorphism} \} \] 
be the character group of $T$. For each $\beta \in X(T)$, define the weight space
\[ \mathfrak{u}_\beta = \{ v \in \mathfrak{u} \mid t \cdot v = \beta(t)v \textrm{ }\forall t \in T \} \leq \mathfrak{u}, \]
and define the set of weights
\[ \Phi = \{ \beta \in X(T) \mid \mathfrak{u}_\beta \neq 0 \}. \]
Let $\mathbb{T}_{\mathbb{U}}$ be the normaliser of $\mathbb{U}$ in $\mathbb{D}_n$, so $T$ is a subgroup of $\mathbb{T}_{\mathbb{U}}(F)$. By Proposition 8.4 of \jcite{Bo91}, $\mathbb{T}_{\mathbb{U}}(F)$ acts diagonalisably on $\mathfrak{u}$, and hence so does $T$. \\
Thus $\mathfrak{u}$ is the direct sum of its weight spaces:
\[ \mathfrak{u} = \bigoplus_{\beta \in \Phi} \mathfrak{u}_\beta .\]
We define the group homomorphism $f$ by
\[ f: T \rightarrow \mathbb{Z}^\Phi, \quad f(t) = \big(v_F( \beta(t) )\big)_{\beta \in \Phi}, \]
where $v_F : F^\times \rightarrow \mathbb{Z}$ is the discrete valuation on $F$.\\

\begin{restatable}{jthm}{generalsolvablegroups} \h \label{Main result general coherence thm}
Let $G = T \ltimes U$ be an upper-triangular subgroup of $GL_n(F)$ as described above, with set of roots $\Phi$ and group homomorphism $f: T \rightarrow \mathbb{Z}^\Phi$. The mod $p$ augmented Iwasawa algebra $kG$ is a coherent ring if and only if the image $f(T) \leq  \mathbb{Z}^\Phi$ is a cyclic subgroup, generated by an element of $(\mathbb{Z}_{\geq 0})^\Phi$. \\
\end{restatable}

See Section \ref{Characterisations of coherence} for the proof of Theorem \ref{Main result general coherence thm}. Then, Corollary \ref{Main result solvable algebraic group non-coherence cor}, Theorem \ref{Main result semisimple algebraic group non-coherence thm}, and Corollary \ref{Main result general linear group coherence cor} are deduced in Section \ref{Coherence for algebraic groups}.

\section{Augmented Iwasawa algebras}

\subsection{Definitions} \label{Definitions}

If $G$ is a profinite group, we can define its completed group ring as an inverse limit of group rings of finite groups.

\begin{jdef} \h \label{completed group algebra def}
Let $G$ be a profinite group and $R$ be a commutative ring. The completed group ring is
\[ RG = \varprojlim_{H \trianglelefteq_o G} R\left[\faktor{G}{H}\right], \]
where the limit is taken over the inverse system of open normal subgroups of $G$, with respect to reverse inclusion.\\
\end{jdef}

In the case when $p$ is a prime, $G$ a compact $p$-adic Lie group, and $k$ is a field of characteristic $p$, the completed group ring $kG$ is known as the (mod $p$) Iwasawa algebra of $G$. \\

Kohlhaase has extended this definition to any locally profinite group, see section 1, page 6 of \jcite{Koh17}. See also section 6.1 of \jcite{Ard21} for a generalised construction. We summarise the definition as the following proposition, which is Proposition 3.2 of \jcite{Sho20}. Recall that a locally profinite group is defined to be a topological group with an open profinite subgroup.

\begin{jprop} \h \label{augmented Iwasawa algebra definition prop}
Let $G$ be a locally profinite group, $R$ be a commutative ring. Let $K \leq G$ be an open profinite subgroup, with completed group ring $RK$. Then there is a unique $R$-algebra structure on
\[ RG = RK \otimes_{R[K]} R[G] \]
such that the natural maps $R[G] \rightarrow RG$ and $RK \rightarrow RG$ are $R$-algebra homomorphisms. This $R$-algebra is independent of the choice of $K$, up to canonical isomorphism. \\
\end{jprop}

This means that $RG$ is generated as an $R$-algebra by $RK$ and the elements of $G$. In the case that $G$ is a $p$-adic Lie group and $k$ is a field of characteristic $p$, we call the $k$-algebra $kG$ the (mod $p$) augmented Iwasawa algebra of $G$. If $G$ is compact, then the Iwasawa algebra and the augmented Iwasawa algebra are isomorphic rings, so our notation is unambiguous.\\

\begin{jex} \hsmall
Suppose $G$ is a direct limit of open profinite subgroups, ${G = \varinjlim_{a \in A} K_a}$, where the direct limit is taken in the category of groups. The functor taking a group to its group $k$-algebra commutes with colimits, such as direct limits, because it is left adjoint to the functor that takes a $k$-algebra to its group of units. Therefore $k[G] = \varinjlim_{a \in A} k[K_a]$.\\
If $K \leq G$ is an open profinite subgroup with $K \leq K_a$, then $kK \otimes_{k[K]} k[K_a]  = kK_a$. Hence
\[ kG = kK \otimes_{k[K]} \Big(\varinjlim_{a \in A} k[K_a]\Big) = \varinjlim_{a \in A} (kK \otimes_{k[K]} k[K_a]) = \varinjlim_{a \in A} kK_a. \]
\end{jex}

\begin{jex} \hsmall
Let $G = \mathbb{Q}_p^\times$ and $k$ be a perfect field of characteristic $p$. Then $G$ has an open profinite subgroup $\mathbb{Z}_p^\times$, and $\mathbb{Q}_p^\times/\mathbb{Z}_p^\times \cong \mathbb{Z}$. It follows that \[ k\mathbb{Q}_p^\times \cong k\mathbb{Z}_p^\times[X, X^{-1}] \cong k[[t]][Y][X,X^{-1}]/(Y^{p-1}), \]
so in fact $k\mathbb{Q}_p^\times$ is a Noetherian ring. Similarly $kF^{\times} \cong k\mathcal{O}_F^\times[X,X^{-1}]$ is Noetherian if $F$ is a finite extension of $\mathbb{Q}_p$.\\
\end{jex}

The nomenclature ``augmented'' comes from the fact that the augmented representations given in Definition 2.1.5 of \jcite{Eme10} are modules over the augmented Iwasawa algebra, see Corollary \ref{augmented rep Noetherian version cor}.

\subsection{Universal property of completed group algebras} \label{Universal property of completed group algebras}

In general, it is useful to consider (augmented) Iwasawa algebras defined over a range of coefficient rings. For example, $\mathbb{F}_p$, its finite extensions and algebraic closure, rings of integers $\mathcal{O}_F$ of a $p$-adic field $F$, profinite commutative rings such as $\hat{\mathbb{Z}}$, or the complete local Noetherian $\mathcal{O}_F$-algebras considered in \jcite{Eme10} and \jcite{Sho20}. The following definition, see the introduction of \jcite{Bru66}, encompasses all of these cases.

\begin{jdef} \h
A pseudocompact ring $R$ is a complete Hausdorff topological ring which has a neighbourhood basis of $0$ consisting of two-sided ideals $I$ with $R/I$ an Artinian ring. \\
\end{jdef}

Throughout this section we will assume that $k$ is a commutative pseudocompact ring. Later on we will specialise to the case of $k$ being a perfect field of characteristic $p$.

\begin{jdef} \h
A pseudocompact $k$-algebra is a complete Hausdorff topological ring $A$, with a continuous ring homomorphism from $k$ to the centre of $A$, and such that $A$ has a neighbourhood basis of $0$ consisting of two-sided ideals $I$ with $A/I$ a finite-length $k$-module.\\
\end{jdef}

The completed group algebra $kG$ of a profinite group is naturally a pseudocompact $k$-algebra (see section 4 of \jcite{Bru66}), and satisfies the following property, which is proved over profinite coefficient rings in Proposition 7.1.2 of \jcite{Wil98}.

\begin{jprop} \h \label{completed group algebra universal property prop}
Let $G$ be a profinite group, $k$ a pseudocompact commutative ring, and $A$ be a pseudocompact $k$-algebra. Let $\phi : G \rightarrow A^{\times}$ be a continuous group homomorphism to the units of $A$. There is a unique continuous $k$-algebra homomorphism $\tilde{\phi} : kG \rightarrow A$ extending $\phi$.\\
\end{jprop}
\begin{jpf} \h
Let $\mathcal{B}$ be a neighbourhood basis of $0$ consisting of two-sided ideals $I$ with $A/I$ a finite-length $k$-module. For such an ideal $I \in \mathcal{B}$, composition of $\phi$ with the natural map $A^\times \rightarrow (A/I)^\times$ gives a continuous group homomorphism
\[ \phi_I: G \rightarrow (A/I)^\times. \]
Because $I$ is open, $N_I = \Ker \phi_I$ is an open normal subgroup of $G$, and by the universal property of group algebras, there is a unique $k$-algebra homomorphism
\[ \tilde{\phi_I'}: k[G/N_I] \rightarrow A/I \]
extending the natural group homomorphism $\phi_I': G/N_I \rightarrow (A/I)^\times$. To check that $\tilde{\phi_I'}$ is continuous, note that the homomorphism from $k$ to the centre of $A$ that makes $A$ a $k$-algebra and $\tilde{\phi_I'}$ agree on the restriction
\[ i_I: k \rightarrow A/I. \]
Because $I$ is open, $A/I$ is discrete and $\alpha_I = \Ker i_I$ is an open ideal of $k$. Then ${k[G/N_I]\alpha_I \subseteq k[G/N_I]}$ is open because $G/N_I$ is finite, and $\Ker \tilde{\phi_I'}$ contains $k[G/N_I]\alpha_I$, hence is open. It follows that $\tilde{\phi_I'}$ is continuous.\\
By composing with the natural homomorphisms $kG \rightarrow k[G/N_I]$, we obtain a system of continuous $k$-algebra homomorphisms
\[ \tilde{\phi_I}: kG \rightarrow A/I, \]
which are compatible under composition with the maps $A/I \rightarrow A/J$ for $I, J \in \mathcal{B}$ with $I \subseteq J$, because $N_I \leq N_J$. Since $A$ is complete, it follows from the universal property of inverse limit that there is a continuous $k$-algebra homomorphism
\[ \tilde{\phi}: kG \rightarrow A \]
extending the $\tilde{\phi_I}$, thus extending $\phi$. Moreover, there is a unique $k$-algebra homomorphism $k[G] \rightarrow A$ extending $\phi$, thus $\tilde{\phi}$ is unique because $k[G]$ is dense in $kG$. $\square$\\
\end{jpf}

This property is enough to define the completed group algebra $kG$ up to canonical isomorphism of pseudocompact $k$-algebras. Moreover, associating the completed group algebra to a profinite group gives a functor.\\

\begin{jprop} \h \label{completed group algebra functor prop}
The mapping $F$ given on objects by $F(G) = kG$, and on morphisms by $F(\phi)= \tilde{\phi}$, is a functor from the category of profinite groups (with continuous group homomorphisms) to the category of pseudocompact $k$-algebras. The functor preserves injectivity and surjectivity of morphisms.\\
\end{jprop}
\begin{jpf} \h
Let $\phi:G_1 \rightarrow G_2$ be a continuous group homomorphism. By Proposition \ref{completed group algebra universal property prop}, there is a unique extension of $\phi: G_1 \rightarrow kG_2^\times$ to a continuous $k$-algebra homomorphism $\tilde{\phi}: kG_1 \rightarrow kG_2$, so $F$ is well-defined. Clearly if $G_2=G_1$ and $\phi=id_{G_1}$, then $\tilde{\phi} = id_{kG_1}$. If $\psi: G_2 \rightarrow G_3$ is another continuous group homomorphism, then $\tilde{\psi} \circ \tilde{\phi}: kG_1 \rightarrow kG_3$ extends $\psi \circ \phi$, and therefore is equal to $\widetilde{\psi \circ \phi}$ by uniqueness. Therefore $F$ is a functor.\\
Now, $\tilde{\phi}$ is constructed by taking the inverse limit of the $k$-algebra homomorphisms
\[ \phi_U: k[G_1/\phi^{-1}(U)] \rightarrow k[G_2/U], \]
as $U$ ranges over the open normal subgroups of $G_2$. If $\phi$ is injective then $\phi_U$ is injective for all $U$, so $\tilde{\phi}$ is injective by left-exactness of the inverse limit.\\
If $\phi$ is surjective, the $\phi_U$ are all surjective. The transition maps that define the inverse limit on each side are all surjective, meaning the strong Mittag-Leffler condition is satisfied and so $\tilde{\phi}$ is surjective, see Definition 4.8.3 and Theorem 4.8.5 of \jcite{EdHa76}. $\square$\\
\end{jpf}

Propositions \ref{completed group algebra universal property prop} and \ref{completed group algebra functor prop} can be used to prove basic facts about maps between completed group algebras. Motivated by this, we will generalise this universal property to augmented Iwasawa algebras and to particular modules.\\

\subsection{Profinite modules for an Iwasawa algebra} \label{Profinite modules for an Iwasawa algebra}

In this subsection, let $H$ be a profinite group.

\begin{jdef} \h \label{profinite H-space def}
A profinite $H$-space is a compact Hausdorff totally disconnected topological space $X$ with a continuous action $H \times X \rightarrow X$, and \textit{finitely many} $H$-orbits.\\
\end{jdef}

Proposition 1.1.7 of \jcite{Wil98} tells us that the topological space $X$ is the inverse limit of its finite discrete quotient spaces, justifying the terminology here.\\

We impose the condition of having finitely many orbits because we will be most interested in the case where $H$ acts transitively on $X$. (When a profinite group acts on a profinite space continuously with infinitely many orbits, the notion of \textit{locally profinite} in Definition \ref{locally profinite def} may instead be applied.)\\

We now show that the topology on $X$ can be defined by an inverse limit of finite $H$-spaces.

\begin{jlem} \h \label{inverse limit of orbits lemma}
In the category of topological spaces,
\[ X \cong \varprojlim_{U \trianglelefteq_o H} \mathrm{Orb}_U(X), \]
where $\mathrm{Orb}_U(X)$ is the space of $U$-orbits in $X$. Moreover, the isomorphism is $H$-equivariant with respect to the natural $H$-actions.\\
\end{jlem}
\begin{jpf} \h
We show that $\{U \cdot x \mid U\trianglelefteq_o H, x \in X \}$ is a basis for the topology on $X$.\\
Let $U$ be an open normal subgroup of $H$. Let $Z \subseteq X$ be an $H$-orbit, and consider its $U$-orbits $\mathrm{Orb}_U(Z)$. Because $H$ acts transitively on $Z$, the finite group $H/U$ acts transitively on $\mathrm{Orb}_U(Z)$, and therefore $\mathrm{Orb}_U(Z)$ is finite. Because $X$ has finitely many $H$-orbits, it follows that $\mathrm{Orb}_U(X)$ is finite. Now, $U$ is profinite and hence compact. Thus for any $x \in X$, the $U$-orbit $U \cdot x$ is a continuous image of a compact set, so is compact. But $X$ is Hausdorff, and therefore $U \cdot x$ is a closed subset of $X$. Because the $U$-orbits are closed and give a finite partition of $X$, it follows that any $U$-orbit is also open.\\
Now, let $Y \subseteq X$ be an open subset and $x \in Y$. Let $f : H \times \{x \} \rightarrow X$ be the restriction of the map defining the $H$-action. Then $f^{-1}(Y) = \{ (h,x) \mid h \cdot x \in Y \}$ is open, and hence $V = \{h \in H \mid h \cdot y \in Y\} \subseteq H$ is open. But $V$ contains the identity of $H$, and therefore $V$ contains an open subgroup of $H$, hence contains an open normal subgroup $U$ because $H$ is profinite. Then $U \cdot x \subseteq Y$. It follows that $Y$ is a union of open sets of this form, and so $\{U \cdot x \mid U\trianglelefteq_o H, x \in X \}$ forms a basis of $X$.\\
Moreover, suppose $Y \subseteq X$ is both open and closed. Then $Y$ is compact, so it is a finite union of open sets in the above basis,
\[ Y = \bigcup_{j=1}^n U_j \cdot y_j. \]
Let $U = \bigcap_{j=1}^n U_j \trianglelefteq_o H$. Then $Y$ is a union of finitely many $U$-orbits, which must all be disjoint.\\
Now suppose we have a partition of $X$ into subsets that are both open and closed,
\[ X = \bigcup_{j=1}^n Y_j. \]
By the above, there exist $U_j \trianglelefteq_o H$ such that $Y_j$ is a union of $U_j$-orbits. Letting $U = \bigcap_{j=1}^n U_j$, we have that each $Y_j$ is a disjoint union of $U$-orbits, and $U \trianglelefteq_o H$.\\
By Proposition 1.1.7 of \jcite{Wil98}, $X$ is the inverse limit of its finite discrete quotient spaces $X/ \sim$. Given such a quotient space, let $x_1, \dots, x_n$ be representatives of the equivalence classes of $\sim$, and $Y_j = \Pi^{-1}(\{\Pi(x_j)\})$, where $\Pi:X \rightarrow X/ \sim$ is the quotient map. Since $X/ \sim$ is discrete, each $Y_j$ is both open and closed, so by the above paragraph there exists $U \trianglelefteq_o H$ such that each $Y_j$ is a disjoint union of $U$-orbits. The orbit space $\mathrm{Orb}_U(X)$ is a quotient space of $X$ which is finite, and discrete, since orbits are open. Thus we have a commuting diagram of continuous maps
\[ \begin{tikzcd}
                                       & \mathrm{Orb}_U(X) \arrow[d, "q_U"] \\
X \arrow[ru, "\Pi_U"] \arrow[r, "\Pi"] & X/ \sim                           
\end{tikzcd} \]
where $\Pi_U, q_U$ are the natural quotient maps. So the inverse system of orbit spaces $\mathrm{Orb}_U(X)$ is cofinal in the inverse system of all discrete finite quotients, so by Proposition 1.1.7 of \jcite{Wil98}, we have that $X \cong \varprojlim_{U \trianglelefteq_o H} \mathrm{Orb}_U(X)$.\\
Moreover, each of the natural maps $\Pi_U : X \rightarrow \mathrm{Orb}_U(X)$, $\Pi_{UV}: \mathrm{Orb}_V(X) \rightarrow \mathrm{Orb}_U(X)$, for $V \leq U$ open normal subgroups of $H$, are easily seen to be $H$-equivariant. Hence the above topological isomorphism is $H$-equivariant with respect to the natural $H$-actions. $\square$\\
\end{jpf}

The above description of a profinite $H$-space $X$ motivates us to define a ``completed module'' associated to $X$.

\begin{jdef} \h
Let $X$ be a profinite $H$-space. The completed module of $X$ is the (left) $kH$-module
\[ kX = \varprojlim_{U \trianglelefteq_o H} k[\mathrm{Orb}_U(X)]. \]
\end{jdef}

The $k$-module $kX$ is indeed a $kH$-module because each $k[\mathrm{Orb}_U(X)]$ is a left $k[H/U]$-module, $kH$ is the inverse limit of these $k$-algebras, and all morphisms are $H$-equivariant. Moreover, the $k[H]$-module $k[X]$ naturally embeds into $kX$, and is easily seen to be dense by definition of the inverse limit topology.

\begin{jdef} \h
Let $A$ be a pseudocompact $k$-algebra, $M$ be a (left) topological $A$-module. $M$ is pseudocompact if and only if it is an inverse limit of discrete $A$-modules of finite length.\\
\end{jdef}

\begin{jlem} \h
Let $X$ be a profinite $H$-space. Then $kX$ is a pseudocompact $kH$-module.\\
\end{jlem}
\begin{jpf} \h
Recall that $X$ has finitely many $H$-orbits. Thus for each $U\trianglelefteq_o H$, the $kH$-module $k[\mathrm{Orb}_U(X)]$ is a free $k$-module of finite rank. Let $\mathcal{B}$ be a basis of open ideals $I$ of $k$ with $k/I$ Artinian. Then
\[ kX = \varprojlim_{U \trianglelefteq_o H} k[\mathrm{Orb}_U(X)] \cong \varprojlim_{U \trianglelefteq_o H} \varprojlim_{I \in \mathcal{B}} (k/I)[\mathrm{Orb}_U(X)] \cong \varprojlim_{U \trianglelefteq_o H, I \in \mathcal{B}} (k/I)[\mathrm{Orb}_U(X)] \]
and each $(k/I)[\mathrm{Orb}_U(X)]$ is of finite length as a $k$-module, hence is a (discrete) $kH$-module of finite length. Therefore $kX$ is a pseudocompact $kH$-module. $\square$\\
\end{jpf}

\begin{jlem} \h \label{completed module basic properties lemma}
1. For $X$ a profinite $H$-space, $kX \cong \bigoplus_{Z \in \mathrm{Orb}_H(X)} kZ$.\\
2. Let $J \leq H$ be an open subgroup, and $\Res_{J}^HX$ be $X$ considered as a profinite $J$-space. There is an isomorphism of $kJ$-modules $i_{J,H}: k(\Res_{J}^HX) \rightarrow \Res_{kJ}^{kH} kX$.\\
\end{jlem}
\begin{jpf} \h
1. For each $U \trianglelefteq_o H$, we have that
\[ \mathrm{Orb}_U(X) = \bigsqcup_{Z \in \mathrm{Orb}_H(X)} \mathrm{Orb}_U(Z), \]
so there is a $k[H/U]$-module isomorphism,
\[ k[\mathrm{Orb}_U(X)] \cong \bigoplus_{Z \in \mathrm{Orb}_H(X)} k[\mathrm{Orb}_U(Z)]. \]
By assumption, $X$ has finitely many $H$-orbits, so this direct sum is also a direct product. Since direct products commute with limits, taking inverse limits gives the result.\\
2. Since $J$ is open, it has finite index in $H$, and $X$ has finitely many $H$-orbits, so $X$ has finitely many $J$-orbits. So $X$ is indeed a profinite $J$-space. For a fixed $V \trianglelefteq_o H$ with $V \leq J$, the inverse system $(V \cap U)_{U \trianglelefteq_o H}$ is cofinal in both $(U)_{U \trianglelefteq_o H}$ and $(U')_{U' \trianglelefteq_o J}$. By considering the natural $k[J/V\cap U]$-module homomorphism
\[ k[\mathrm{Orb}_{V \cap U}(\Res_{J}^HX)] \rightarrow \Res_{k[J/V\cap U]}^{k[H/V\cap U]} k[\mathrm{Orb}_{V \cap U}(X)], \]
and taking inverse limits, we obtain the desired isomorphism. $\square$ \\
\end{jpf}

We now prove a universal property of these modules which is analogous to Proposition \ref{completed group algebra universal property prop}.

\begin{jprop} \h \label{completed module universal property prop}
Let $X$ be a profinite $H$-space and $M$ be a pseudocompact $kH$-module. Let $\phi : X \rightarrow M$ be a continuous $H$-equivariant map. There is a unique extension of $\phi$ to a continuous $kH$-module homomorphism $\tilde{\phi}: kX \rightarrow M$.\\
\end{jprop}
\begin{jpf} \h
Since $M$ is a pseudocompact $kH$-module, let $M = \varprojlim_{i \in I} M_i$ where each $M_i$ is a discrete finite-length $kH$-module. The image of the natural $H$-equivariant continuous map $\phi_i: X \rightarrow M_i$ is both discrete and compact, hence finite. Let $U_i \trianglelefteq H$ be the kernel of the natural continuous group homomorphism
\[H \rightarrow \mathrm{Sym}(\phi_i(X)), \]
so $U_i$ is an open normal subgroup of $H$. It follows that $\phi_i$ has a reduction $\bar{\phi_i}: \mathrm{Orb}_{U_i}(X) \rightarrow M_i$. This naturally extends to a continuous $k[H/U_i]$-module homomorphism $\bar{\phi_i}': k[\mathrm{Orb}_{U_i}(X)] \rightarrow M_i$, giving a continuous $kH$-module homomorphism $\tilde{\phi_i}: kX \rightarrow M_i$. By construction, the $\tilde{\phi_i}$ agree under the transition morphisms on the $M_i$. Thus we define
\[ \tilde{\phi}: kX \rightarrow M = \varprojlim_{i \in I} M_i \]
to be the continuous $kH$-module homomorphism constructed from the $\tilde{\phi_i}$.\\
To show uniqueness, note that there is a unique $k$-module homomorphism $k[X] \rightarrow M$ extending $\phi$. Since $k[X]$ is dense in $kX$, it follows there is a unique continuous extension $kX \rightarrow M$. $\square$ \\
\end{jpf}

As for completed group algebras, it follows that the construction of the completed module gives a functor.\\

\begin{jprop} \h \label{completed module functor prop}
The mapping $F$ given on objects by $F(X) = kX$, and on morphisms by $F(\phi)= \tilde{\phi}$, is a functor from the category of profinite $H$-spaces (with continuous $H$-equivariant maps) to the category of pseudocompact $kH$-modules. The functor preserves injectivity and surjectivity of morphisms.\\
\end{jprop}
\begin{jpf} \h
Let $\phi:X_1 \rightarrow X_2$ be a continuous $H$-equivariant map. By Proposition \ref{completed module universal property prop}, there is a unique extension of $\phi: X_1 \rightarrow kX_2$ to a continuous $kH$-module homomorphism $\tilde{\phi}: kX_1 \rightarrow kX_2$, so $F$ is well-defined. Clearly if $X_2=X_1$ and $\phi=id_{X_1}$, then $\tilde{\phi} = id_{kX_1}$. If $\psi: X_2 \rightarrow X_3$ is another continuous $H$-equivariant map, then $\tilde{\psi} \circ \tilde{\phi}: kX_1 \rightarrow kX_3$ extends $\psi \circ \phi$, and therefore is equal to $\widetilde{\psi \circ \phi}$ by uniqueness. Therefore $F$ is a functor.\\
Now, $\tilde{\phi}$ is constructed by taking the inverse limit of the $kH$-algebra homomorphisms
\[ \phi_U: k[\mathrm{Orb}_U(X_1)] \rightarrow k[\mathrm{Orb}_U(X_2)], \]
as $U$ ranges over the open normal subgroups of $H$. If $\phi$ is injective then $\phi_U$ is injective for all $U$, so $\tilde{\phi}$ is injective by left-exactness of the inverse limit.\\
If $\phi$ is surjective, the $\phi_U$ are all surjective. The transition maps that define the inverse limit on each side are all surjective, meaning the strong Mittag-Leffler condition is satisfied and so $\tilde{\phi}$ is surjective, see Definition 4.8.3 and Theorem 4.8.5 of \jcite{EdHa76}. $\square$\\
\end{jpf}

\subsection{Universal property of augmented Iwasawa algebras} \label{Universal property of augmented Iwasawa algebras}

Motivated by subsections \ref{Universal property of completed group algebras} and \ref{Profinite modules for an Iwasawa algebra}, we now prove a similar universal property for the augmented Iwasawa algebra. First, we must define an appropriate topology.

\begin{jdef} \h
Let $G$ be a locally profinite group, with open profinite subgroup $H$. Then
\[ kG = \bigoplus_{g \in H \backslash G} kH \otimes g. \]
The topology on $kG$ is given by the direct sum topology for the topological abelian groups $kH \otimes g$, where the topology on $kH \otimes g$ is defined by requiring the natural right multiplication map $kH \rightarrow kH \otimes g$ to be a homeomorphism.\\
\end{jdef}

A description of the direct sum topology on abelian topological groups can be found in \jcite{ChaDo03}, Proposition 5, and a more special case in section 5.E of \jcite{Schn02}.\\
Then, we can show $kG$ has the following properties.

\begin{jprop} \h \label{properties of kG prop}
1. The topology on $kG$ is independent of choice of open profinite subgroup $H \leq G$, and is Hausdorff.\\
2. The natural map $G \rightarrow kG$ is a homeomorphism onto its image.\\
3. If $K \leq G$ is an open profinite subgroup, then $K$ lies in the group of units of a pseudocompact $k$-subalgebra of $kG$.\\
4. The subring $k[G]$ is dense in $kG$.\\
\end{jprop}
\begin{jpf} \h
1. Let $H' \leq H \leq G$ be open profinite subgroups of $G$. We first show that
\[ \bigoplus_{h \in H'\backslash H} kH' \otimes h \cong kH \]
is a topological isomorphism, where the the left hand side has the direct sum topology and the right hand side the inverse limit topology.\\
Now, the inverse system $X = \{ U \mid U \trianglelefteq_o H, U \trianglelefteq_o H' \}$ (under reverse inclusion) is cofinal in both $\{U \trianglelefteq_o H\}$ and $\{ U \trianglelefteq_o H'\}$. This follows from the fact that if $U$ is normal in $H'$, then the group
\[V = \bigcap_{h \in H} hUh^{-1} \]
is normal in $H$, and is a finite intersection of conjugates of $U$, since $H'$ has finite index in $H$. Therefore we have topological isomorphisms
\[ kH \cong \varprojlim_{U \in X} k\left[\faktor{H}{U}\right], \quad kH' \cong \varprojlim_{U \in X}k\left[\faktor{H'}{U}\right]. \]
Then, we have topological isomorphisms
\[ \bigoplus_{h \in H' \backslash H} kH' \otimes h \cong \bigoplus_{h \in H' \backslash H} \bigg( \varprojlim_{U \in X}k\left[\faktor{H'}{U}\right] \otimes h \bigg) \cong  \varprojlim_{U \in X} \bigg(\bigoplus_{h \in H' \backslash H}k\left[\faktor{H'}{U}\right] \otimes h \bigg), \]
since the direct sums are finite. Then consider the natural isomorphism
\[ \bigoplus_{h \in H' \backslash H}k\left[\faktor{H'}{U}\right] \otimes h \cong k\left[\faktor{H}{U}\right]. \]
The left hand side has the direct sum topology, which coincides with the product topology because the direct sum is finite. Since $k[H'/U]$, $k[H/U]$ are free $k$-modules of finite rank with the corresponding product topology, it follows that the above isomorphism is a topological isomorphism. Thus
\[ \bigoplus_{h \in H' \backslash H} kH' \otimes h \cong \varprojlim_{U \in X}k\left[\faktor{H}{U}\right] \cong kH \]
is a topological isomorphism.\\
It then follows that the topologies on $kG$ defined by $H', H$ are equivalent, due to the isomorphism
\[ \bigoplus_{g \in H \backslash G} kH \otimes g \cong \bigoplus_{g \in H \backslash G}  \bigoplus_{h \in H'\backslash H} kH' \otimes hg \cong \bigoplus_{g \in H' \backslash G} kH' \otimes g. \]
For arbitrary $H, H' \leq G$ open profinite, we obtain an isomorphism by considering $H \cap H' \leq G$ and applying the result twice.\\
The topology on $kG$ is Hausdorff since the topology on each $kH \otimes g$ is Hausdorff. \vspace{2pt}\\
2. First we show this result in the profinite case: let $H$ be a profinite group, and consider the natural map $i=i_H: H \rightarrow kH$. The algebra $kH$ has the inverse limit topology with projection maps
\[ \Pi_{U,I}: kH \rightarrow \varprojlim_{U \trianglelefteq_o H, I \in \mathcal{B}} (k/I)[H/U]. \]
Here $\mathcal{B}$ is an open basis of ideals of $k$ with each $k/I$ Artinian. The map $i$ is continuous if and only if $\phi_{U,I}= \Pi_{U,I} \circ i$ is continuous for all $U,I$. Now, the image of $\phi_{U,I}: H \rightarrow (k/I)[H/U]$ is $H/U$, a finite set. Thus for any subset $X \subseteq (k/I)[H/U]$, the inverse image $\phi_{U,I}^{-1}(X)$ is a union of finitely many cosets of $U$ in $H$, which is open. Hence, $\phi_{U,I}$ is continuous. So, $i$ is continuous.\\
Moreover, $i$ is clearly injective, whilst $H$ is compact and $kH$ Hausdorff. Thus $i: H \rightarrow i(H)$ is a homeomorphism.\\
Now, let $G$ be a locally profinite group, and consider the natural map
\[ i_G: G \rightarrow kG = \bigoplus_{g \in H \backslash G} kH \otimes g.\]
Clearly $i_G(x) = i_H(xg^{-1})\otimes g$ when $x \in Hg$, for $g \in H \backslash G$. So, for each $g \in H \backslash G$, the restriction
\[ j_g={i_G}|_{Hg}: Hg \rightarrow i_H(H) \otimes g \]
is a homeomorphism. Clearly, the cosets $Hg$ form a disjoint open cover of $G$, whilst the sets $i_H(H) \otimes g$ are also disjoint. We show that each $i_H(H) \otimes g$ is also open. For a fixed $g' \in H \backslash G$, let $U_{g'} =kH$. When $g' \not \in Hg$, let $U_g \subseteq kH$ be a proper two-sided ideal of $kH$ -- which exists since $kH$ is a pseudocompact ring. In this case, $U_g \cap i_H(H)$ is empty. Let
\[ U = \bigoplus_{g \in H\backslash G} U_g \otimes g \leq kG, \]
so $U$ is an open subgroup of $kG$. Therefore $U \cap i_G(G) = i_H(H) \otimes g'$ is open in $i_G(G)$. Thus $\{ i_H(H) \otimes g \mid g \in H \backslash G \}$ is a disjoint open cover of $i_G(G)$.\\
The restrictions $j_g={i_G}|_{Hg}$ above give homeomorphisms between corresponding parts of the disjoint open covers of $G$ and of $i_G(G)$, from which it follows straightforwardly that $i_G$ is a homeomorphism onto its image. \vspace{2pt}\\
3. This is trivial, as $K \leq (kK)^{\times}$ and $kK$ is a pseudocompact subalgebra of $kG = kK \otimes_{k[K]}k[G]$. \vspace{2pt}\\
4. We have that 
\[ k[G] = \bigoplus_{g \in H \backslash G} k[H] \otimes g \leq kG, \]
and $k[H]$ is dense in $kH$. If $C$ is a closed set containing $k[G]$, then for all $g \in H \backslash G$, $C \cap (kH \otimes g)$ is closed in $kH \otimes g$, but contains $k[H] \otimes g$. Hence $C \cap (kH \otimes g) = kH \otimes g$, and hence $C = kG$. Thus, $k[G]$ is dense in $kG$. $\square$ \\
\end{jpf}

To prove a universal property for $kG$ similar to that of Proposition \ref{completed group algebra universal property prop}, we will need a class of algebras to replace ``pseudocompact $k$-algebras'', and which contains all algebras $kG$ where $G$ is locally profinite. For this reason, and using the proposition just proved, we define the following class of algebras.

\begin{jdef} \h \label{Iwasawa-like algebra def}
Let $A$ be a topological $k$-algebra, that is, a topological ring with a continuous homomorphism from $k$ to the centre of $A$, and assume $A$ is Hausdorff. Let $G \leq A^\times$ be a subgroup of the units of $A$. The subgroup $G$ makes $A$ Iwasawa-like if and only if the subspace topology makes $G$ a topological group (that is, the inversion map on $G$ is continuous), and there is an open profinite subgroup $K \leq G$ and a pseudocompact $k$-subalgebra $B \leq A$ such that $K$ is a subgroup of $B^\times$.\\
\end{jdef}

If $G$ is a locally profinite group, the subgroup $G \leq (kG)^{\times}$ makes $kG$ an Iwasawa-like topological $k$-algebra, by the second and third statements of Proposition \ref{properties of kG prop}. Moreover, the ring of $k$-linear endomorphisms of a smooth representation is also Iwasawa-like.

\begin{restatable}{jprop}{smoothrepsprop} \h \label{endomorphisms of smooth rep is Iwasawa-like prop}
Let $V$ be a smooth representation of $G$. The ring of $k$-linear endomorphisms of $V$ is naturally an Iwasawa-like algebra.\\
\end{restatable}

This is proved in subsection \ref{Smooth representations of G}, where our notion of smooth representation is also defined. We can prove the following universal property for augmented Iwasawa algebras.

\begin{jprop} \h \label{universal property prop}
Let $G$ be a locally profinite group. Let $A$ be an Iwasawa-like topological $k$-algebra, via subgroup $L \leq A^\times$. Let $\phi: G \rightarrow L$ be a continuous group homomorphism. There is a unique continuous $k$-algebra homomorphism $\tilde{\phi} : kG \rightarrow A$ extending $\phi$.\\
\end{jprop}

An application of Proposition \ref{universal property prop} to the result of Proposition \ref{endomorphisms of smooth rep is Iwasawa-like prop} gives the following result on smooth representations.

\begin{restatable}{jcor}{smoothrepscor} \h \label{smooth reps are kG-modules cor}
Let $V$ be a smooth representation of $G$. There is a unique $kG$-module structure extending the $k[G]$-module action on $V$.\\
\end{restatable}

The proof of Corollary \ref{smooth reps are kG-modules cor} can also be found in subsection \ref{Smooth representations of G}. We now prove Proposition \ref{universal property prop}, which requires the following weaker statement.

\begin{jlem} \h \label{profinite universal property lemma}
Let $H$ be a profinite group. Let $A$ be an Iwasawa-like topological $k$-algebra, via subgroup $L \leq A^\times$. Let $\phi: H \rightarrow L$ be a continuous group homomorphism. There is a unique continuous $k$-algebra homomorphism $\tilde{\phi} : kH \rightarrow A$ extending $\phi$.\\
\end{jlem}

\begin{jpfof}[Proposition \ref{universal property prop}] \h
For any open profinite $H \leq G$, the restriction of $\phi: G \rightarrow L$ to $\phi|_H : H \rightarrow L$ is a continuous group homomorphism, and hence by Lemma \ref{profinite universal property lemma}, there is a unique continuous $k$-algebra homomorphism $\psi_H : kH \rightarrow A$ extending $\phi|_H$. There is also a unique $k$-algebra homomorphism $\psi_G : k[G] \rightarrow A$ extending $\phi$, by the universal property of the usual group algebra. Then, the restrictions $\psi_G|_{k[H]} = \psi_H|_{k[H]}$, and $\psi_H|_{kJ} = \psi_J$ if $J \leq H$ are open profinite subgroups of $G$.\\
For any open profinite subgroup $H \leq G$, define the following $k$-bilinear map:
\[ \phi_H': kH \times k[G] \rightarrow A, \quad \phi_H'(x,y) = \psi_H(x)\psi_G(y). \]
The map $\phi_H'$ is also $(kH, k[G])$-bilinear, where $kH, k[G]$ act on $A$ via $\psi_H, \psi_G$. Moreover if $r \in k[H]$, then because $\psi_G, \psi_H$ are ring homomorphisms and agree on $k[H]$, 
\[ \phi_H'(xr, y) = \psi_H(xr)\psi_G(y) = \psi_H(x)\psi_H(r)\psi_G(y) = \psi_H(x) \psi_G(r)\psi_G(y) = \psi_H(x)\psi_G(ry) = \phi_H'(x, ry). \]
Thus, using the universal property of the tensor product, $\phi_H'$ induces a homomorphism of $(kH, k[G])$-bimodules
\[ \tilde{\phi}_H: kH \otimes_{k[H]} k[G] \rightarrow A, \]
where $\tilde{\phi}_H(x \otimes y) = \psi_H(x)\psi_G(y)$. \vspace{4pt}\\
We now show that the homomorphisms $\tilde{\phi}_H$ for varying $H$ agree, under appropriate identifications.\\
Let $J \leq H$ be open profinite subgroups of $G$. Let $i_{J,H}': kJ \rightarrow kH$ be the natural inclusion of their completed group algebras, and define
\[ i_{J,H}: kJ \otimes_{k[J]}k[G] \rightarrow kH \otimes_{k[H]} k[G], \quad i_{J,H}(x \otimes z) = i_{J,H}'(x) \otimes z. \]
We also define a map $\Pi_{J,H}$ in the opposite direction, in the following way. Let
\[ \alpha : kH \rightarrow kJ \otimes_{k[J]} k[H] \]
be the natural $(kJ, k[H])$-bimodule isomorphism with $\alpha \circ i_{J,H}' = id_{kJ} \otimes 1$. Let $s$ be the natural isomorphism of $(kJ, k[G])$-bimodules,
\[ s: (kJ \otimes_{k[J]} k[H] ) \otimes_{k[H]} k[G]  \rightarrow kJ \otimes_{k[J]} (k[H]  \otimes_{k[H]} k[G]),  \]
and $\beta$ be the natural $(k[H], k[G])$-bimodule isomorphism
\[ \beta: k[H] \otimes_{k[H]} k[G] \rightarrow k[G]. \]
Then we have the maps
\[ kH \otimes_{k[H]} k[G] \xrightarrow{\alpha \otimes 1} (kJ \otimes_{k[J]} k[H] ) \otimes_{k[H]} k[G]  \xrightarrow{s}  kJ \otimes_{k[J]} (k[H] \otimes_{k[H]} k[G]) \xrightarrow{1 \otimes \beta} kJ \otimes_{k[J]} k[G], \]
so let
\[ \Pi_{J,H} = (1 \otimes \beta) \circ s \circ (\alpha \otimes 1). \]
Then $\Pi_{J,H}$ is an isomorphism since $\beta, s, \alpha$ are, and
\begin{align*} \Pi_{J,H} \circ i_{J,H}(x \otimes z) &= \Pi_{J,H}( i_{J,H}'(x) \otimes z)\\ &=(1 \otimes \beta) \circ s ( \alpha( i_{J,H}'(x) \otimes z))\\ &= (1 \otimes \beta) \circ s( (x \otimes 1) \otimes z)\\ &= 1 \otimes \beta (x \otimes (1 \otimes z))\\ &= x \otimes z \end{align*}
for all $x \in kJ, z \in k[G]$. Thus $i_{J,H}$, $\Pi_{J,H}$ are mutually inverse isomorphisms. We will also write $i_{H,J} = \Pi_{J,H}$.\\
Now, the following diagram,
\[ \begin{tikzcd}
{kH \otimes_{k[H]} k[G]} \arrow[r, "\tilde{\phi}_H"]                          & A \\
{kJ \otimes_{k[J]} k[G]} \arrow[ru, "\tilde{\phi}_J"'] \arrow[u, "{i_{J,H}}"] &  
\end{tikzcd}
\]
commutes, because for all $x \in kJ, z \in k[G]$,
\[ \tilde{\phi}_H \circ i_{J,H}(x \otimes z) = \tilde{\phi}_H(i_{J,H}'(x) \otimes z) =  \psi_H(i_{J,H}'(x))\psi_G(z) =\psi_J(x) \psi_G(z) = \tilde{\phi}_J(x \otimes z). \]
This uses that $\psi_H \circ i_{J,H}' = \psi_J$, by uniqueness of the extensions $\psi_J, \psi_H$. So, if $J \leq H$, then ${\tilde{\phi}_H \circ i_{J,H} = \tilde{\phi}_J}$, and so $\tilde{\phi}_H = \tilde{\phi}_J \circ i_{J,H}^{-1} = \tilde{\phi}_J \circ i_{H,J}$ also.\\
For arbitrary open profinite subgroups $H, H' \leq G$, define $i_{H,H'} = i_{H \cap H', H'} \circ i_{H, H \cap H'}$. It then follows that $\tilde{\phi}_H \circ i_{H',H} = \tilde{\phi}_{H'}$. \vspace{4pt}\\
We now use these identifications to show that $\tilde{\phi}_H$ is a ring homomorphism, for any $H \leq G$ open profinite.\\
Let $g \in G, x \in kH$. We consider the product $(1 \otimes g)(x \otimes 1) \in kG = kH \otimes_{k[H]}k[G]$. The multiplication on $kG$ is such that
\[ (1 \otimes g)(x \otimes 1) = gx = i_{K,H}(gxg^{-1} \otimes g), \]
where $K = gHg^{-1}$. Thus
\begin{align*} \tilde{\phi}_H\big( (1 \otimes g)(x \otimes 1) \big) &= \tilde{\phi}_H \circ i_{K,H}(gxg^{-1} \otimes g)\\ &= \tilde{\phi}_K(gxg^{-1} \otimes g)\\ &= \psi_K(gxg^{-1})\psi_G(g). \end{align*}
It remains to determine $\psi_K(gxg^{-1})$. We show that $\psi_K(gxg^{-1}) = \psi_G(g)\psi_H(x) \psi_G(g)^{-1}$. Let
\[ d: kH \rightarrow A, \quad d(t) = \psi_K(gtg^{-1}) - \psi_G(g)\psi_H(t) \psi_G(g)^{-1}. \]
Now, $\psi_K$, $\psi_H$ are continuous and multiplication in $A$ is continuous, so $d$ is a continuous function. If $t \in k[H]$, then $gtg^{-1} \in k[G]$, and 
\[ \psi_K(gtg^{-1}) = \psi_G(gtg^{-1}) = \psi_G(g)\psi_G(t) \psi_G(g)^{-1} = \psi_G(g)\psi_H(t) \psi_G(g)^{-1}, \]
so $d(t)=0$. Thus $d(k[H]) = \{0\}$, but $k[H]$ is dense in $kH$ and $A$ is Hausdorff, so $\{0\}$ is closed and $d(kH) = \{0\}$. Thus,
\[ \tilde{\phi}_H\big( (1 \otimes g)(x \otimes 1) \big) = \psi_G(g)\psi_H(x) \psi_G(g)^{-1} \psi_G(g) = \psi_G(g) \psi_H(x). \]
Because $\tilde{\phi}_H$ is a $(kH, k[G])$-bimodule homomorphism, for all $z \in kH \otimes_{k[H]} k[G]$,
\[ \tilde{\phi}_H(zg) = \tilde{\phi}_H(z)\psi_G(g), \quad \tilde{\phi}_H(yz) = \psi_H(y)\tilde{\phi}_H(z). \]
So, let $x_1, x_2 \in kH$ and $g_1, g_2 \in G$. Then
\begin{align*} \tilde{\phi}_H\big( (x_1 \otimes g_1) (x_2 \otimes g_2) \big) &= \psi_H(x_1) \tilde{\phi}_H\big((1 \otimes g_1) (x_2 \otimes 1) \big) \psi_G(g_2)\\ &= \psi_H(x_1)\psi_G(g_1)\psi_H(x_2)\psi_G(g_2)\\ &= \tilde{\phi}_H(x_1 \otimes g_1) \tilde{\phi}_H(x_2 \otimes g_2). \end{align*}
Therefore $\tilde{\phi}_H$ is a ring (and hence $k$-algebra) homomorphism. Now, $\psi_H$ is continuous, so $\tilde{\phi}_H$ is continuous on $kH \otimes g$, for any $g \in H \backslash G$. Because
\[ \tilde{\phi}_H = \bigoplus_{g \in H \backslash G} \tilde{\phi}_H|_{kH\otimes g} :  \bigoplus_{g \in H \backslash G} kH \otimes g \rightarrow  \bigoplus_{g \in H \backslash G} kH \otimes g, \]
it follows that $\tilde{\phi}_H$ is a continuous $k$-algebra homomorphism extending $\phi : G \rightarrow A^{\times}$. \vspace{4pt}\\
To show uniqueness, recall from Proposition \ref{properties of kG prop} that $k[G]$ is dense in $kG$. Thus any two continuous maps $kG \rightarrow A$ agreeing on $k[G]$ are equal. But $\psi_G :k[G] \rightarrow A$ is the unique $k$-algebra homomorphism extending $\phi$ to $k[G]$, therefore $\tilde{\phi}$ is the unique continuous $k$-algebra homomorphism extending $\phi$. $\square$\\
\end{jpfof}

We now prove Lemma \ref{profinite universal property lemma}.\\

\begin{jpfof}[Lemma \ref{profinite universal property lemma}] \h
We have a continuous group homomorphism $\phi:H \rightarrow L \leq A^\times$. Since $L$ makes $A$ Iwasawa-like, let $K \leq L$ be an open profinite subgroup, and $B \leq A$ be a pseudocompact $k$-subalgebra containing $K$ in its units. Then, $\phi(H) \cap K$ is an open subgroup of $\phi(H)$, and so $\phi^{-1}(\phi(H) \cap K)$ is an open subgroup of $H$. Let $K' \leq H$ be an open normal subgroup of $H$ such that $K' \leq \phi^{-1}(\phi(H) \cap K)$, which exists because $H$ is profinite. So, we have a continuous group homomorphism
\[ \phi|_{K'} : K' \rightarrow K \leq B^\times. \]
By the universal property for completed group algebras in Proposition \ref{completed group algebra universal property prop}, $\phi|_{K'}$ extends uniquely to a continuous $k$-algebra homomorphism $\psi : kK' \rightarrow B$.\\
Using the universal property of ordinary group algebras, let $f: k[H] \rightarrow A$ be the unique $k$-algebra extension of $\phi$. By uniqueness of these maps, we have that $f|_{k[K']} = \psi|_{k[K']}$. Thus, define
\[ \tilde{\phi}: kH \rightarrow A, \quad \tilde{\phi}(x \otimes h) = \psi(x)f(h), \]
where we identify $kH \cong kK' \otimes_{k[K']} k[H]$, and transport the ring structure.\\
Then if $x_1, x_2 \in kK'$ and $h_1, h_2 \in H$,
\begin{align*} \tilde{\phi}\big((x_1 \otimes h_1) (x_2 \otimes h_2) \big) &= \tilde{\phi}\big( (x_1 \otimes 1)(h_1x_2h_1^{-1} \otimes h_1h_2) \big)\\ &= \tilde{\phi}(x_1h_1x_2h_1^{-1} \otimes h_1h_2)\\ &= \psi(x_1h_1x_2h_1^{-1})f(h_1h_2)\\ &= \psi(x_1)\psi(h_1x_2h_1^{-1})f(h_1)f(h_2). \end{align*}
Because $k[K']$ is dense in $kK'$, very similar reasoning as in the proof of Proposition \ref{universal property prop} shows that
\[ \psi(h_1x_2h_1^{-1}) = f(h_1)\psi(x_2)f(h_1)^{-1}, \]
and so
\[ \tilde{\phi}\big((x_1 \otimes h_1) (x_2 \otimes h_2) \big) = \psi(x_1)f(h_1)\psi(x_2)f(h_2) = \tilde{\phi}(x_1 \otimes h_1) \tilde{\phi}(x_2 \otimes h_2) . \]
Thus $\tilde{\phi}$ is a $k$-algebra homomorphism that extends $\phi$. Because $\psi$ is continuous and $kH$ can be identified with $kK \otimes_{k[K]} k[H]$, identical reasoning as in the proof of Proposition \ref{universal property prop} shows that $\tilde{\phi}$ is continuous.\\
Since $k[H]$ is dense in $kH$ and $f$ is the unique $k$-algebra homomorphism $k[H] \rightarrow A$ extending $\phi$, it follows that $\tilde{\phi}$ must be the unique continuous $k$-algebra homomorphism $kH \rightarrow A$ extending $\phi$. $\square$\\
\end{jpfof}

It follows from Proposition \ref{universal property prop} that construction of the augmented Iwasawa algebra gives a functor.

\begin{jprop} \h \label{augmented Iwasawa algebra functor prop} 
The mapping $F$ given on objects by $F(G) = kG$ and on morphisms by $F(\phi) = \tilde{\phi}$, is a functor from the category of locally profinite groups (with continuous group homomorphisms) to the category of Iwasawa-like topological $k$-algebras (with continuous $k$-algebra homomorphisms). The functor preserves surjectivity of morphisms and injectivity of closed morphisms.\\
\end{jprop}
\begin{jpf} \h
That $F$ is a well-defined functor follows directly from Proposition \ref{universal property prop}, using an identical argument as found in the proof of Proposition \ref{completed group algebra functor prop}.\\
Let $\phi: G_1 \rightarrow G_2$ be a continuous group homomorphism between locally profinite groups.\\
Suppose $\phi$ is injective and a closed map. Then $\phi(G_1)$ is a closed subgroup of $G_2$. Let $H_2 \leq G_2$ be an open profinite subgroup, and $H_1 = \phi^{-1}(H_2)$. Then $H_1$ is an open subgroup of $G_1$, and is isomorphic to $\phi(H_1) = \phi(G_1) \cap H_2$ as a topological group, hence is profinite. Now, $\tilde{\phi}|_{kH_1} = \widetilde{\phi|_{H_1}}$ is injective by Proposition \ref{completed group algebra functor prop}, and the natural map induced by $\phi$ on cosets,
\[ H_1 \backslash G_1 \rightarrow H_2 \backslash G_2, \quad H_1g \mapsto H_2\phi(g), \]
is injective. Now, $\tilde{\phi}$ is given by
\[ \tilde{\phi} :  \bigoplus_{g \in H_1 \backslash G_1} kH_1 \otimes g \rightarrow \bigoplus_{g \in H_2 \backslash G_2} kH_2 \otimes g, \quad \tilde{\phi}(x \otimes g) = \tilde{\phi}|_{kH_1}(x) \otimes \phi(g),\]
therefore $\tilde{\phi}$ is injective.\\
Suppose instead $\phi$ is a surjective continuous group homomorphism. Let $H_1 \leq G_1$ be an open profinite subgroup. Because surjective maps between topological groups are open, $\phi(H_1) \leq G_2$ is an open profinite subgroup. By Proposition \ref{completed group algebra functor prop}, $\tilde{\phi}|_{kH_1} = \widetilde{\phi|_{H_1}}$ is surjective, and so 
\[ \tilde{\phi} : \bigoplus_{g \in H_1 \backslash G_1} kH_1 \otimes g \rightarrow \bigoplus_{g \in \phi(H_1) \backslash G_2} k\phi(H_1) \otimes g \]
 is surjective. $\square$\\
\end{jpf}

Proposition \ref{augmented Iwasawa algebra functor prop}, and its proof, implies that if $G' \leq G$ is a closed subgroup, then $kG'$ is identified with the subring of $kG$,
\[ \bigoplus_{ g \in H \backslash HG'} k(G' \cap H) \otimes g, \]
where $H \backslash HG'$ denotes (a set of representatives for) the set $\{Hg \mid g \in G'\} \subseteq H \backslash G$, for $H \leq G$ open profinite. The proposition also implies that if $L$ is a quotient of $G$, $kL$ is a quotient ring of $kG$.

\subsection{Modules for an augmented Iwasawa algebra} \label{Modules for an augmented Iwasawa algebra}

In this subsection, let $G$ be a locally profinite group. We define modules associated to locally profinite $G$-spaces and prove a similar universal property for certain $kG$-modules as in subsection \ref{Profinite modules for an Iwasawa algebra}.\\
Before doing this, it will be useful to prove the characterisation of $kG$-modules as being the augmented representations defined in \jcite{Eme10}.

\begin{jprop} \h \label{augmented rep topology version prop} 
Let $M$ be a $k[G]$-module equipped with a $kH$-module structure for some open profinite subgroup $H \leq G$ such that the two induced $k[H]$-actions coincide. Suppose that $M$ has a Hausdorff topology such that each element of $G$ acts continuously on $M$, and the action $kH \times M \rightarrow M$ is continuous. Then $M$ has a unique $kG$-module structure which extends the action of $k[G]$ and $kH$.\\
\end{jprop}
\begin{jpf} \h
Let $\alpha_H : kH \times M \rightarrow M$ denote the $kH$-module action on $M$ and $\beta: k[G] \times M \rightarrow M$ denote the $k[G]$-module structure. The actions agree on $k[H]$, that is, $\alpha_H|_{k[H] \times M} = \beta|_{k[H] \times M}$.\\
Then the map
\[ \phi_H': kH \times k[G] \times M \rightarrow M, \quad \phi_H'(x, y, m) = \alpha_H(x, \beta(y,m)) \]
is a $k[H]$-balanced map and so we get an induced map 
\[ \phi_H: kH \otimes_{k[H]}k[G] \times M \rightarrow M, \quad \phi_H(x \otimes y, m ) = \alpha_H(x, \beta(y,m)). \]
Note $\phi_H$ extends $\alpha_H$ and $\beta$, and is continuous on $(kH\otimes g) \times M$ for all $g \in G$, thus is continuous.\\
Let $g \in G$ and $m \in M$. Consider the function
\[ d: kH \rightarrow M, \quad d(x) = \phi_H(gxg^{-1}, \beta(g,m)) - \beta(g, \alpha_H(x,m)). \]
Now, $\alpha_H$ is continuous, $G$ acts continuously on $kH$ by conjugation, and $g$ acts continuously on $M$ (via $\beta$). Therefore $d$ is a continuous function. If $x \in k[H]$, we have that $x , gxg^{-1} \in k[G]$, so it follows by the module axioms for $k[G]$ that $d(x)=0$. But, $k[H]$ is dense in $kH$, therefore as $M$ is Hausdorff, we have that $d(kH)=0$.\\
It follows that
\[ \phi_H(gxg^{-1}, \beta(g,m)) = \beta(g, \alpha_H(x,m)), \]
that is,
\[ \phi_H(gxg^{-1}, \phi_H(g,m)) = \phi_H(g, \phi_H(x,m)), \]
for all $x \in kH$, $g \in G$, $m \in M$. Moreover, by definition $\phi_H$ has the property that 
\[ \phi_H(x'zg', m) = \phi_H(x', \phi_H(z, \phi_H(g', m))), \]
for any $x' \in kH$, $z \in kG$, $g' \in G$. Combining these statements, it is straightforward to show
\[ \phi_H((x_1 \otimes g_1)(x_2 \otimes g_2), m) = \phi_H(x_1 (g_1x_2g_1^{-1})g_1g_2, m) = \phi_H(x_1 \otimes g_1, \phi_H(x_2 \otimes g_2, m)), \]
for any $x_1, x_2 \in kH$, $g_1, g_2 \in G$. By $k$-linearity, it follows that for any $z_1, z_2 \in kG$, we have the module axiom
\[ \phi_H(z_1, \phi_H(z_2, m)) = \phi_H(z_1z_2, m). \]
Thus $\phi_H$ defines a module action of $kG = kH \otimes_{k[H]}k[G]$ on $M$. This extends the actions of $kH$ and of $k[G]$, and moreover must be the unique such action because these subrings generate $kG$. $\square$\\
\end{jpf}

\begin{jcor} \h \label{augmented rep Noetherian version cor}
Let $k$ be a profinite commutative ring. Let $M$ be a $k[G]$-module equipped with a $kH$-module structure for some open profinite subgroup $H \leq G$ such that the two induced $k[H]$-actions coincide. If $kH$ is a Noetherian ring, then $M$ has a unique $kG$-module structure that extends the action of $k[G]$ and of $kH$.\\
\end{jcor}
\begin{jpf} \h
Since $k$ is profinite, $kH$ is a profinite Noetherian ring. Thus every finitely-generated $kH$-submodule of $M$ can be given the ``canonical'' topology, see Proposition 2.1.3 of \jcite{Eme10}. We give $M$ the final topology corresponding to the inclusions $N \rightarrow M$ of all finitely-generated submodules $N$ of $M$. That is, $U \subseteq M$ is open if and only if $U \cap N \subseteq N$ is open for every finitely-generated submodule $N$. \\
For $N \leq M$ a finitely-generated $kH$-submodule, the $kH$-action $kH \times N \rightarrow N$ is continuous, again by Proposition 2.1.3 of \jcite{Eme10}. That is, letting $\alpha_H: kH \times M \rightarrow M$ denote the $kH$-module action, the restriction $\alpha_H|_{kH \times N}$ is continuous for any finitely-generated $N$. It is also straightforward to check that a subset $U' \subseteq kH \times M$ is open if and only if $U' \cap (kH \times N) \subseteq kH \times N$ is open for all $N$, using the fact that the inclusions $N \rightarrow N'$ of one finitely-generated submodule to another are open continuous maps, so for any open subset $V \subseteq N$ of any finitely-generated submodule, $V$ is open in $M$.\\
Now, let $U \subseteq M$ be an open subset of $M$. Then
\[ \alpha_H^{-1}(U) \cap (kH \times N) = (\alpha_H|_{kH \times N})^{-1}(U \cap N) \]
is open in $kH \times N$, for any finitely-generated submodule $N$. Therefore $\alpha_H^{-1}(U)$ is open, by the above. So $\alpha_H$ is continuous.\\
Let $g \in G$. The action of $g$ on $M$ induces a natural bijection $g: N \rightarrow gN$ for any finitely-generated submodule $N$. We show it is a homeomorphism.\\
Now, $gN$ is a finitely-generated $k(gHg^{-1} \cap H)$-submodule of $M$, and $N$ is a finitely-generated $k(H \cap g^{-1}Hg)$-submodule of $M$. The canonical topology on $N, gN$ is determined by these respective module structures. The map $g: N \rightarrow gN$ naturally induces a bijective correspondence between $k(H \cap g^{-1}Hg)$-submodules of $N$ and $k(gHg^{-1} \cap H)$-submodules of $gN$. The open submodules of $N, gN$ are those of finite index, and since the action of $g$ gives a bijection between $N$ and $gN$, it must therefore induce a bijective correspondence between open submodules, $(U \leq N) \mapsto (gU \leq gN)$. Thus, $g:N \rightarrow gN$ is a homeomorphism.\\
Now, let $U \subseteq M$ be open, and $N$ be a finitely-generated $kH$-submodule of $M$. Then $U \cap (kH\cdot gN)$ is open in $kH \cdot gN$ and hence $U \cap gN$ is open in $gN$. Then $g^{-1}(U) \cap N = g^{-1}(U \cap gN)$ is open, by the above. Thus by definition $g^{-1}(U)$ is open. So the action $g:M \rightarrow M$ is continuous for any $g \in G$.\\
Moreover the topology on $M$ is Hausdorff, since it is Hausdorff on any finitely-generated submodule $N$. Thus, the topology on $M$ satisfies the properties of Proposition \ref{augmented rep topology version prop}. Therefore $M$ has a unique $kG$-module structure extending the $k[G]$-module and $kH$-module actions. $\square$\\
\end{jpf}

We now generalise Definition \ref{profinite H-space def} to locally profinite groups and spaces.

\begin{jdef} \h \label{locally profinite def}
A locally profinite $G$-space is a locally compact Hausdorff totally disconnected topological space $X$ with a continuous $G$-action $G \times X \rightarrow X$.\\
\end{jdef}

Notice that for $H$ an open profinite subgroup, the $H$-orbits of a locally profinite $G$-space $X$ are compact, Hausdorff and totally disconnected, and (trivially) each has finitely many $H$-orbits, thus are profinite $H$-spaces. Unlike subsection \ref{Profinite modules for an Iwasawa algebra}, we make no assumption on the number of orbits of $G$ on $X$ here. (In particular we are free to take $G$ and $X$ to be profinite, with infinitely many $G$-orbits.)

\begin{jdef} \h
Let $X$ be a locally profinite $G$-space, $H$ be an open profinite subgroup of $G$. The augmented completed module of $X$ is
\[ kX = \bigoplus_{Z \in \mathrm{Orb}_H(X)} kZ, \]
with the direct sum topology, where each $kZ$ has its natural inverse limit topology.\\
\end{jdef}

We now show that this definition gives a $kG$-module which does not depend on choice of $H$.

\begin{jprop} \h \label{augmented completed module basic properties prop}
Let $X$ be a locally profinite $G$-space. For any $H \leq G$ an open profinite subgroup, there is a $kG$-action extending the left $kH$-module action on $(kX)_H = \bigoplus_{Z \in \mathrm{Orb}_H(X)} kZ$. Moreover, each of these $kG$-modules are naturally isomorphic and the topologies agree.\\
\end{jprop}
\begin{jpf} \h
Let $J \leq H$ be open profinite subgroups of $G$. For $Z$ any $H$-orbit, there is a $kJ$-module isomorphism
\[ i_{J,H}^Z: \bigoplus_{Z' \in \mathrm{Orb}_J(Z)} kZ' \rightarrow kZ, \]
using both parts of Lemma \ref{completed module basic properties lemma}. Thus, define the $kJ$-module isomorphism
\[ j_{J,H} = \bigoplus_{Z \in \mathrm{Orb}_H(X)} i_{J,H}^Z : (kX)_J \rightarrow (kX)_H. \]
We can extend this definition to arbitrary $J, H \leq G$ by $j_{J,H} = j_{J \cap H, H} \circ j_{J \cap H, J}^{-1}$. \vspace{2pt}\\
We define a $kG$-module action on $(kX)_H$, by first defining a $G$-action on $kX$. For $Z$ an $H$-orbit of $X$ and $g \in G$, we have the natural continuous action $g: Z \rightarrow gZ$. We can consider both $Z$ and $gZ$ as profinite $H$-spaces by letting $H$ act on $gZ$ via the isomorphism $H \cong gHg^{-1}$. Thus we obtain a (bijective) map $\tilde{g}: kZ \rightarrow k(gZ)$, by Proposition \ref{completed module universal property prop}. Note that $gZ$ is a profinite $gHg^{-1}$-space, so linearly extending, we obtain a (continuous) map
\[ \tilde{g}_H: (kX)_H \rightarrow (kX)_{gHg^{-1}}. \]
Define the action of $g \in G$ on $m \in (kX)_H$ by $g \cdot m = j_{gHg^{-1}, H}(\tilde{g}_H(m))$. It is easy to check that
\[ \tilde{g}_{H}(j_{J, H}(m))= j_{gJg^{-1},gHg^{-1}}(\tilde{g}_{J}(m)) \]
for any open profinite subgroups $J,H \leq G$ (by first considering the case $J \leq H$). Then for any $g, g' \in G$ and $m \in (kX)_H$, it follows that $(gg') \cdot m = g \cdot (g' \cdot m)$, because $\tilde{g}_{g'H(g')^{-1}}\circ \tilde{g'}_H = \widetilde{gg'}_{H}$ by Proposition \ref{completed module functor prop} and because $j_{H_2, H_3} \circ j_{H_1, H_2} = j_{H_1, H_3}$ for any open profinite subgroups $H_1, H_2, H_3 \leq G$. So, the action of $G$ on $(kX)_H$ is a group action.\\
Moreover, each isomorphism $j_{J,H}$ is $G$-equivariant. The isomorphisms of Lemma \ref{completed module basic properties lemma} are easily seen to be topological isomorphisms, so each $i_{J,H}^Z$, therefore each $j_{J,H}$, is also a topological isomorphism. Hence, every $g \in G$ acts continuously on $(kX)_H$, and $(kX)_H$ is Hausdorff because each $kZ$ is Hausdorff. The $G$-action is obviously $k$-linear, and agrees with the $H$-action induced from the natural $kH$-action. Thus we have a $k[G]$-module action on $(kX)_H$ which agrees with the $kH$-module action on $k[H]$. By Proposition \ref{augmented rep topology version prop}, it follows that there is a unique $kG$-module action on $(kX)_H$ extending these actions.\\
Since each $j_{J,H}$ is a topological $kJ$-module isomorphism, and is $G$-equivariant, it follows that $(kX)_J \cong (kX)_H$ as topological $kG$-modules whenever $J \leq H$, and this obviously extends to arbitrary open profinite subgroups $J, H$. $\square$ \\
\end{jpf}

\begin{jex} \hsmall
Let $G = \mathbb{Q}_p^\times$ and $X = \mathbb{Q}_p$ with $G$ acting by multiplication. A compact open subgroup of $\mathbb{Q}_p^\times$ is $H = \mathbb{Z}_p^\times$. The $H$-orbits in $X$ are $\{p^n \mathbb{Z}_p^\times \mid n \in \mathbb{Z}\} \cup \{ \{0\} \}$, and so
\[ kX = \bigoplus_{n \in \mathbb{Z}} k\mathbb{Z}_p^\times \oplus k \]
as a $k\mathbb{Z}_p^\times$-module, where $k$ is the trivial module. The element $p^m \in \mathbb{Q}_p^\times$ acts on the infinite direct sum by the natural "shift by $m$" among the isomorphic summands, and acts trivially on $k$. Since $k\mathbb{Q}_p^\times \cong k\mathbb{Z}_p^\times[X, X^{-1}]$, we have a $k\mathbb{Q}_p^\times$-module isomorphism
\[k\mathbb{Q}_p \cong k\mathbb{Q}_p^\times\oplus k. \]
\end{jex}

We now prove a universal property analogous to those in the previous three subsections. Similarly to subsection \ref{Universal property of augmented Iwasawa algebras}, this requires us to define a class of modules which generalises ``pseudocompact $kH$-modules'', but which may look unnatural at first sight.

\begin{jdef} \h \label{Iwasawa-like module def}
Let $M$ be a Hausdorff topological left $kG$-module. Let $L \subseteq M$ be a subset closed under the action of $G$. The subset $L$ makes $M$ an Iwasawa-like $kG$-module if and only if there is an open profinite subgroup $H \leq G$ such that every $H$-orbit in $L$ is contained in a pseudocompact $kH$-submodule of $M$.\\
\end{jdef}

For any locally profinite $G$-set $X$ the $kG$-module $kX$ is Iwasawa-like via $X$, because for any open profinite $H \leq G$ and any $H$-orbit $Z$ in $X$, we have that $kZ \leq kX$ is a pseudocompact $kH$-submodule of $kX$.\\
In this case, any open profinite subgroup $H \leq G$ may be considered, and we now show this holds for a general Iwasawa-like module.

\begin{jlem} \h \label{Iwasawa-like via any subgroup lemma}
Let $M$ be an Iwasawa-like $kG$-module via $L \subseteq M$. For any open profinite subgroup $J \leq G$, every $J$-orbit in $L$ is contained in a pseudocompact $kJ$-submodule of $M$.\\
\end{jlem}
\begin{jpf} \h
Let $H \leq G$ be an open profinite subgroup such that every $H$-orbit of $L$ is contained in a pseudocompact $kH$-submodule. Clearly if $J$ is a subgroup of $H$, this property also holds for $J$.\\
Let $J \leq G$ be open profinite with $H \leq J$. Without loss of generality, we may assume $H$ is small enough such that $H$ is normal in $J$. Let $Y$ be a $J$-orbit of $L$ and $Z \subseteq Y$ be an $H$-orbit. By our assumption on $H$, there is a pseudocompact $kH$-submodule $N$ such that $Z \subseteq N$. Then, because $H$ is normal in $J$, we have a $kH$-submodule
\[ P = \sum_{j \in J} jN = \sum_{j \in J/H} jN \leq M \]
which must contain
\[ Z' = \bigcup_{j \in J} jZ \]
Moreover, the $G$-action on $M$ restricts to a $J$-action on $P$ extending the $H$-action. Because $kJ = kH \otimes_{k[H]} k[J]$, it follows that $P$ is a $kJ$-module. Moreover, where $N = \varprojlim_{i \in I} N/N_i$ is an inverse limit of finite-length $kH$-modules, then $P = \varprojlim_{i \in I} \sum_{j \in J/H} j(N/N_i)$ is an inverse limit of finite-length $kJ$-submodules, so $P$ is a pseudocompact $kJ$-submodule.\\
So if $H \leq G$ has the property of Definition \ref{Iwasawa-like module def}, so does any open profinite subgroup contained in or containing $H$. Thus, by considering the inclusions $J \cap H \leq J$ and $J\cap H \leq H$, the statement follows.  $\square$\\
\end{jpf}

That is, for any $kG$-module $M$ and subset $L \subseteq M$ closed under the action of $G$, the condition of Definition \ref{Iwasawa-like module def} holds for a particular open profinite subgroup $H \leq G$ if and only if it holds for all open profinite subgroups.

\begin{jprop} \h \label{augmented completed module universal property prop}
Let $X$ be a locally profinite $G$-space, $M$ be an Iwasawa-like $kG$-module via $L \subseteq M$, and $\phi: X \rightarrow L$ be a continuous $G$-equivariant map. There is a unique continuous $kG$-module homomorphism $\tilde{\phi} : kX \rightarrow M$ extending $\phi$.\\
\end{jprop}
\begin{jpf} \h
By Lemma \ref{Iwasawa-like via any subgroup lemma}, any open profinite subgroup of $G$ verifies the condition of Definition \ref{Iwasawa-like module def}, of $M$ being Iwasawa-like via $L$. Let $H \leq G$ be open profinite, and let $Z$ be an $H$-orbit of $X$.\\
Then $\phi(Z)$ is an $H$-orbit of $L$, so there exists a pseudocompact $kH$-submodule $M_Z \leq M$ containing $\phi(Z)$. By Proposition \ref{completed module universal property prop}, the continuous $H$-equivariant map $\phi|_Z: Z \rightarrow M_Z$ extends to a continuous $kH$-module homomorphism $\psi_Z: kZ \rightarrow M_Z$. Define
\[ \tilde{\phi}_H = \sum_{Z \in \mathrm{Orb}_H(X)} \psi_Z : kX \rightarrow M. \]
Clearly $\tilde{\phi}_H$ is a $kH$-module homomorphism. We show that $\tilde{\phi}_H$ is $G$-equivariant. \\
Recall the definition of the isomorphisms $j_{J,H}$ from the proof of Proposition \ref{augmented completed module basic properties prop}. By Lemma \ref{Iwasawa-like via any subgroup lemma}, Proposition \ref{completed module functor prop}, and the second part of Lemma \ref{completed module basic properties lemma}, $\tilde{\phi}_H \circ j_{J,H} = \tilde{\phi}_J$ for any open profinite subgroups $H, J \leq G$.\\
Now, let $Z \in \mathrm{Orb}_H(X)$, $m \in kZ$, and $g \in G$. Then
\[ \tilde{\phi}_H (g \cdot m) = \tilde{\phi}_{gHg^{-1}} \circ j_{H, gHg^{-1}} (g \cdot m) = \psi_{gZ}(gm). \]
If $m \in k[Z]$ then $\psi_{gZ}(gm) = g \psi_Z(m)$ because $\phi$ is $G$-equivariant and the $\psi$ are $k$-linear. Because the $\psi$ are continuous, $k[Z]$ is dense in $kZ$, and $M$ is Hausdorff, it follows that $\psi_{gZ}(gm) = g \psi_Z(m)$ for all $m \in kZ$. Thus, by linearity, $\tilde{\phi}_H (g \cdot m) = g \tilde{\phi}_H(m)$, for all $m \in kX$.\\
So $\tilde{\phi}_H$ is $G$-equivariant, and is a $kH$-module homomorphism by Lemma \ref{completed module universal property prop}, hence is a $kG$-module homomorphism, which is continuous since each $\psi_Z$ is continuous.\\
Finally, $k[X]$ is dense in $kX$ by similar reasoning to the fourth part of Proposition \ref{properties of kG prop}, using that each $k[Z]$ is dense in $kZ$, and
\[ k[X] = \bigoplus_{Z \in \mathrm{Orb}_H(X)} k[Z]. \]
Now, $M$ is Hausdorff, and clearly $\phi$ has a unique $k$-linear extension to $k[X]$. Therefore $\tilde{\phi} = \tilde{\phi}_H : kX \rightarrow M$ must be the unique continuous $kG$-module homomorphism extending $\phi$. $\square$\\
\end{jpf}

\begin{jprop} \h \label{augmented completed module functor prop}
The mapping $F$ given on objects by $F(X) = kX$ and on morphisms by $F(\phi) = \tilde{\phi}$, is a functor from the category of locally profinite $G$-spaces (with continuous $G$-equivariant maps) to the category of Iwasawa-like topological $kG$-modules (with continuous $kG$-module homomorphisms). The functor preserves surjectivity of morphisms and injectivity of morphisms.\\
\end{jprop}
\begin{jpf} \h
By Proposition \ref{augmented completed module universal property prop}, and identical reasoning to that given in the proof of Proposition \ref{completed module functor prop}, $F$ is a functor.\\
Let $\phi: X \rightarrow Y$ be a continuous $G$-equivariant map between locally profinite $G$-spaces, and let $H \leq G$ be an open profinite subgroup. If $Z$ is an $H$-orbit of $X$, then $\phi(Z)$ is an $H$-orbit of $Y$.\\
If $\phi$ is injective, then each $\tilde{\phi}|_{kZ} = \widetilde{\phi|_Z}$ is injective by Proposition \ref{completed module functor prop}, and so $\tilde{\phi}$ is injective since it is the direct sum of such maps.\\
If $\phi$ is surjective, then for all $H$-orbits $Z' \subseteq Y$, there is an $H$-orbit $Z \subseteq X$ such that $\phi(Z) = Z'$. Again by Proposition \ref{completed module functor prop}, each $\tilde{\phi}|_{kZ}$ is surjective, so has image $kZ' \leq kY$. So $\tilde{\phi}(kX)$ contains the (direct) sum of all such $kZ'$, which is equal to $kY$. Thus $\tilde{\phi}$ is surjective. $\square$\\
\end{jpf}

Notice that a potential ambiguity appears given the two main constructions detailed in this section. Namely, when $G$ is a locally profinite group, $G$ can be considered as a locally profinite $G$-space, or simply as a locally profinite group, and the resulting constructions of the object $kG$ may differ, \textit{a priori}. However, in fact the constructions agree. Let $G, X$ be locally profinite groups and $\phi: G \rightarrow X$ be a continuous group homomorphism. We then have a group action of $G$ on $X$ given by
\[ \phi \times {id}_X : G \times X \rightarrow X, \quad (g,x) \mapsto \phi(g)x. \]
This is easily seen to be a continuous group action, and hence by Proposition \ref{augmented completed module basic properties prop}, we obtain a continuous $kG$-module action on $kX$, which we denote
\[ \widetilde{\phi \times {id}_{X}} : kG \times kX \rightarrow kX. \]
Alternatively, by Proposition \ref{universal property prop}, $\phi$ extends to a continuous $k$-algebra homomorphism
\[ \tilde{\phi} : kG \rightarrow kX, \]
and this gives a continuous $kG$-module action on $kX$ via
\[ \tilde{\phi} \times {id}_{kX}: kG \times kX \rightarrow kX, \quad (y, x) \mapsto \tilde{\phi}(y)x . \]
Now, $\widetilde{\phi \times {id}_{X}}$ and $\tilde{\phi} \times {id}_{kX}$ agree on $k[G] \times k[X]$. Thus, they must agree on $kG \times kX$  because $k[G] \times k[X]$ is dense and $kX$ is Hausdorff.\\
In particular, the construction of $kG$ from the left multiplication action of $G$ on itself is the same as constructing the augmented Iwasawa algebra $kG$ and considering the rank one free left module.

\subsection{Properties of augmented Iwasawa algebras} \label{Properties of augmented Iwasawa algebras}

We now deduce some basic properties of augmented Iwasawa algebras from the universal property of Proposition \ref{universal property prop}.

\begin{jprop} \h
Let $G$ be a locally profinite group. Then $kG$ is isomorphic to its opposite ring $(kG)^{op}$.\\
\end{jprop}
\begin{jpf} \h
Let $A$ be an Iwasawa-like $k$-algebra, via the subgroup $L \leq A^\times$. The opposite algebra $A^{op}$ has the same underlying set as $A$, but reversed multiplication. Then, $A^{op}$ is an Iwasawa-like algebra, via the subgroup $L^{op} \leq (A^{op})^\times$, where $L^{op}$ is the set $L$, with reversed group operation $g \cdot h = hg$.\\
Let $\phi: G \rightarrow L$ be a continuous group homomorphism. Then, define $\phi': G \rightarrow L^{op}$ by ${\phi'(g) = \phi(g)^{-1}}$. Then, $\phi'$ is also a continuous group algebra homomorphism, and so by the universal property of Proposition \ref{universal property prop}, $\phi'$ extends to a continuous $k$-algebra homomorphism $\tilde{\phi'}: kG \rightarrow A^{op}$. Taking opposite algebras obtains a continuous $k$-algebra homomorphism
\[ \tilde{\phi} = \tilde{\phi'}^{op} : (kG)^{op} \rightarrow A. \]
Now, $(kG)^{op}$ is an Iwasawa-like algebra, via the subgroup $G^{op} \leq ((kG)^{op})^\times$. This is because the topology on $(kG)^{op}$ is the same as the topology on $kG$, and because $(kK)^{op}$ is a pseudocompact subalgebra of $(kG)^{op}$ for any open profinite $K \leq G$.\\
The group $G^{op}$ is isomorphic (and homeomorphic) to $G$, via $G \rightarrow G^{op}, g \mapsto g^{-1}$. Thus we have an injection $j: G \rightarrow (kG)^{op}$, and $\tilde{\phi} \circ j = \phi$. So, $\tilde{\phi}$ is a continuous $k$-algebra homomorphism extending $\phi$. It also must be the unique such map, because the values on $k[G^{op}] = k[G]$ are determined by the homomorphism property, and $k[G^{op}]$ is dense in $(kG)^{op}$.\\
So, $j$ embeds the group $G$ into the units of the Iwasawa-like algebra $(kG)^{op}$, such that the universal property of Proposition \ref{universal property prop} is satisfied. Therefore, $(kG)^{op}$ is (canonically) isomorphic to $kG$. $\square$\\
\end{jpf}

As a consequence, $kG$ is left coherent (see Definition \ref{coherent ring def}) if and only if $kG$ is right coherent. Thus we will often refer to augmented Iwasawa algebras as being coherent or not coherent, without a prefix ``left'' or ``right''. It also follows that $kG \cong k[G] \otimes_{k[H]} kH$ as a $(k[G], kH)$-bimodule.\\

To prove statements about quotients of augmented Iwasawa algebras, it is necessary to define the augmentation ideal of a subgroup.

\begin{jdef} \h \label{augmentation ideal def}
Let $G$ be a locally profinite group, $H \leq G$ be a closed subgroup. Consider the locally profinite $G$-space $G/H$ where $G$ acts by left multiplication, the $G$-equivariant quotient map $q_H :G \rightarrow G/H$, and its unique extension $\tilde{q_H}: kG \rightarrow k(G/H)$ to a continuous $kG$-module homomorphism, given by Proposition \ref{augmented completed module universal property prop}. The augmentation ideal of $H$ in $G$ is $\epsilon_G(H) = \Ker \tilde{q_H}$.\\  
\end{jdef}

Notice that $\epsilon_G(H)$ is a closed ideal because $\tilde{q_H}$ is continuous and $k(G/H)$ is Hausdorff. By Proposition \ref{augmented completed module functor prop}, $\tilde{q_H}$ is surjective, giving a natural isomorphism of left $kG$-modules ${kG/\epsilon_{G}(H) \cong k(G/H)}$ for any closed subgroup $H \leq G$. We may occasionally write $\epsilon(H) = \epsilon_G(H)$, when $G$ is clear from context.

\begin{jprop} \h \label{surjection of groups prop}
Let $G$ be a locally profinite group, $N$ a closed normal subgroup. There is a natural surjective ring homomorphism $kG \rightarrow k(G/N)$, with kernel the augmentation ideal $\epsilon_G(N)$.\\
\end{jprop}
\begin{jpf} \h
Note that the $G$-equivariant quotient map $q_N: G \rightarrow G/N$ is a continuous group homomorphism. We can extend this to a continuous $k$-algebra homomorphism $kG \rightarrow k(G/N)$ and to a continuous $kG$-module homomorphism $kG \rightarrow k(G/N)$, by Propositions \ref{universal property prop} and \ref{augmented completed module universal property prop}. By the remark at the end of subsection \ref{Modules for an augmented Iwasawa algebra}, these homomorphisms may be identified, so we have a natural continuous ring homomorphism $kG \rightarrow k(G/N)$ with kernel the (two-sided) ideal $\Ker \tilde{q_N} = \epsilon_G(N)$. $\square$\\
\end{jpf}

We now prove some results on augmentation ideals which will be of later use.

\begin{jlem} \h \label{tensor augmentation lemma}
Let $G' \leq G$ be a closed subgroup. The augmentation ideal $\epsilon_G(G') = kG \otimes_{kG'} \epsilon_{G'}(G')$.\\
\end{jlem}
\begin{jpf} \h
Consider the trivial $kG'$-module $k \cong kG'/\epsilon_{G'}(G')$. We show that $kG \otimes_{kG'} k$ is isomorphic to $k(G/G')$. Suppose $\phi: G/G' \rightarrow M$ is a $G$-equivariant map to an Iwasawa-like $kG$-module as in Proposition \ref{augmented completed module universal property prop}. There is a unique extension of the composition $\phi \circ q_{G'}:G \rightarrow M$ to $\widetilde{\phi \circ q_{G'}}: kG \rightarrow M$. Now, $(\phi \circ q_{G'})|_{G'}$ is constant, and therefore $\widetilde{\phi \circ q_{G'}}|_{kG'}: kG' \rightarrow M$ factors through a $kG'$-module homomorphism $\theta: k \rightarrow M$. Then
\[ f: kG \times k \rightarrow M, \quad f(x,y) = x\theta(y), \]
is a $kG'$-balanced $k$-bilinear map, and $f(x,1) = x\widetilde{\phi \circ q_{G'}}(1) = \widetilde{\phi \circ q_{G'}}(x)$. By the universal property of tensor product, there is a unique extension $ \tilde{f}: kG \otimes_{kG'} k \rightarrow M$, which extends the original map $\phi$, where $gG' \in G/G'$ is identified as $g \otimes 1$. Any other extension of $\phi$ to $kG \otimes_{kG'} k \rightarrow M$ must restrict to $\widetilde{\phi \circ q_{G'}}$ on $kG$ and $\theta$ on $k$, so $\tilde{f}$ is unique. Therefore $kG \otimes_{kG'} k$ satisfies the universal property of Proposition \ref{augmented completed module universal property prop}, and we have a natural isomorphism $k(G/G') \cong kG \otimes_{kG'} k$.\\
By Proposition \ref{locally profinite flatness prop} (which does not depend on this lemma), $kG$ is a flat right $kG'$-module, so there is a commutative diagram with exact rows as follows.
\[ 
\begin{tikzcd}
0 \arrow[r] & \epsilon_G(G') \arrow[r]                               & kG \arrow[r]                    & k(G/G') \arrow[r]                               & 0 \\
0 \arrow[r] & kG \otimes_{kG'} \epsilon_{G'}(G') \arrow[r] \arrow[u] & kG \arrow[r] \arrow[u, "\cong"] & kG \otimes_{kG'} k \arrow[r] \arrow[u, "\cong"] & 0
\end{tikzcd}
\]
Therefore the natural (inclusion) map $kG \otimes_{kG'} \epsilon_{G'}(G') \rightarrow \epsilon_G(G')$ is an isomorphism. $\square$\\
\end{jpf}

Clearly, the augmentation ideal of a closed subgroup $H$ must contain the set $\{h-1 \mid h \in H\}$. In the cases of most interest, this set is a generating set.

\begin{jlem} \h \label{augmentation generation lemma}
Let $G$ be a locally profinite group. Suppose $G' \leq G$ is a closed subgroup such that there is an open profinite $H \leq G'$ with $kH$ Noetherian. Then $\epsilon_G(G') = kG\{g'-1 \mid g' \in G'\}$.\\
\end{jlem}
\begin{jpf} \h
By Lemma \ref{tensor augmentation lemma}, it is enough to prove this in the case $G'=G$. Since $H \leq G$ is open, we have that $G/H$ is discrete and $k(G/H) \cong k[G/H]$. Considering the chain of surjective homomorphisms
\[ kG \rightarrow k(G/H) \rightarrow k = kG/\epsilon_G(G), \]
it follows that $\epsilon_G(G) = \epsilon_G(H) + kG\{g-1 \mid g \in G/H \}$.\\
Now consider the case $G=G'=H$. Since $kH$ is Noetherian, every left ideal of $kH$ is finitely-generated, hence closed by Corollary 22.4 of \jcite{Schn11}. From the definition of the (inverse limit) topology on $kH$, the ideal $kH\{h-1 \mid h \in H\}$ is dense in $\epsilon_H(H)$, but must be closed. So $\epsilon_H(H) = kH\{h-1 \mid h \in H\}$.\\
Thus, by Lemma \ref{tensor augmentation lemma},
\[ \epsilon_G(H) = kG \otimes_{kH} \epsilon_H(H)= kG \otimes_{kH} (kH\{h-1 \mid h \in H\}) = kG\{h-1 \mid h \in H\}. \]
It follows that $\epsilon_G(G) = kG\{h-1 \mid h \in H\} + kG\{g-1 \mid g \in G/H \} \subseteq kG\{g-1 \mid g \in G \}$. The reverse inclusion is obvious, completing the proof. $\square$\\
\end{jpf}

If $\mathcal{O}$ is a complete discrete valuation ring with $p \in \mathcal{O}$ a prime element, and $k=\mathcal{O}$ or $k=\mathcal{O}/p\mathcal{O}$, then Lemma \ref{augmentation generation lemma} holds whenever $G$ is a $p$-adic Lie group. This follows from Theorem 33.4 of \jcite{Schn11}.

\begin{jlem} \h \label{augmentation endomorphism lemma}
Let $G$ be a $p$-adic Lie group, $k=\mathcal{O}$ or $k=\mathcal{O}/p\mathcal{O}$ as above, and $G' \leq G$ be a closed subgroup. Let $\sigma$ be a ring endomorphism of $kG$ which restricts to a group homomorphism $G \rightarrow G$. Then $kG\sigma(\epsilon_G(G')) = \epsilon_G(\sigma(G'))$. \\
\end{jlem}
\begin{jpf} \h
By Lemma \ref{augmentation generation lemma}, any augmentation ideal $\epsilon_G(G')$ is generated abstractly by $\{g-1 \mid g \in G' \}$. The left ideal $kG\sigma(\epsilon_G(G'))$ clearly contains $\{\sigma(g)-1 \mid g \in G' \}$, hence contains $\epsilon_G(\sigma(G'))$. Conversely, if $x = \sum_{g \in G'} c_g (g-1) \in \epsilon_G(G')$, then $\sigma(x) = \sum_{g \in G'} \sigma(c_g) (\sigma(g)-1) \in \epsilon_G(\sigma(G'))$, thus $\epsilon_G(\sigma(G'))$ must contain $kG\sigma(\epsilon_G(G'))$. $\square$\\
\end{jpf}

In particular, Lemmas \ref{augmentation generation lemma} and \ref{augmentation endomorphism lemma} hold whenever $G$ is a $p$-adic Lie group and $k$ is a discrete perfect field of characteristic $p$, because by Theorem II.6.8 of \jcite{Ser79}, the ring of Witt vectors $W(k)$ is a complete discrete valuation ring with residue field $W(k)/pW(k) =k$.

\subsection{Smooth representations of $G$} \label{Smooth representations of G}

To complete this section, we give a proof that any smooth representation is a module for the augmented Iwasawa algebra $kG$, using Proposition \ref{universal property prop}. We give a definition of a smooth representation of a locally profinite group, over any commutative pseudocompact ring, naturally generalising Definitions 2.2.1 and 2.2.5 of \jcite{Eme10}.

\begin{jdef} \h
Let $G$ be a locally profinite group, and $k$ be a commutative pseudocompact ring. A smooth representation of $G$ over $k$ is a $k[G]$-module $V$ such that for any $v \in V$, there exists an open profinite subgroup $H \leq G$ and an open ideal $I \subseteq k$, such that $H \cdot v = \{ v \}$ and $I \cdot v = \{0\}$.\\
\end{jdef}

Showing that any smooth representation carries a $kG$-module structure requires working with finite-length $k$-submodules, necessitating the following lemmas.

\begin{jlem} \h \label{smooth rep k submodule lemma}
Let $V$ be a smooth representation of $G$ over $k$ and $M \leq V$ be a $k$-submodule of $V$. The following are equivalent:
\begin{itemize}
\item $M$ is a finitely-generated $k[H]$-module for some open profinite subgroup $H \leq G$.
\item $M$ is a finite-length $k[H]$-module for some open profinite subgroup $H \leq G$.
\item $M$ is a finitely-generated $k$-module.
\item $M$ is a finite-length $k$-module.\\
\end{itemize}
\end{jlem}

\begin{jlem} \h \label{endomorphism finite length module lemma}
If $M$ is a finite-length $k$-module, then $\mathrm{End}_k(M)$ is a finite length $k$-module.\\
\end{jlem}

With input from these lemmas, we show that the endomorphism ring of any smooth representation is an Iwasawa-like algebra, and deduce that there is a natural and unique $kG$-module structure on any smooth representation of $G$.

\smoothrepsprop*

\begin{jpf} \h
Let $A = \mathrm{End}_k(V)$ be the ring of $k$-linear endomorphisms of $V$. For any $M$ a finite-length $k$-submodule of $V$, define the left ideal of $A$,
\[ I^M = \{f \in A \mid f(M)=0 \}. \]
The collection of such ideals, $\{I^M \mid \text{$M$ a finite-length $k$-submodule of $V$} \}$, is closed under finite intersections. Thus we give $A$ the linear topology with an open neighbourhood basis of 0 given by this collection, and we now show that $A$ is a topological ring. We show that the multiplication of $A$, which is given by composition of functions,
\[ c: A \times A \rightarrow A, \quad c(f,h) = f \circ h ,\]
is continuous. Given $M$ a finite-length $k$-submodule of $V$, consider
\[c^{-1}(I^M) = \{(f,h) \in A \times A \mid f(h(M))=0 \}. \]
Note that for any $h \in A$, if $f \in I^{h(M)}$ and $h' \in I^M$, then $(f,h+h') \in c^{-1}(I^M)$, so
\[I^{h(M)} \times (h+I^M) \subseteq c^{-1}(I^M). \]
Moreover, $c^{-1}(I^M) = \{(f,h) \in A \times A \mid f \in I^{h(M)} \}$. Thus we have a chain of inclusions
\[ c^{-1}(I^M) \subseteq \bigcup_{h \in A} I^{h(M)} \times \{h\} \subseteq \bigcup_{h \in A} I^{h(M)} \times (h+I^M) \subseteq c^{-1}(I^M), \]
in which equality must hold throughout. It follows $c^{-1}(I_M)$ is open in $A \times A$. Thus $c$ is continuous and $A$ is a topological ring. Also, $A$ is Hausdorff since by Lemma \ref{smooth rep k submodule lemma}, $V$ is the union of its finite-length $k$-submodules, so the intersection of all the $I^M$ is zero. \\
Moreover, $A^\times$ is a topological group under the subspace topology. To show this, let $a \in A^\times$, $M$ be a finite-length $k$-submodule of $V$, and consider the open subset $U_{a,M} = (a+I^M) \cap A^\times \subseteq A^\times$. Then
\[ U_{a,M}^{-1} = \{b \in A^\times \mid b(a+f)=1 \text{ for some }f \in I^M \} = \{b \in A^\times \mid (1-ba)(M) =0 \}, \]
with inclusion of the third term in the second given by considering the choice $f=b^{-1}(1-ba)$. It follows that
\[ U_{a,M}^{-1} = \{b \in A^\times \mid (a^{-1}-b)(a(M)) = 0 \} = (a^{-1} + I^{a(M)}) \cap A^\times \]
is open, and so the inversion map is continuous and $A^\times$ is a topological group.\\
Now, for any open profinite $H \leq G$, define the subring of $A$,
\[ B_H = \{ f \in A \mid f(N) \subseteq N \textrm{ for all finite-length $k[H]$-submodules }N \leq V\}. \]
Then $B_H$ is a topological ring, which we show is a pseudocompact $k$-algebra. An open neighbourhood basis of $0 \in B_H$ is given by the collection of ideals of the form $I^M \cap B_H$ for $M$ a finite-length $k$-submodule of $V$. But $I^{k[H]M} \cap B_H \subseteq I^M \cap B_H$, and $k[H]M$ is a finitely-generated $k[H]$-module, hence is a finite-length $k[H]$-module and $k$-module by Lemma \ref{smooth rep k submodule lemma}. Thus an open neighbourhood basis of $0 \in B_H$ is given by the collection of ideals,
\[ \{ I_H^M = I^M \cap B_H \mid \text{$M$ is a finite-length $k[H]$-submodule of $V$}\}. \]
Now, $I_H^M$ is the kernel of the natural $k$-algebra homomorphism from $B_H$ to $\mathrm{End}_{k}(M)$ given by restriction. This gives a natural 
$k$-algebra homomorphism
\[ \psi_H: B_H \rightarrow \varprojlim_{M} \mathrm{End}_{k}(M) = \{(f_M)_M \mid f_M|_N = f_N \text{ if } N \leq M\} \subseteq \prod_{M} \mathrm{End}_{k}(M). \]
Now, $\psi_H$ is injective because $V$ is the union of its finitely-generated, hence by Lemma \ref{smooth rep k submodule lemma}, finite-length $k[H]$-submodules. Given an element $(f_M)_M$ of the inverse limit, we can define an element of $B_H$ by
\[ f: V = \bigcup_{M} M \rightarrow V, \quad f(x) = f_M(x) \text{ if }x \in M,\]
and then $\psi_H(f) = (f_M)_M$. So, $\psi_H$ is an isomorphism.\\
By Lemma \ref{endomorphism finite length module lemma}, $\mathrm{End}_k(M)$ is a finite-length $k$-module for any finite-length $k$-submodule $M \leq V$. Therefore, $B_H$ is a pseudocompact $k$-algebra under the subspace topology induced from $A$.\\
Let $\phi:G \rightarrow A^\times$ be the natural group homomorphism given by $\phi(g)(v) = g \cdot v$.
Because $A^\times$ is a topological group under the subspace topology from $A$, the continuity of $\phi$ can be checked on a basis of open neighbourhoods of $1 \in A^\times$. Let $U_M = U_{1,M} \subseteq A^\times$ for $M$ a finite-length $k$-submodule of $V$. Then
\[ \phi^{-1}(U_M) = \{ g \in G \mid (g-1)(M) = 0 \} = \{g \in G \mid g \cdot m =m, \, \forall m \in M \}. \]
Let $S$ be a finite generating set of $M$ over $k$, and for each $s \in S$, choose an open profinite subgroup $H_s \leq G$ that fixes $s$. Then the open profinite subgroup
\[ H = \bigcap_{s \in S}H_s \]
fixes all elements of $M$, because the actions of $G$ and of $k$ commute. So $H \leq \phi^{-1}(U_M)$, which implies $\phi^{-1}(U_M)$ is an open subgroup of $G$. Thus $\phi$ is continuous.\\
Thus, $\phi: G \rightarrow \phi(G)$ is a surjective continuous map of topological groups, hence is an open map. Therefore $\phi(G)$ is a locally profinite group with open profinite subgroup $\phi(H)$. Moreover, $\phi(H)$ is contained in the units of $B_H$, which is a pseudocompact subalgebra of $A$. This shows that $A$ is an Iwasawa-like algebra, via the subgroup $\phi(G) \leq A^\times$. $\square$\\
\end{jpf}

\smoothrepscor*

\begin{jpf} \h
Recall from the proof of Proposition \ref{endomorphisms of smooth rep is Iwasawa-like prop} the natural group homomorphism
\[ \phi:G \rightarrow \mathrm{End}_k(V)^\times, \]
given by $\phi(g)(v) = g \cdot v$. In that proof we showed that $\mathrm{End}_k(V)$ is Iwasawa-like via $\phi(G)$, and that $\phi$ is continuous. By Proposition \ref{universal property prop}, there is a unique continuous $k$-algebra homomorphism $\tilde{\phi}: kG \rightarrow \mathrm{End}_k(V)$ extending $\phi$. This defines a $kG$-module structure via $x \cdot v = \tilde{\phi}(x)(v)$ for all $x \in kG$, $v \in V$.
Hence, there is a unique $kG$-module structure on the smooth representation $V$ extending the $k[G]$-module structure. $\square$\\
\end{jpf}

We now give the straightforward proofs of Lemmas \ref{smooth rep k submodule lemma} and \ref{endomorphism finite length module lemma}.\\

\begin{jpfof}[Lemma \ref{smooth rep k submodule lemma}] \h
Suppose $M$ is a finitely-generated $k[H]$-module, and let $S$ be a finite set of generators. For each $s \in S$, there exists an open normal subgroup $K_s \trianglelefteq_o H$ that fixes $s$, and an open ideal $I_s \subseteq k$ that annihilates $s$. Thus $M$ is a finitely-generated $(k/I)[H/K]$-module, where
\[ I = \bigcap_{s \in S} I_s, \quad K = \bigcap_{s \in S} K_s. \]
Then, $k/I$ is Artinian and $H/K$ is a finite group, so $(k/I)[H/K]$ is an Artinian ring, thus $M$ has finite length as a $k[H]$-module. Moreover, $(k/I)[H/K]$ is finitely-generated as a $k$-module, so $M$ is finitely-generated as a $k$-module.\\
If $M$ is any finitely-generated $k$-module, there exist open ideals $I_s \subseteq k$ annihilating each of the elements of a finite generating set $S$ for $M$. Then $M$ is a finitely-generated $k/I$-module, and $k/I$ is Artinian (where $I$ is as above), so $M$ is of finite length as a $k$-module.\\
Thus, it remains to show that any finite-length $k$-submodule $M$ is a finitely-generated $k[H]$-submodule for some open profinite subgroup $H \leq G$. Clearly such an $M$ is finitely-generated, let $m_1, \dots, m_n$ be generators for $M$ as a $k$-module. Because $V$ is smooth, there exist open profinite subgroups $H_j$ such that $H_j$ fixes $m_j$, and so the open profinite subgroup
\[H = \bigcap_{j=1}^n H_j\]
fixes each $m_j$, and so $M$ is naturally a finitely-generated $k[H]$-submodule of $V$. $\square$\\
\end{jpfof}

\begin{jpfof}[Lemma \ref{endomorphism finite length module lemma}] \h
We show that for any finite-length $k$-modules $M, N$, that the $k$-module $\jhom_k(M,N)$ is of finite length, inducting on the length of $M$. If $M$ has length 1, then $M$ is simple and $M \cong k/I$ for some maximal ideal $I \trianglelefteq k$. So $\jhom_k(M,N)$ is isomorphic to the submodule $\{x \in N \mid Ix=0 \}$ of $N$, and hence is of finite length. If $M$ is of length greater than 1, let $M' \leq M$ be a submodule such that $M/M'$ is simple. Then $\jhom_k(M/M',N)$ is of finite length, and $M'$ has length strictly less than that of $M$, so $\jhom_k(M',N)$ is of finite length by induction. By considering the exact sequence
\[ 0 \rightarrow \jhom_k(M', N) \rightarrow \jhom_k(M,N) \rightarrow \jhom_k(M/M', N), \]
it follows that $\jhom_k(M,N)$ is a finite-length $k$-module. Putting $M=N$ gives the result. $\square$\\
\end{jpfof}

\section{Coherent rings and skew polynomial rings}

\subsection{Coherent rings} \label{Coherent rings}

We recall the definition of a coherent ring.

\begin{jdef} \h \label{coherent ring def}
Let $R$ be a ring. $R$ is left coherent if and only if every finitely-generated left ideal of $R$ is finitely-presented.\\
\end{jdef}

Note that any left Noetherian ring is left coherent. We can also define the notion of a coherent module.

\begin{jdef} \h \label{coherent module def}
Let $R$ be a ring. A left $R$-module $M$ is coherent if and only if $M$ is finitely-generated, and every finitely-generated submodule of $M$ is finitely-presented.\\
\end{jdef}

A ring $R$ is left coherent if and only if it is coherent as a left $R$-module. A core motivation for considering coherent rings in the context of representation theory comes from the following proposition, see \jcite{Glaz89}, Theorems 2.5.1 and 2.1.2. (Technically only commutative rings are considered, but the proof works identically for non-commutative rings.)

\begin{jprop} \h \label{coherent abelian category prop}
Let $R$ be a ring. Then, $R$ is left coherent if and only if the category of finitely-presented left $R$-modules is an abelian subcategory of the category of left $R$-modules.\\
\end{jprop}

Shotton makes use of Proposition \ref{coherent abelian category prop} to deduce that the finitely-presented smooth mod $p$ representations of $SL_2(F)$ form an abelian category, see Theorems 1.2 and 4.5 of \jcite{Sho20}.

\subsection{Skew polynomial rings}

See Chapter 2 of \jcite{GoWa04} for an extensive treatment of the theory of skew polynomial rings.

\begin{jdef} \h
Let $R$ be a ring, and $\sigma_X: R \rightarrow R$ be a ring endomorphism of $R$. The skew polynomial ring $R[X; \sigma_X]$ is the free left $R$-module $\bigoplus_{j \geq 0} RX^j$, given a ring structure by addition being $R$-module addition, and multiplication given by
\[ r_1X^{n_1} \cdot r_2 X^{n_2} = r_1 \sigma_X^{n_1}(r_2) X^{n_1+n_2}, \] 
extended suitably.\\
\end{jdef}

In this article we consider only the skew polynomial rings of the above type. Discussion of the more general theory can be found be in \jcite{GoWa04}, Chapter 2, page 32. We will often omit the endomorphism from our notation and write simply $R[X; \sigma_X] = R[X]$.\\

\section{Faithful flatness for augmented Iwasawa algebras} \label{Faithful flatness for augmented Iwasawa algebras}

In this section we show that the mod $p$ augmented Iwasawa algebra of any $p$-adic Lie group is faithfully flat over the augmented Iwasawa algebra of any closed subgroup. 

\subsection{The compact case} \label{The compact case}

Let $k$ be a commutative pseudocompact ring throughout this section.

\begin{jthm} \label{profinite flatness thm} \h 
Let $G$ be a profinite group, and $H \leq G$ be a closed subgroup. Then $kG$ is a flat right $kH$-module.\\
\end{jthm}
\begin{jpf} \h
In Lemma 4.5 of \jcite{Bru66}, Brumer proved that $kG$ is a projective object in the category of pseudocompact $kH$-modules. This does not imply that $kG$ is a projective $kH$-module in the usual sense, so we cannot immediately deduce flatness. However, Lemma 2.1 of \jcite{Bru66} implies that the completed tensor product $kG \widehat\otimes_{kH}$ is an exact functor and $kG \widehat\otimes_{kH}M = kG \otimes_{kH} M$ for any finitely-generated left $kH$-module. It follows that if $I' \subseteq kH$ is a finitely-generated left ideal, applying $kG\otimes_{kH}$ to the short exact sequence
\[ 0 \rightarrow I' \rightarrow kH \rightarrow kH/I' \rightarrow 0 \]
gives a short exact sequence
\[ 0 \rightarrow kG \otimes_{kH} I' \rightarrow kG \rightarrow kG \otimes_{kH}(kH/I') \rightarrow 0. \]
Let $I \subseteq kH$ be an arbitrary ideal of $kH$, and let $I = \varinjlim_{\alpha \in A} I_\alpha$ be the direct limit of the finitely-generated ideals it contains. By the above, there is a short exact sequence
\[ 0 \rightarrow kG \otimes_{kH} I_\alpha \rightarrow kG \rightarrow kG \otimes_{kH}(kH/I_\alpha) \rightarrow 0, \]
for each $\alpha \in A$. By Theorem 2.6.15 of \jcite{Wei94}, direct limits of $kH$-modules are exact, so we have a short exact sequence
\[ 0 \rightarrow \varinjlim_{\alpha \in A}( kG \otimes_{kH} I_\alpha) \rightarrow kG \rightarrow \varinjlim_{\alpha \in A}\left( kG \otimes_{kH}(kH/I_\alpha)\right) \rightarrow 0. \]
Since tensor products commute with direct limits, we have that
\[ 0 \rightarrow kG \otimes_{kH} I \rightarrow kG \rightarrow kG \otimes_{kH}(kH/I) \rightarrow 0\]
is exact. Therefore $\Tor_1^{kH}(kG, kH/I)=0$ for any left ideal $I \subseteq kH$, so by Proposition 3.2.4 of \jcite{Wei94}, $kG$ is a flat right $kH$-module. $\square$\\
\end{jpf}

In this section we will use a version of Mackey's restriction formula to generalise this result and deduce various corollaries.

\subsection{Mackey's restriction formula}

If $G$ is a finite group, Mackey's restriction formula tells us how to compute restrictions of a representation induced from a subgroup of $G$.

\begin{jthm}[Mackey's restriction formula] \h
Let $K$ be any field, $G$ be a finite group, $H_1, H_2 \leq G$ be subgroups. Let $W$ be a left $K[H_2]$-module. Then
\[ \Res_{H_1}^{G} \Ind_{H_2}^G W \cong \bigoplus_{g \in H_1 \backslash G / H_2} \Ind_{g(H_2)g^{-1} \cap H_1}^{H_1} W_g, \]
as left $K[H_1]$-modules.\\
\end{jthm}
Here $W_g$ is the $K$-vector space of symbols $g \otimes W$, with action $ghg^{-1} \cdot (g \otimes w) = g \otimes hw$ for all $h \in H_2$. $\Ind_A^B$ means the functor $K[B] \otimes_{K[A]}$. We interpret the sum over $g \in H_1 \backslash G / H_2$ to mean, ``choose a set $X \subseteq G$ such that $X$ is a set of representatives for the double cosets $H_1 \backslash G / H_2$, and sum over this set". \\ See \jcite{Ser77}, Proposition 22, for a proof of this result.\\

In the next subsection, we show that an analogous statement to Mackey's restriction formula for finite groups holds for locally profinite groups. To be precise, we will prove the following theorem.\\

\begin{jthm} \h \label{Mackey augmented Iwasawa thm}
Let $G$ be a locally profinite group. Let $G_1 \leq G$ be a closed subgroup, and $H \leq G$ be an open profinite subgroup. Let $M$ be a left $kG_1$-module. Then
\[ \Res_{H}^{G} \Ind_{G_1}^{G} M \cong \bigoplus_{g \in H \backslash G / G_1} \Ind_{gG_1g^{-1} \cap H}^{H} M_g, \]
as left $kH$-modules. \\
\end{jthm}

Here $\Ind_{A}^{B}$ means the functor $kB \otimes_{kA}$. We interpret $M_g$ similarly to the finite group case, as the $k$-vector space of symbols $g \otimes M$, with action $gxg^{-1} \cdot(g \otimes m) = g \otimes xm$ for all $x \in kG_1$, noticing that $k(gG_1g^{-1})$ is equal to $g(kG_1)g^{-1}$ as subrings of $kG$.\\

\subsection{Proof of the restriction formula}

In this subsection let us fix $G$ a locally profinite group, $G_1 \leq G$ a closed subgroup, and $H \leq G$ an open profinite subgroup.

Consider the functors which make up each side of the Mackey formula:
\begin{jdef} \h
The functors $\textrm{Mack}_1, \textrm{Mack}_2 : kG_1\textrm{-Mod} \rightarrow kH\textrm{-Mod}$ are
\[ \textrm{Mack}_1 = \Res_{H}^{G} \Ind_{G_1}^{G}, \quad \textrm{Mack}_2 = \bigoplus_{g \in H \backslash G / G_1} \Ind_{gG_1g^{-1} \cap H}^{H} (\textrm{  })_g. \]
\end{jdef}
\begin{jdef} \h
For any $g \in G$, $\textrm{Mack}_{2,g}: kG_1\textrm{-Mod} \rightarrow kH\textrm{-Mod}$ is the functor $\Ind_{gG_1g^{-1} \cap H}^{H} (\textrm{  })_g$ .\\
\end{jdef}

The statement of Theorem \ref{Mackey augmented Iwasawa thm} is then that $\textrm{Mack}_1(M) \cong \textrm{Mack}_2(M)$ for any $M$. We will in fact show that there is a canonical such isomorphism, defined next. Recall from Proposition \ref{augmented Iwasawa algebra definition prop} that the augmented Iwasawa algebra $kG \cong kH \otimes_{k[H]} k[G]$ as a left $kH$-module.

\begin{jdef} \h
Let $M$ be a left $kG_1$-module. We define the $kH$-module homomorphism
\[ \psi_M:  \textrm{Mack}_2(M) \rightarrow \textrm{Mack}_1(M) \] to be \[ \psi_M = \bigoplus_{g \in H\backslash G / G_1} \psi_{M,g}\] where
\[ \psi_{M,g} : \Ind_{gG_1g^{-1} \cap H}^{H} M_g \rightarrow \Res_{H}^{G} \Ind_{G_1}^{G} M \]
is given by \[\psi_{M,g}(r \otimes_{k(gG_1g^{-1} \cap H)} (g \otimes m)) = (r \otimes_{k[H]} g) \otimes_{kG_1} m, \] for all $r \in kH, m \in M$. \\
\end{jdef}

\begin{jprop} \h \label{Mackey kG_1 case prop}
The $kH$-module homomorphism $\psi_{kG_1}: \textrm{Mack}_2(kG_1) \rightarrow \textrm{Mack}_1(kG_1)$ is an isomorphism.\\
\end{jprop}
\begin{jpf} \h
Throughout, let $M = kG_1$, a left $kG_1$-module under left multiplication. Fix a set $X \subseteq G$ of representatives for the double cosets $H \backslash G / G_1$.\\
Let $g \in X$. The left $k(gG_1g^{-1})$-module $M_g$ is given by the symbols $g \otimes M = g \otimes kG_1$, and so as a left module may be considered as $g \otimes kG_1 = k(gG_1g^{-1}) \otimes g$.\\ Then, because $gG_1g^{-1} \cap H$ is an open profinite subgroup of $gG_1g^{-1}$,
\[ k(gG_1g^{-1}) = k(gG_1g^{-1} \cap H) \otimes_{k[gG_1g^{-1} \cap H]} k[gG_1g^{-1}], \]
as a left $k(gG_1g^{-1} \cap H)$-module. Therefore 
\begin{align*} \Ind_{gG_1g^{-1} \cap H}^{H} k(gG_1g^{-1})  &= kH \otimes_{k(gG_1g^{-1} \cap H)} k(gG_1g^{-1} \cap H) \otimes_{k[gG_1g^{-1} \cap H]} k[gG_1g^{-1}]\\ &= kH  \otimes_{k[gG_1g^{-1} \cap H]} k[gG_1g^{-1}]\\ &= \bigoplus_{ a \in (gG_1g^{-1} \cap H) \backslash gG_1g^{-1} } kH \otimes a . \end{align*}
Hence, as $M_g = g \otimes kG_1 = k(gG_1g^{-1}) \otimes g$,
\[ \textrm{Mack}_{2,g}(M) = \Ind_{gG_1g^{-1} \cap H}^{H} M_g = \bigoplus_{ a \in (gG_1g^{-1} \cap H) \backslash gG_1g^{-1} } kH \otimes ag .\]
Note that each $kH \otimes ag$ is a free rank 1 left $kH$-module. \vspace{10pt}\\
For each $g \in X$, choose a set of representatives $Y_g \subseteq gG_1g^{-1}$ for the right cosets\\ $(gG_1g^{-1} \cap H) \backslash gG_1g^{-1}$.
It then follows that 
\[ \textrm{Mack}_2(M) = \bigoplus_{g \in H \backslash G / G_1} \Ind_{gG_1g^{-1} \cap H}^{H} M_g   =  \bigoplus_{(g,a) \in S} kH \otimes ag,\]
where $S$ is the set
\[ S = \{ (g,a) \in G \times G \mid g \in X, a \in Y_g \} .\vspace{15pt}\]
Then, $\textrm{Mack}_1(M)$ can be straightforwardly seen to be
\[ \textrm{Mack}_1(M) = \Res_{H}^{G} \Ind_{G_1}^{G} kG_1 = \Res_{H}^{G} kG = \Res_{H}^{G}( kH \otimes_{k[H]}k[G]) = \bigoplus_{b \in H \backslash G} kH \otimes b .\]
We thus have the homomorphisms (for each $g \in X$),
\[ \psi_{M,g}: \bigoplus_{ a \in Y_g } kH \otimes ag \rightarrow \bigoplus_{b \in H \backslash G} kH \otimes b, \]
and their sum
\[ \psi_M :  \bigoplus_{(g,a) \in S} kH \otimes ag \rightarrow \bigoplus_{b \in H \backslash G} kH \otimes b. \]
Let $g \in X$, and $a \in Y_g$. Let $r \in kH$. Using the identifications described above, we have
\[ \psi_{M,g}(r \otimes ag) = \psi_{M,g} (r \otimes_{k(gG_1g^{-1} \cap H)} (g \otimes g^{-1}ag)) = (r \otimes_{k[H]}g) \otimes_{kG_1} g^{-1}ag. \]
Then $r \otimes_{k[H]}g \in kG$ and $g^{-1}ag \in kG_1$, so their tensor product lies in $kG$, and we have
\[ (r \otimes_{k[H]}g) \otimes_{kG_1} g^{-1}ag = r \otimes_{k[H]}gg^{-1}ag = r \otimes_{k[H]}ag \in kH \otimes ag  .\]
So $\psi_M$ is the identity when restricted to $kH \otimes ag \rightarrow kH \otimes ag \leq \textrm{Mack}_1(M)$. \vspace{2pt}\\
Thus, to show that $\psi_M$ is an isomorphism, it is enough to check that $Z = \{ag \mid (a,g) \in S \} \subseteq G$ is a set of right coset representatives for $H$ in $G$, that is, $Z$ is a set of representatives for $H\backslash G$.\\
Because $Y_g$ is a set of representatives for $(gG_1g^{-1} \cap H) \backslash gG_1g^{-1}$, it is easy to see that $g^{-1}Y_g g$ is a set of representatives for $(G_1 \cap g^{-1}Hg) \backslash G_1$. By Exercise I.8 of \jcite{Lang02} (interchanging left and right), it follows that
\[ \{ gy \mid g \in X, b \in g^{-1}Y_gg \} = \{ ag \mid g \in X, a \in Y_g \} = Z \]
is a set of representatives for $H \backslash G$. Therefore $\psi_M$ is an isomorphism. $\square$\\
\end{jpf}

Next we extend the result of Proposition \ref{Mackey kG_1 case prop} to arbitrary free modules.

\begin{jlem} \h
Let $M$ be a left $kG_1$-module, and $M' = \bigoplus_{J} M$, for some index set $J$. If $\psi_M$ is an isomorphism, then $\psi_{M'}$ is an isomorphism.\\
\end{jlem}
\begin{jpf} \h
Because $\textrm{Mack}_1, \textrm{Mack}_2$ are functors composed of direct sums of tensor products, and because tensor products commute with arbitrary direct sums, we have that \begin{align*} \textrm{Mack}_1(M') &= \textrm{Mack}_1\left( \bigoplus_J M\right) = \bigoplus_J \textrm{Mack}_1(M), \\ \textrm{Mack}_2(M') &= \textrm{Mack}_2\left( \bigoplus_J M\right) = \bigoplus_J \textrm{Mack}_2(M), \end{align*}
and
\[ \psi_{M'} = \psi_{\oplus_ J M} = \bigoplus_J \psi_M: \textrm{Mack}_2(M') \rightarrow \textrm{Mack}_1(M'). \]
It follows that if $\psi_M$ is an isomorphism, so is $\psi_{M'}$. $\square$\\
\end{jpf}

\begin{jcor} \h \label{Mackey free cor}
If $F$ is a free $kG_1$-module, then $\psi_F$ is an isomorphism.\\
\end{jcor}

We can now extend to arbitrary modules, providing a proof of Theorem \ref{Mackey augmented Iwasawa thm}.\\

\begin{jpfof}[Theorem \ref{Mackey augmented Iwasawa thm}] \h
Let $M$ be a left $kG_1$-module. Then $M$ has a free presentation
\[ F_1 \rightarrow F_0 \rightarrow M \rightarrow 0. \]
Now, each of the functors $\textrm{Mack}_1, \textrm{Mack}_2$ are (direct sums of) tensor products, and hence are right exact functors. Therefore we have exact sequences 
\[ \textrm{Mack}_1(F_1) \rightarrow \textrm{Mack}_1(F_0) \rightarrow \textrm{Mack}_1(M) \rightarrow 0 \]
and
\[ \textrm{Mack}_2(F_1) \rightarrow \textrm{Mack}_2(F_0) \rightarrow \textrm{Mack}_2(M) \rightarrow 0 .\]
Moreover, the homomorphisms $\psi_M, \psi_{F_0}, \psi_{F_1}$ give us the following commutative diagram, which has exact rows.
\[ \begin{tikzcd}
\textrm{Mack}_1(F_1) \arrow[r]                         & \textrm{Mack}_1(F_0) \arrow[r]                         & \textrm{Mack}_1(M) \arrow[r]                     & 0 \arrow[r]                     & 0                     \\
\textrm{Mack}_2(F_1) \arrow[r] \arrow[u, "\psi_{F_1}"] & \textrm{Mack}_2(F_0) \arrow[r] \arrow[u, "\psi_{F_0}"] & \textrm{Mack}_2(M) \arrow[r] \arrow[u, "\psi_M"] & 0 \arrow[r] \arrow[u, "\psi_0"] & 0 \arrow[u, "\psi_0"]
\end{tikzcd} \]
By Corollary \ref{Mackey free cor}, $\psi_{F_1}, \psi_{F_0}$ are isomorphisms and $\psi_0 $ is clearly an isomorphism. Thus by the five lemma, $\psi_M$ must be an isomorphism. So, $\textrm{Mack}_1(M) \cong \textrm{Mack}_2(M)$. $\square$\\
\end{jpfof}

\subsection{Faithful flatness of augmented Iwasawa algebras} \label{Faithful flatness of augmented Iwasawa algebras}

The formula in Theorem \ref{Mackey augmented Iwasawa thm} now allows us to prove results similar to Theorem \ref{profinite flatness thm}, for locally profinite groups.

\begin{jprop} \h \label{locally profinite flatness prop}
Let $G$ be a locally profinite group, $G_1 \leq G$ be a closed subgroup. Then $kG$ is a flat right $kG_1$-module.\\
\end{jprop}
\begin{jpf} \h
To show that $kG$ is a flat $kG_1$-module, we show that for any injection $i : M \rightarrow N$ of $kG_1$-modules, the induced map
\[j : kG \otimes_{kG_1} M \rightarrow kG \otimes_{kG_1} N\]
is an injection. Let $H$ be an open profinite subgroup of $G$.\\
Note that $j = \textrm{Mack}_1(i)$. For each $g \in H \backslash G / G_1$, let $j_g = \textrm{Mack}_{2,g}(i)$, so that
\[ j_g : \Ind_{gG_1g^{-1} \cap H}^{H} M_g \rightarrow \Ind_{gG_1g^{-1} \cap H}^{H} N_g .\]
Since $gG_1g^{-1}$ is a closed subgroup of $G$, we have that $gG_1g^{-1} \cap H$ is a closed subgroup of $H$, and thus by Theorem \ref{profinite flatness thm}, $\Ind_{gG_1g^{-1} \cap H}^{H}$ is a (left) exact functor, as is $(\textrm{ })_g$. It follows that $j_g$ is injective.\\
By the proof of Theorem \ref{Mackey augmented Iwasawa thm}, we have that
\[ j = \textrm{Mack}_1(i) = \psi_{N} \circ \textrm{Mack}_2(i) \circ \psi_M^{-1}. \]
Then
\[ \textrm{Mack}_2(i) = \bigoplus_{g \in H \backslash G / G_1} \textrm{Mack}_{2,g}(i) = \bigoplus_{g \in H \backslash G / G_1} j_g \]
is injective because each $j_g$ is injective. Because $\psi_M, \psi_N$ are isomorphisms, $j$ must be injective. $\square$\\
\end{jpf}

For the remainder of this subsection, let $\mathcal{O}$ be a complete discrete valuation ring with residue field $k_{\mathcal{O}}$ of characteristic $p$, and assume $k=\mathcal{O}$ or $k = k_{\mathcal{O}}$. For instance, $k$ may be a perfect field of characteristic $p$, with the discrete topology.

\begin{jthm} \h \label{non-compact locally pro-p faithfully flat thm}
Let $G$ be a locally pro-$p$ group, $G_1 \leq G$ be a closed subgroup. Then $kG$ is a faithfully flat right $kG_1$-module.\\
\end{jthm}

Given the proposition above, clearly only the faithful part of the theorem requires proving, and to do this we first show it for pro-$p$ groups.

\begin{jlem} \h \label{pro-p faithful lemma}
Let $U$ be a pro-$p$ group, $H \leq U$ be a closed subgroup. Then $kU$ is a faithfully flat right $kH$-module.\\
\end{jlem}
\begin{jpf} \h
We show that $kU \otimes_{kH} M \neq 0$ for any non-zero left $kH$-module. Because $kU$ is a flat right $kH$-module by Theorem \ref{profinite flatness thm}, it is enough to show this for finitely-generated modules $M$. Now, there is an obvious projection map from $kU\otimes_{kH} M$ onto
\[ \faktor{kU\otimes_{kH} M}{ \epsilon_U(U)(kU\otimes_{kH} M) } \cong \faktor{kU}{\epsilon_U(U)} \otimes_{kH} M \cong k \otimes_{kH} M \cong \faktor{M}{\epsilon_H(H)M}. \]
But $H$ is pro-$p$, since it is a closed subgroup of a pro-$p$ group, so $\epsilon_H(H)$ is the Jacobson radical of $kH$, by Proposition 19.7 of \jcite{Schn11}. By Nakayama's lemma, $M/\epsilon_H(H)M$ is non-zero, hence $kU \otimes_{kH} M \neq 0$. $\square$\\
\end{jpf}

We can now prove the result for all locally pro-$p$ groups.\\

\begin{jpfof}[Theorem \ref{non-compact locally pro-p faithfully flat thm}] \h
Let $M$ be a non-zero left $kG_1$-module. Let $H \leq G$ be an open pro-$p$ subgroup. Then by Theorem \ref{Mackey augmented Iwasawa thm},
\[ \Res_{H}^{G} \Ind_{G_1}^{G} M \cong \bigoplus_{g \in H \backslash G / G_1} \Ind_{gG_1g^{-1} \cap H}^{H} M_g .\]
But taking $g=e \in G$, we have that $\Ind_{H \cap G_1}^{H} M$ is identified with a submodule of $\Ind_{G_1}^{G} M$, and is non-zero by Lemma \ref{pro-p faithful lemma}. Thus $kG \otimes_{kG_1} M \neq 0$. So, $kG$ is a faithfully flat right $kG_1$-module, by Proposition \ref{locally profinite flatness prop}. $\square$\\
\end{jpfof}

Any $p$-adic Lie group is locally pro-$p$, so the following result is immediate.

\begin{jthm} \h \label{non-compact faithfully flat thm}
Let $G$ be a $p$-adic Lie group, $G_1 \leq G$ be a closed subgroup. Then $kG$ is a faithfully flat right $kG_1$-module.\\
\end{jthm}

\section{Coherent rings}

\subsection{Extensions of coherent rings}

It will be important to us to be able to transfer the property of coherence when extending a ring in some fashion. The following was proved by Harris in Corollary 1.2 of \jcite{Har66}.

\begin{jlem} \h \label{fp extension coherence lemma}
Let $R$ be a left coherent ring, and $\phi: R \rightarrow S$ be a ring homomorphism such that $S$ is finitely-presented as a left $R$-module. Then $S$ is a left coherent ring.\\
\end{jlem}

One consequence of this lemma is an analogue, for coherent rings, of the fact that any quotient of a Noetherian ring is Noetherian. See Theorem 2 of \jcite{Har67}.

\begin{jlem} \h \label{finite ideal quotient coherence lemma}
Let $R$ be a left coherent ring, $I \subseteq R$ be a two-sided ideal. Suppose that $I$ is finitely-generated as a left ideal. Then the ring $R/I$ is left coherent.\\
\end{jlem}

The next result can be found as Theorem 2.3.3 of \jcite{Glaz89} -- the proof is identical for commutative or non-commutative rings.

\begin{jlem}\label{direct limit coherence lemma} \h
Let $R$ be a ring such that $R = \varinjlim_{a \in A} R_a$ is a directed colimit of some left coherent subrings $R_a$, ordered by inclusion, and such that $R$ is a flat right $R_a$-module, $\forall a \in A$. Then R is left coherent. \\
\end{jlem}

\begin{jex} \hsmall
Let $G = \mathbb{Q}_p$. Then $G$ is the direct limit of its compact open subgroups $H_n = p^{-n}\mathbb{Z}_p$, and so the augmented Iwasawa algebra of $G$ is a direct limit,
\[ kG = \varinjlim_{n \geq 0} kH_n. \]
If $k$ is a (discrete) perfect field of characteristic $p$, then $kH_n$ is a left Noetherian ring, hence left coherent. By Proposition \ref{locally profinite flatness prop}, $kG$ is a flat right $kH_n$-module, since $H_n \leq G$ is a closed subgroup. So by Lemma \ref{direct limit coherence lemma}, $kG$ is a left coherent ring. \\
In this case $kG \cong \varinjlim_{n \geq 0} k[[X^{\frac{1}{p^n}}]]$, which we write as $kG \cong k[[X]]^{\frac{1}{p^\infty}}$.\\
\end{jex}
\begin{jex} \hsmall
More generally, if $G = \mathbb{Q}_p^m$, then
\[ kG = \varinjlim_{n \geq 0} k(p^{-n}\mathbb{Z}_p^m) \cong \varinjlim_{n \geq 0} k[[X_1^{\frac{1}{p^n}}, \dots, X_m^{\frac{1}{p^n}}]], \]
which we write as $kG \cong k[[X_1, \dots, X_m]]^{\frac{1}{p^\infty}}$, and this ring is also coherent. \\
\end{jex} 

The following lemma was first proved, for commutative rings only, by Harris in \jcite{Har66}. For the convenience of the reader, we have replicated this proof, with the appropriate qualifications needed in the non-commutative case.

\begin{jlem} \label{localisations are coherent modules lemma} \h
Let $R$ be a left coherent ring. Let $X \subseteq R$ be a left denominator set (a left reversible Ore set). Let $M$ be a left coherent $R$-module. Then $X^{-1}M$ is left coherent as a $X^{-1}R$-module.\\
\end{jlem}
\begin{jpf} \h
By Proposition 10.12 and Corollary 10.13 of \jcite{GoWa04}, $X^{-1}R\otimes_R A \cong X^{-1}A$ for any $R$-module $A$, and $X^{-1}R$ is flat as a (right) $R$-module. \\
Let $N' \leq X^{-1}M$ be a finitely-generated $X^{-1}R$-module. Since $X^{-1}M \cong X^{-1}R \otimes_R M$, and by Lemma 10.2 of \jcite{GoWa04}, the generators of $N'$ can be written with a common (left) denominator, hence there exist $m_1, m_2, \dots m_r \in M$ such that
\[N' =  X^{-1}R \cdot1^{-1}m_1 +  X^{-1}R \cdot 1^{-1}m_2 + \dots + X^{-1}R \cdot 1^{-1}m _r  .\]
Let 
\[ N =  R m_1 + R m_2 + \dots + Rm _r  \leq M. \] 
Let $i: N \rightarrow M$ be the natural $R$-module injection. Tensoring, we have the induced map \[j:X^{-1}R \otimes_R N \rightarrow X^{-1}R \otimes_R M, \]and because $X^{-1}R$ is flat as an $R$-module, $j$ is injective.\\
So we have an injection $j: X^{-1}N \rightarrow X^{-1}M $, and clearly $\textrm{Im }j = N'$. So $X^{-1}N \cong N'$ as\\$X^{-1}R$-modules. \vspace{4pt}\\
Now, since $M$ is coherent and $N$ is finitely-generated, $N$ is finitely-presented. Let
\[ 0 \rightarrow K \rightarrow R^n \rightarrow N \rightarrow 0\] be an exact sequence of $R$-modules with $K$ finitely-generated. Again, because $X^{-1}R$ is flat,
\[ 0 \rightarrow X^{-1}R \otimes_R K \rightarrow {(X^{-1}R)}^n \rightarrow X^{-1}N \rightarrow 0 \]
is exact. Then $X^{-1}R \otimes_R K$ is finitely-generated, so $N' \cong X^{-1}N$ is finitely-presented. \vspace{2pt}\\
So $X^{-1}M$ is a coherent $X^{-1}R$-module. $\square$\\
\end{jpf}

\begin{jcor} \h \label{localisations are coherent rings cor}
Let $R$ be a left coherent ring, and $X \subseteq R$ be a left denominator set. If $R$ is left coherent, then $X^{-1}R$ is left coherent.\\
\end{jcor}

See Theorem 2 and Corollary 2.1 of \jcite{Har66} for a proof of the following result.

\begin{jprop} \h \label{coherence faithfully flat descent prop}
Let $\phi : R \rightarrow S$ be a ring homomorphism such that $S$ is a faithfully flat right $R$-module, via $\phi$. If $S$ is left coherent, then $R$ is left coherent.\\
\end{jprop}

We can combine this proposition with Theorem \ref{non-compact faithfully flat thm} to deduce the following. For the rest of this section, let $\mathcal{O}$ be a complete discrete valuation ring with $p \in \mathcal{O}$ a prime element, and assume $k = \mathcal{O}$ or $k=\mathcal{O}/p\mathcal{O}$.

\begin{jprop} \h \label{closed subgroup coherence prop}
Let $G$ be a $p$-adic Lie group, $H$ be a closed subgroup of $G$. If $kG$ is a coherent ring, then $kH$ is a coherent ring.\\
\end{jprop}

Notice that the contrapositive to this proposition says that $kG$ cannot be coherent if even one of its closed subgroups has a non-coherent augmented Iwasawa algebra.\\

\subsection{Coherence for unipotent $p$-adic Lie groups}

We can use the results above to show that unipotent $p$-adic Lie groups always have a coherent augmented Iwasawa algebra.  

\begin{jcor} \h \label{unipotent coherent cor}
Let $F$ be a finite field extension of $\mathbb{Q}_p$, $\mathbb{U}_n$ be the affine group scheme of upper unitriangular matrices in $\mathbb{GL}_n$. Let $U$ be a closed subgroup of $\mathbb{U}_n(F)$ with the $p$-adic topology. Then the augmented Iwasawa algebra $kU$ is a coherent ring.\\
\end{jcor}
\begin{jpf} \h
First we show that $\mathbb{U}_n(F)$ is a direct limit of compact subgroups.\\
Let $\pi \in F$ be a uniformiser, and $v_F$ be the discrete valuation on $F$. For each $j \in \mathbb{Z}_{\geq 0}$, let
\[ G_j = \begin{pmatrix} 1 & \pi^{-j}\mathcal{O}_F & \pi^{-2j}\mathcal{O}_F & \hdots & \pi^{-(n-1)j}\mathcal{O}_F\\ & 1 & \pi^{-j}\mathcal{O}_F & \ddots  & \vdots\\ && \ddots & \ddots & \pi^{-2j}\mathcal{O}_F\\ &&& 1 & \pi^{-j}\mathcal{O}_F \\ &&&&1 \end{pmatrix} \subseteq \mathbb{U}_n(F). \]
Then it is easy to check that $G_j$ is a subgroup of $\mathbb{U}_n(F)$, and clearly it is compact and closed. We also have that $G_{j} \leq G_{j+1}$ for all $j$. If $A \in \mathbb{U}_n(F)$, let $N = -\min\{v_F(A_{ij}) \mid 1 \leq i \leq j \leq n \}$. Then $A \in G_N$. So, $\mathbb{U}_n(F)$ is the direct limit of the $G_j$.\\
Then $k\mathbb{U}_n(F) = \varinjlim_{j \geq 0} kG_j$. Because the $G_j$ are compact, each $kG_j$ is a left Noetherian (from Theorem 33.4 of \jcite{Schn11}), hence left coherent, subring of $k\mathbb{U}_n(F)$. Because the $G_j$ are closed, $k\mathbb{U}_n(F)$ is a flat right $kG_j$-module for all $j \geq 0$, by Proposition \ref{locally profinite flatness prop}. Thus $k\mathbb{U}_n(F)$ is a direct limit of left coherent subrings, satisfying the properties in Lemma \ref{direct limit coherence lemma}, hence is left coherent.\\
Let $U \leq \mathbb{U}_n(F)$ be a closed subgroup. Then by Proposition \ref{closed subgroup coherence prop}, $kU$ is also left coherent. $\square$\\
\end{jpf}

\section{Coherence of a skew polynomial ring}

\subsection{The non-commutative Hilbert basis theorem}

In commutative algebra, the Hilbert basis theorem tells us that if $R$ is a Noetherian ring, then the polynomial ring $R[X]$ is also Noetherian. There is a corresponding result for a nice class of skew polynomial rings, which we call the non-commutative Hilbert basis theorem.

\begin{jthm}[non-commutative Hilbert basis theorem] \h \label{non-commutative Hilbert basis thm}
Let $R$ be a ring. Let $\sigma_X : R \rightarrow R$ be a ring automorphism. If $R$ is left Noetherian, then the skew polynomial ring $R[X ; \sigma_X]$ is left Noetherian. If $R$ is right Noetherian, then $R[X; \sigma_X]$ is right Noetherian.\\
\end{jthm}

A proof is given in \jcite{GoWa04}, Theorem 1.14. Note that in Theorem \ref{non-commutative Hilbert basis thm}, it is enough in fact to assume that $\sigma_X$ is a surjective endomorphism, since any surjective endomorphism of a Noetherian ring is also injective.\\
If we have a non-surjective ring endomorphism, then Theorem \ref{non-commutative Hilbert basis thm} can fail, see Exercise 2P of \jcite{GoWa04}, page 38. However, with appropriate hypotheses on a non-surjective endomorphism, we can conclude that the skew polynomial ring is left coherent.

\begin{jthm} \label{coherent 1 variable thm} \h 
Let $A$ be a left Noetherian ring, and $\sigma_F: A \rightarrow A$ be an injective ring endomorphism, such that $A$ is flat as a right $\sigma_F(A)$-module.
Then, the skew polynomial ring $A[F; \sigma_F]$ is left coherent.\\
\end{jthm}

This result is due to Emerton in \jcite{Eme08}. For the reader's convenience, we provide a proof of Theorem \ref{coherent 1 variable thm} in the next subsection.

\subsection{Coherence of a skew polynomial ring, proof}

In this subsection we prove Theorem \ref{coherent 1 variable thm}. We will often use the convention of omitting the endomorphism, so writing $A[F]$ for $A[F; \sigma_F]$. \vspace{4pt}\\

Let $A$ be a left Noetherian ring and $\sigma_F : A \rightarrow A$ be an injective ring endomorphism. $A$ is naturally a right $\sigma_F(A)$-module, under right multiplication, and we suppose that this right module is flat. We write $R=A[F] = A[F ; \sigma_F]$ for the corresponding skew polynomial ring. \vspace{5pt}\\
The ring $R$ is a free left $A$-module on the basis $\{1,F, F^2, \dots \}$. So
\[ R = A[F] = \bigoplus_{j =0}^{\infty} AF^j . \]
We define an increasing filtration on $R$ given by $F$-degree,
\[R^{\leq k} = \bigoplus_{j =0}^{k} AF^j. \]
This also defines a filtration on any (left) ideal $I$ of $R$, as well as on the free module $R^n$, and any left submodule $M \leq R^n$. 
Let \[ J = RF = \bigoplus_{j =1}^{\infty} AF^j \subseteq R, \] so $J$ is the left ideal generated by $F$. Note that $J$ is also a right ideal of $R$.\\

\begin{jlem} \h \label{J flat lemma}
The ideal $J$ is a flat right $R$-module.\\
\end{jlem}
\begin{jpf} \h
Notice that the endomorphism $\sigma_F$ can be naturally extended to $R$ by defining $\sigma_F(F)=F$, and then $\sigma_F(R) = \sigma_F(A)[F]$. Then, there is an $(R,R)$-bimodule isomorphism 
\[ J=RF \cong R \otimes_{\sigma_F(R)} (FR), \quad rF \mapsto r\otimes F \cdot 1, \quad r\sigma_F(s)F \mapsfrom r \otimes Fs.  \]
Moreover, for any left $R$-module M, there is a natural left $A$-module isomorphism,
\[ R \otimes_{\sigma_F(R)}M \cong A \otimes_{\sigma_F(A)} M, \quad \Big( \sum_{j=0}^\infty a_jF^j \Big)\otimes m \mapsto \sum_{j=0}^\infty a_j \otimes (F^jm), \quad a\otimes m \mapsfrom a \otimes m. \]
Because $A$ is a flat right $\sigma_F(A)$-module, it follows that $R$ is a flat right $\sigma_F(R)$-module. Also, $FR$ is a free right $R$-module because $\sigma_F$ is injective, and so the functor
\[ J \otimes_R = R \otimes_{\sigma_F(R)} (FR) \otimes_R \]
is a composition of two exact functors, hence is exact. So $J$ is a flat right $R$-module. $\square$\\
\end{jpf}

Proving that $R$ is coherent relies on describing the finitely-generated ideals of $R$, and more generally, the finitely-generated submodules of $R^n$.

\begin{jlem}\h \label{M/JM lemma}
Let $M \leq R^n$ be a left $R$-submodule of $R^n$. Then, $M$ is finitely-generated as an $R$-module if and only if $\faktor{M}{JM}$ is finitely-generated as an $A$-module.\\
\end{jlem}

This lemma allows us to show that $R$ is coherent.\\

\begin{jpfof}[Theorem \ref{coherent 1 variable thm}] \h
Let $I \subseteq R$ be a finitely-generated left ideal of $R$. Then we have a short exact sequence
\[ 0 \rightarrow M \rightarrow R^n \rightarrow I \rightarrow 0 .\]
Now, $I$ is finitely-presented if and only if $M$ is finitely-generated (see for example Lemma 2.1.1 of \jcite{Glaz89}). By Lemma \ref{M/JM lemma}, $M$ is finitely-generated if and only if $M/JM = R/J \otimes_R M$ is finitely generated as an $A$-module. Applying $(R/J) \otimes_R$ to the short exact sequence, we obtain the exact sequence
\[ 0 \rightarrow \Tor_1^R(R/J, I) \rightarrow M/JM \rightarrow A^n \rightarrow I/JI \rightarrow 0. \]
Since $A$ is Noetherian, any submodule of $A^n$ is finitely-generated. Thus $M/JM$ is finitely generated as an $A$-module if and only if $\Tor_1^R(R/J, I)$ is.
Now, by dimension shifting,
\[ \Tor_1^R(R/J, I) = \Tor_2^R(R/J, R/I) = \Tor_1^R(J, R/I) .\]
But $J$ is a flat right $R$-module by Lemma \ref{J flat lemma}. Thus $\Tor_1^R(J, R/I) = 0$. Hence $M$ is finitely-generated and $I$ is finitely-presented. \\
So every finitely-generated left ideal of $R$ is finitely-presented -- therefore $R$ is left coherent. $\square$\\
\end{jpfof}

It remains, therefore, to prove Lemma \ref{M/JM lemma}. This will rely on the following description of filtrations.

\begin{jlem} \h \label{filtration JM lemma}
Let $M$ be a left submodule of $R^n$. Then $(JM)^{\leq k} = AFM^{\leq k-1}$ for all $k \geq 1$. \\
\end{jlem}
\begin{jpf}\h
Note $JM = RFM = AFM$ because $M$ is an $A[F]$-module. Let $m \in (JM)^{\leq k}$. Then there are $x_1, x_2, \dots, x_d \in A$ and $m_1, m_2, \dots, m_d \in M$ such that
\[ m = \sum_{j=1}^d x_jFm_j .\]
Let $K \geq k$ be such that $m_1, m_2, \dots, m_d \in M^{\leq K}$.
Denote the coefficient of $F^r$ in the $i$th component of an element $z \in R^n$ by $z^{(i,r)} \in A$. For any $r \in \{k, \dots, K\}$, $i \in \{1, \dots, n\}$,
\[ \sum_{j=1}^d x_jFm_j^{(i,r)}F^{r} = \sum_{j=1}^d x_j\sigma_F(m_j^{(i,r)})F^{r+1} = 0, \]
so \[\sum_{j=1}^d x_j\sigma_F(m_j^{(i,r)}) =0. \]
Now, $A$ is a flat right $\sigma_F(A)$-module. By Lemma 3.65 of \jcite{Rot09}, there exist $b_{qj} = \sigma_F(a_{qj}) \in \sigma_F(A)$, $y_q \in A$ such that for all $j \in \{1,2, \dots, d\}$,
\[ x_j = \sum_{q=1}^N y_qb_{qj}, \]
and for all $q \in \{1,2,\dots, N\}$,
\[ \sum_{j=1}^d b_{qj}\sigma_F\big(m_j^{(i,r)}\big) = \sigma_F\left( \sum_{j=1}^d a_{qj}m_j^{(i,r)}\right)=0. \]
Thus
\[ m = \sum_{j=1}^d x_jFm_j = \sum_{j=1}^d \sum_{q=1}^N y_qb_{qj}Fm_j = \sum_{q=1}^N y_qF \left( \sum_{j=1}^d a_{qj}m_j \right). \]
We define
\[ m_q' = \sum_{j=1}^d a_{qj}m_j \in M, \]
for all $q \in \{1,2,\dots, N\}$. Then, by the properties of the $b_{qj}$ and because $\sigma_F$ is injective,
\[{m_q'}^{(i,r)} = \sum_{j=1}^d a_{qj}m_j^{(i,r)} = 0, \]
for all $i \in \{1,2, \dots, n\}$ and $r \in \{k, \dots, K\}$. Also, $m_q' \in \sum_{j=1}^d Am_j \leq M^{\leq K}$, and therefore $m_q' \in M^{\leq k-1}$, for all $q$. Hence
\[ m = \sum_{q=1}^N y_qFm_q' \in AFM^{\leq k-1}. \]
Therefore, $(JM)^{\leq k} \subseteq AFM^{\leq k-1}$, and the reverse inclusion is obvious. $\square$\\
\end{jpf}

Now we can prove Lemma \ref{M/JM lemma}.\\

\begin{jpfof}[Lemma \ref{M/JM lemma}] \h
Let $M$ be a submodule of $R^n$.\\
If $M$ is finitely generated, there is a surjection $R^N \rightarrow M$. Applying $(R/J)\otimes_R$, we obtain a surjection ${(R/J)}^N = A^N \rightarrow M/JM$, and hence $M/JM$ is finitely generated as an $A$-module. \vspace{3pt}\\
Conversely, suppose that $M/JM$ is finitely generated as an $A$-module. This means that there is a finite set $S \subseteq M$, such that
\[ M = AS + JM .\]
Then, since $S$ is finite, there exists $d \in \mathbb{N}$ such that $S \subseteq M^{\leq d}$. Then
\[ M = M^{\leq d} + JM .\]
We claim that $M = RM^{\leq d}$. If $k \leq d$ then trivially $M^{\leq k} \subseteq RM^{\leq d}$.\\
Let $k >d$, and $m \in M^{\leq k}$. Let $y \in M^{\leq d}$ and $z \in JM$ such that $m=y+z$. Then $z=m-y \in M^{\leq k}$ and $z \in JM$, so $z \in (JM)^{\leq k} =  AFM^{\leq k-1}$, by Lemma \ref{filtration JM lemma}. \\
By induction, we may assume $M^{\leq k-1} \subseteq RM^{\leq d}$, and therefore $z \in RM^{\leq d}$. So $m = y +z \in RM^{\leq d}$. Hence $M^{\leq k} \subseteq RM^{\leq d}$, for all $k$. Therefore $M = RM^{\leq d}$. \\
Now, $M^{\leq d}$ is a submodule of the finitely-generated $A$-module $R^{\leq d}$, hence is finitely generated as an $A$-module, since $A$ is Noetherian. Therefore $M = RM^{\leq d}$ is finitely-generated as an $R$-module. $\square$\\
\end{jpfof}
We conclude that $A[F]$ is a left coherent ring.\\

\subsection{Application to a $p$-adic Lie group}

Recall that in Corollary \ref{unipotent coherent cor}, we showed that all unipotent $p$-adic Lie groups have a coherent augmented Iwasawa algebra. In this subsection, we give an example of a solvable, non-unipotent $p$-adic Lie group with coherent augmented Iwasawa algebra. The reasoning here is a precursor to that for Theorem \ref{general T coherence thm}.\\

Let $G$ be the subgroup of $GL_2(\mathbb{Q}_p)$,
\[ G = \begin{pmatrix} p^{\mathbb{Z}} & \mathbb{Q}_p \\ 0 & 1 \end{pmatrix} = \left\{\begin{pmatrix} p^m & x\\ 0 & 1 \end{pmatrix} \mid m \in \mathbb{Z}, x \in \mathbb{Q}_p \right\}. \]
Let $k$ be a discrete perfect field of characteristic $p$. In this subsection we show that the augmented Iwasawa algebra $kG$ is a coherent ring. We do this by looking explicitly at the ring structure of $kG$.\\

\begin{jlem} \h
The augmented Iwasawa algebra $kG$ is a skew-Laurent polynomial ring
\[ kG = B[F, F^{-1} ; \sigma_F], \]
where $B$ is the subring
\[B =kH = k\begin{pmatrix} 1 & \mathbb{Q}_p \\ 0 & 1 \end{pmatrix}  = \varinjlim_{n \geq 0} k[[t^{\frac{1}{p^n}}]] = k[[t]]^{\frac{1}{p^{\infty}}}, \]
and $F$ corresponds to the group element $\begin{pmatrix} p & 0\\ 0 & 1 \end{pmatrix} \in G$. Then $\sigma_F$ is the $k$-algebra automorphism of $B$ given by $\sigma_F(f(t^a))  = f(t^{pa})$.\\
\end{jlem}

This ring structure is determined by first noting $B = kH$ is naturally a subring of $kG$, by Proposition \ref{augmented Iwasawa algebra functor prop}. Then, we note $kG$ is generated as a $k$-algebra by $kH$ and coset representatives for $G/H$.\\

Let $A$ be the Iwasawa algebra of the compact subgroup $H_0 = \begin{pmatrix} 1 & \mathbb{Z}_p \\ 0 & 1 \end{pmatrix}$, so $A = k[[t]] \leq B$. Note that $FgF^{-1} = g^p$ for all $g \in H_0$ (this corresponds to multiplication by $p$ on $\mathbb{Z}_p$). Thus $\sigma_F$ restricts to a ring endomorphism $\sigma_F'$ of $A$, given by $\sigma_F'(f(t)) = f(t^p)$. Thus the subring of $kG$ generated by $A$ and $F$ is a skew polynomial ring $A[F]=A[F ; \sigma_F']$.

\begin{jprop} \h \label{A[F] coherent example prop}
The skew polynomial ring $A[F]$ is left coherent.\\
\end{jprop}
\begin{jpf} \h
Note that $A=k[[t]]$ is left Noetherian. Because $\sigma_F'= \sigma_F|_A$ and $\sigma_F$ is an automorphism, $\sigma_F'$ is injective. Then, $\sigma_F'(A) = k[[t^p]]$, and so $A$ is free of rank $p$ as a $\sigma_F'(A)$-module, hence flat. Therefore $A[F]$ is left coherent, by Theorem \ref{coherent 1 variable thm}. $\square$\\
\end{jpf}

From this we can use facts about coherent rings to show that $kG$ is coherent.

\begin{jcor} \h
The augmented Iwasawa algebra $kG$ is coherent.\\
\end{jcor}
\begin{jpf} \h
For each $n \geq 0$, let $A_n = kH_n$, where $H_n = \begin{pmatrix} 1 & p^{-n}\mathbb{Z}_p \\ 0 & 1 \end{pmatrix}$ is a compact subgroup. Then, $A_n = k[[t^{\frac{1}{p^n}}]]$, and by similar reasoning to that above, $A_n[F]$ is a skew polynomial ring. Then, $\sigma_F^n|_{A_n}$ is a ring isomorphism onto $A$, and commutes with the action of $F$, therefore $A_n[F] \cong A[F]$ is left coherent by Proposition \ref{A[F] coherent example prop}.\\
Now, $H = \varinjlim_{n \geq 0} H_n$, therefore $B = kH= \varinjlim_{n \geq 0} A_n$. Now, each $H_n$ is a closed subgroup of $H$, therefore $B$ is a flat right $A_n$-module by Proposition \ref{locally profinite flatness prop}. It follows that $B[F]$ is a flat right $A_n[F]$-module, since there is a natural isomorphism of left $B$-modules,
\[ B[F] \otimes_{A_n[F]} M \cong B \otimes_{A_n} M, \quad \Big( \sum_{j=0}^\infty b_jF^j \Big) \otimes m \mapsto \sum_{j=0}^\infty b_j \otimes (F^jm), \quad b\otimes m \mapsfrom b \otimes m, \] 
for any left $A_n[F]$-module $M$. Because $B\otimes_{A_n}$ is an exact functor, so is $B[F] \otimes_{A_n[F]}$.\\
So $B[F] = \varinjlim_{n \geq 0} A_n[F]$ is left coherent, by Lemma \ref{direct limit coherence lemma}. Then, $kG = X^{-1}B[F]$, where ${X = \{F^n \mid n \geq 0 \}}$ is a left denominator set, therefore $kG$ is left coherent, by Corollary \ref{localisations are coherent rings cor}. $\square$\\
\end{jpf}

\section{Non-coherence for two groups}
For the rest of this article, let $k$ be a perfect field of characteristic $p$, given the discrete topology (so $k$ is pseudocompact). For example, $k$ may be a finite extension or an algebraic closure of $\mathbb{F}_p$.\\

In this section we show that the augmented Iwasawa algebras of two particular $p$-adic Lie groups are not coherent. These examples are crucial in proving our characterisation of when augmented Iwasawa algebras are coherent in Theorem \ref{Main result general coherence thm}.\\
Throughout, let $F$ be a finite extension of $\mathbb{Q}_p$ of degree $n$, with discrete valuation $v_F: F^\times \rightarrow \mathbb{Z}$, and ring of integers $\mathcal{O}_F$.

\subsection{Statements}

\begin{jthm} \h \label{kG_3 not coherent thm}
Let $\mathbb{U}_3'$ be the affine group scheme of upper unitriangular matrices in $\mathbb{GL}_3$ with $(1,2)$-entry zero, and
\[U_3' = \mathbb{U}_3'(F) = \begin{pmatrix} 1& 0& F \\ 0& 1 & F \\ 0 &0& 1 \end{pmatrix}. \]
Let $u,v \in F^\times$ be such that $v_F(u) > 0, v_F(v) < 0$, and
\[t' = \begin{pmatrix} u &&\\ & v &\\ && 1 \end{pmatrix} \in GL_3(F). \]
Let $G_{3,t'} = \langle t', U_3' \rangle \leq GL_3(F)$. Then $kG_{3,t'}$ is not coherent.\\
\end{jthm}

\begin{jthm} \h \label{kH_3 not coherent thm}
Let $\mathbb{U}_3$ be the affine group scheme of upper unitriangular matrices in $\mathbb{GL}_3$. Let
\[ U_3 = \mathbb{U}_3(F) = \begin{pmatrix} 1& F& F \\ 0& 1 & F \\ 0 &0& 1 \end{pmatrix}. \]
Let $u, v \in F^\times$ be such that $v_F(u) > 0, v_F(v) < 0$, and
\[ t = \begin{pmatrix} u &&\\ & 1 &\\ && v^{-1} \end{pmatrix} \in GL_3(F). \]
\\
Let $H_{3,t}  =\langle t, U_3 \rangle \leq GL_3(F)$. Then $kH_{3,t}$ is not coherent.\\
\end{jthm}

\subsection{Overview of the proof}

We provide a summary of the arguments used to prove Theorem \ref{kG_3 not coherent thm} and \ref{kH_3 not coherent thm}.\\

We first show Theorem \ref{kG_3 not coherent thm}. We wish to show that the augmented Iwasawa algebra of the $p$-adic Lie group $G_{3,t'} = \langle t', U_3' \rangle$ is not (left) coherent.
\begin{itemize}
\item First, we choose a closed subgroup $G_{3,g} = \langle g, U_3' \rangle \leq G_{3,t'}$, where $g$ is a power of $t'$, the choice made to simplify later calculations.
\end{itemize}
By Proposition \ref{closed subgroup coherence prop}, it is enough to show that $kG_{3,g}$ is not left coherent. \begin{itemize} \item We choose a particular finitely-generated left ideal $I_g \subseteq kG_{3,g}$. We will show it is not finitely-presented. 
\item Next, we find the subgroup $G_{3,D,E} = \langle D,E, U_3' \rangle \leq GL_3(F)$, where $D,E \in GL_3(F)$ are particular diagonal elements such that $DE^{-1}=g$.
\item Then $G_{3,g}$ is a closed subgroup of $G_{3,D,E}$, and $kG_{3,D,E}$ is a (faithfully) flat right $kG_{3,g}$-module, by Theorem \ref{non-compact faithfully flat thm}. \end{itemize}
By the following lemma, to show $I_g$ is not finitely-presented, it will be enough to show that the corresponding left ideal of $kG_{3,D,E}$, which is $I_{D,E} =  kG_{3,D,E} \otimes_{kG_{3,g}} I_g$, is not finitely-presented.
\begin{jlem} \h \label{subring presentation lemma}
Let $R \leq S$ be rings, and $I$ be a finitely-generated left ideal of $R$. If $I$ is finitely-presented as an $R$-module, then $S \otimes_R I$ is finitely-presented as an $S$-module.\\
Moreover, if $S$ is a flat right $R$-module, then $S \otimes_R I$ is naturally identified with a left ideal of $S$.\\
\end{jlem}
\begin{jpf} \h
The first part follows because $S \otimes_R$ is a right exact functor. If $S$ is a flat right $R$-module, then tensoring the natural injection of left $R$-modules $I \rightarrow R$ gives an injection $S \otimes_R I \rightarrow S$, of left $S$-modules. $\square$\\
\end{jpf}
\begin{itemize}
\item Consider the compact unipotent group $\mathbb{U}_3'(\mathcal{O}_F) \leq G_{3,D,E}$, and let $A = k\mathbb{U}_3'(\mathcal{O}_F)$ be its Iwasawa algebra. The subring of $kG_{3,D,E}$ generated by $A$ and the elements $D,E$ is a skew polynomial ring $A[D,E]$. \item Then $kG_{3,D,E}$ is a flat right $A[D,E]$-module, and $I_{D,E}$ can be generated by finitely many elements of $A[D,E]$. \end{itemize} This gives a finitely-generated left ideal $I_A \subseteq A[D,E]$ such that $I_{D,E} = kG_{3,D,E} \otimes_{A[D,E]} I_A$.
\begin{itemize} \item We then compute an infinite set of relations $S$ for $I_A$, over the ring $A[D,E]$. \end{itemize}
The set of relations $S$ for $I_A$ bijectively corresponds to a set of relations $X$ for the ideal $I_{D,E}$.
\begin{itemize} \item We show that no finite subset of $X$ can be a set of relations for $I_{D,E}$. This proves that $I_{D,E}$ cannot be finitely-presented, by the following lemma.\\ \end{itemize}

\begin{jlem} \h \label{finite subset generation lemma}
Let $R$ be a ring, $M$ be a left $R$-module, $X$ be a generating set of $M$. Suppose that $M$ is finitely-generated. Then $M$ is generated by a finite subset of $X$.\\
\end{jlem}
\begin{jpf} \h
Let $M$ have finite generating set $Y$. Each $y \in Y$ can be written as a finite $R$-linear combination of elements in $X$, let
\[ y = r_{1y}x_{1y} + r_{2y}x_{2y} + \dots + r_{a_y y}x_{a_y y}, \qquad a_y \in \mathbb{N}, r_{iy} \in R, x_{iy} \in X. \]
Let $X' = \{ x_{iy} \mid y \in Y, 1 \leq i \leq a_y \} \subseteq X$. Then, $X'$ is a finite subset of $X$, and $Y$ is contained in the module generated by $X'$, so $X'$ must generate $M$. $\square$\\
\end{jpf}
As $I_{D,E}$ is not finitely-presented, it then follows that $I_g$ is not finitely-presented, hence $kG_{3,g}$ is not left coherent, so $kG_{3,t'}$ is not coherent, proving Theorem \ref{kG_3 not coherent thm}.\\

We then deduce Theorem \ref{kH_3 not coherent thm} as follows. \\

We wish to show that $H_{3,t} = \langle t, U_3 \rangle$ is not left coherent.
\begin{itemize} \item We note that $Z(U_3)$ is a normal subgroup of $H_{3,t}$, and that the quotient $H_{3,t}/Z(U_3) \cong G_{3,t'}$. \end{itemize} Then, the augmentation ideal $J = \epsilon(Z(U_3)) \subseteq kH_{3,t}$ is two-sided, and $kH_{3,t}/J \cong kG_{3,t'}$.
\begin{itemize}
\item Because $I_g \subseteq kG_{3,g}$ is not finitely-presented, the corresponding left ideal
\[ I' = kG_{3,t'} \otimes_{kG_{3,g}} I_g \subseteq kG_{3,t'} \]
is not finitely-presented, by faithful flatness. \item We find a finitely-generated left ideal $I \subseteq kH_{3,t}$ such that $I$ contains $J$, we have an isomorphism $I/J \cong I'$, and $JI=J$.
\end{itemize}
From this, we can show that since $I'$ is not finitely-presented, $I$ is not finitely-presented. Therefore $kH_{3,t}$ is not left coherent, and Theorem \ref{kH_3 not coherent thm} is proved.\\

\subsection{Proof of Theorem \ref{kG_3 not coherent thm}: change of rings} \label{Proof of Theorem kG_3 not coherent thm: change of rings}

As in the statement of Theorem \ref{kG_3 not coherent thm}, let $\mathbb{U}_3'$ be the affine group scheme of upper unitriangular matrices in $\mathbb{GL}_3$ with $(1,2)$-entry zero, and
\[U_3' = \mathbb{U}_3'(F) = \begin{pmatrix} 1& 0& F \\ 0& 1 & F \\ 0 &0& 1 \end{pmatrix}. \]
Then let $u,v \in F^\times$ be such that $v_F(u) > 0, v_F(v) < 0$, and
\[t' = \begin{pmatrix} u &&\\ & v &\\ && 1 \end{pmatrix} \in GL_3(F). \]
We define $G_{3,t'} = \langle t', U_3' \rangle \leq GL_3(F)$.\\

Define $n' = v_F(p)$. The field $F$ is a finite extension of $\mathbb{Q}_p$, so $\mathbb{Z}_p$ is a subring of $\mathcal{O}_F$, therefore $v_F|_{\mathbb{Q}_p^\times} = n'v_{\mathbb{Q}_p}$. It follows that
\[u^{n'} = p^{n_u} u_0, \quad v^{n'} = (p^{n_v} v_0)^{-1},\]
for some positive integers $n_u, n_v$ and elements $u_0, v_0 \in \mathcal{O}_F^\times$, and we define
\[ g = {(t')}^{n'} =  \begin{pmatrix} p^{n_u}u_0 &&\\ & (p^{n_v} v_0)^{-1} &\\ && 1 \end{pmatrix}. \]

\begin{jdef} \h \label{G_{3,g} def}
$G_{3,g} = \langle g, U_3' \rangle \leq G_{3,t'}$.\\
\end{jdef}
We describe the corresponding augmented Iwasawa algebra.\\

Now, $\mathcal{O}_F$ is a free $\mathbb{Z}_p$-module of rank $n$, so fix a basis $\{x_0, \dots, x_{n-1}\}$ with $x_0=1$. Then the Iwasawa algebra of $\mathcal{O}_F$ is
\[ k\mathcal{O}_F = k[[ X_0, \dots, X_{n-1}]], \]
where $1+X_i = x_i \in \mathcal{O}_F$.
Then the augmented Iwasawa algebra of $F$ is \[kF = \varinjlim_{r \geq 0} k[[X_0^{\frac{1}{p^r}}, \dots, X_{n-1}^{\frac{1}{p^r}} ]]= k[[X_0, \dots, X_{n-1}]]^{\frac{1}{p^\infty}} , \]
where $1+X_i^{\frac{1}{p^r}} = p^{-r} x_i \in p^{-r}\mathcal{O}_F \leq F$. \\
Now, $U_3' \cong F \oplus F$, and thus we write the augmented Iwasawa algebra as
\[ kU_3' = k[[ s_0, \dots, s_{n-1}, t_0, \dots, t_{n-1} ]]^{\frac{1}{p^\infty}}, \] 
where
\[ 1+s_i = \begin{pmatrix} 1& 0& x_i \\ 0& 1 & 0 \\ 0 &0& 1 \end{pmatrix}, \quad 1+t_i = \begin{pmatrix} 1& 0& 0 \\ 0& 1 & x_i \\ 0 &0& 1 \end{pmatrix}. \]
Then the augmented Iwasawa algebra of $G_{3,g}$ is a skew-Laurent ring $kU_3'[g, g^{-1}]$. We define a finitely-generated left ideal of this ring.
\begin{jdef} \h \label{I_g def}
$ I_g = kG_{3,g} s_0 + \dots + kG_{3,g} s_{n-1} + kG_{3,g} t_0 + \dots + kG_{3,g} t_{n-1} + kG_{3,g}(g-1) \subseteq kG_{3,g}$.\\
\end{jdef}

We will show that $I_g$ is not finitely-presented, by considering the ideal it generates in a larger ring. This larger ring is given by considering a group containing $G_{3,g}$, splitting the element $g$ into its diagonal components.

\begin{jdef} \h
The group $G_{3,D,E} = \langle D, E, U_3' \rangle \leq GL_3(F)$, where
\[ D = \begin{pmatrix} p^{n_u}u_0 &&\\ &1 &\\ && 1 \end{pmatrix}, \quad E = \begin{pmatrix} 1 &&\\ & p^{n_v} v_0 &\\ && 1 \end{pmatrix} \in GL_3(F). \]
\end{jdef}

Notice that $g=DE^{-1}$, and so $G_{3,g}$ is a closed subgroup of $G_{3,D,E}$. We have that $kG_{3,D,E}$ is a skew-Laurent ring in two variables, that is,
\[ kG_{3,D,E} = kU_3'[D,D^{-1}, E, E^{-1}] = kU_3'[D,D^{-1}, E, E^{-1} ; \sigma_D, \sigma_E ], \]
where $\sigma_D, \sigma_E$ are the automorphisms of $kU_3'$ given by $\sigma_D(x) = DxD^{-1}, \sigma_E(x) = ExE^{-1}$. Note that the group elements $D,E$ commute, therefore $\sigma_D, \sigma_E$ are commuting endomorphisms.

\begin{jdef} \h
The left ideal $I_{D,E} = kG_{3,D,E} \otimes_{kG_{3,g}} I_g \subseteq kG_{3,D,E}$.\\
\end{jdef}

By Lemma \ref{subring presentation lemma}, if $I_g$ is finitely-presented, then so is $I_{D,E}$. Since $D,E$ are units of $kG_{3,D,E}$, we have that $I_{D,E}$ is generated by the set $\{s_0, \dots, s_{n-1}, t_0, \dots, t_{n-1}, D-E \}$. Let the natural presentation corresponding to these generators be
\[ 0 \rightarrow K_{D,E} \rightarrow (kG_{3,D,E})^{2n+1} \rightarrow I_{D,E} \rightarrow 0 .\]

We will determine a set of generators for $K_{D,E}$, and then demonstrate that this set cannot be reduced to a finite one. To determine a set of generators for $K_{D,E}$, we will restrict the calculations to a skew polynomial subring of $kG_{3,D,E}$.
\begin{jdef} \h
The ring $A = k\mathbb{U}_3'(\mathcal{O}_F)  = k \begin{pmatrix} 1& 0& \mathcal{O}_F \\ 0& 1 & \mathcal{O}_F \\ 0 &0& 1 \end{pmatrix}$.\vspace{4pt}\\
\end{jdef}

Then $A$ is the Noetherian subring $k[[s_0, \dots, s_{n-1}, t_0, \dots, t_{n-1} ]] \leq kG_{3,D,E}$. Now, it can be easily seen that $\sigma_D(A), \sigma_E(A) \subseteq A$, therefore $\sigma_D|_A, \sigma_E|_A$ are (injective) ring endomorphisms of $A$, and so the subring of $kG_{3,D,E}$ generated by $D,E$ and $A$ is a skew polynomial ring $A[D,E] = A[D, E ; \sigma_D|_A, \sigma_E|_A]$. Also, $kG_{3,D,E}$ is flat over this subring. To prove this requires a minor result about localisations.

\begin{jlem} \h \label{localisation of flat is flat lemma}
Let $R$ be a ring, $X \subseteq R$ be a left denominator set. Let $M$ be a $(R,R)$-bimodule such that $M$ is  a flat right $R$-module. Then $X^{-1}M$ is a flat right $R$-module.\\
\end{jlem}
\begin{jpf} \h
By Proposition 10.12 of \jcite{GoWa04}, the functor $X^{-1}M \otimes_R$ is equal to the composition of functors $X^{-1}R \otimes_R M \otimes_R$. Since $M$ is a flat right $R$-module, $M \otimes_R$ is an exact functor. By Corollary 10.13 of \jcite{GoWa04}, $X^{-1}R$ is a flat right $R$-module, so $X^{-1}R \otimes_R$ is an exact functor. Thus $X^{-1}M \otimes_R$ is an exact functor, hence $X^{-1}M$ is a flat right $R$-module. $\square$\\
\end{jpf}

\begin{jprop} \h
The augmented Iwasawa algebra $kG_{3,D,E}$ is a flat right $A[D,E]$-module. \\
\end{jprop}
\begin{jpf} \h
Because $\mathbb{U}_3'(\mathcal{O}_F) \leq U_3'$ is a closed subgroup, $kU_3'$ is a flat $A$-module by Theorem \ref{non-compact faithfully flat thm}. For any left $A[D,E]$-module $M$, there is a natural isomorphism of left $kU_3'$-modules,
\[ kU_3'[D,E] \otimes_{A[D,E]} M \cong kU_3' \otimes_A M, \quad \Big( \sum_{i,j \geq 0} b_{ij}D^iE^j \Big) \otimes m \mapsto \sum_{i,j \geq 0} b_{ij} \otimes (D^iE^jm), \; b\otimes m \mapsfrom b \otimes m. \]
Since $kU_3' \otimes_A$ is exact, so is $kU_3'[D,E] \otimes_{A[D,E]}$, thus $kU_3'[D,E]$ is a flat right $A[D,E]$-module. By Lemma \ref{localisation of flat is flat lemma}, the localisation $X^{-1}kU_3'[D,E]$ is also a flat right $A[D,E]$-module, where ${X = \{D^aE^b \mid a,b \geq 0 \}}$ is a left denominator set. But $kG_{3,D,E} \cong X^{-1}kU_3'[D,E]$. $\square$ \\
\end{jpf}

Now, notice that $I_{D,E}$ is generated by elements which all lie in $A[D,E]$. So we can consider the corresponding left ideal of $A[D,E]$.

\begin{jdef} \h
$I_A = A[D,E]s_0 + \dots A[D,E]s_{n-1} + A[D,E]t_0 + \dots + A[D,E]t_{n-1} + A[D,E](D-E) \subseteq A[D,E]$. \\
\end{jdef}

Then, we have that $kG_{3,D,E} \otimes_{A[D,E]} I_A = I_{D,E}$, since $kG_{3,D,E}$ is a flat $A[D,E]$-module. Let the natural presentation of $I_A$ be
\[ 0 \rightarrow K_A \rightarrow A[D,E]^{2n+1} \rightarrow I_A \rightarrow 0. \]
Again by flatness, it follows that $kG_{3,D,E} \otimes_{A[D,E]} K_A = K_{D,E}$. So any set of generators of $K_A$ gives a corresponding set of generators of $K_{D,E}$.\\
Let
\[ \Pi_{D,E}: K_{D,E} \rightarrow (kG_{3,D,E})^{2n}, \quad \Pi_{A}: K_{A} \rightarrow A[D,E]^{2n}\]
be the maps giving projection onto the first $2n$ coordinates, and let
\[L_{D,E} = \textrm{Im }\Pi_{D,E}, \quad L_{A} = \textrm{Im }\Pi_{A}. \]
Then $\Pi_{D,E}, \Pi_A$ are injective, since $D-E$ is not a zero divisor in $kG_{3,D,E}$. So $K_{D,E} \cong L_{D,E}$ and $K_A \cong L_A$. It is clear that $kG_{3,D,E} \otimes_{A[D,E]} \Pi_A = \Pi_{D,E}$, and thus\\ $kG_{3,D,E} \otimes_{A[D,E]} L_A = L_{D,E}$. So any set of generators of $L_A$ bijectively corresponds to a set of generators of $L_{D,E}$, in turn corresponding to a set of generators for $K_{D,E}$.

\[
\begin{tikzcd}
{K_{D,E}} \arrow[r, "{\Pi_{D,E}}"]                                                   & {L_{D,E}}                                                             \\
K_A \arrow[r, "\Pi_A"] \arrow[u, "{kG_{3,D,E} \otimes_{A[D,E]}}", dotted, bend left] & L_{A} \arrow[u, "{kG_{3,D,E} \otimes_{A[D,E]}}"', dotted, bend right]
\end{tikzcd}
\]

In the next subsection we calculate a set of generators for $L_A$, and use this information to show that in fact $K_{D,E}$ cannot be finitely-generated.\\

\subsection{Proof of Theorem \ref{kG_3 not coherent thm}: calculation of generators}

The following proposition gives a set of generators for $L_A$. This information allows us to prove Theorem \ref{kG_3 not coherent thm}, see subsection \ref{Proof of Theorem kG_3 not coherent thm: conclusion}.

\begin{jprop} \h \label{L_A generating set prop}
The left $A[D,E]$-module $L_A$ has a generating set $S = S_1 \cup S_2 \cup S_3$, where $S_1, S_2, S_3$ are as follows.\\
There exist elements $b_{ji}, c_{ji} \in A$ for $i,j \in \{0, \dots, n-1 \}$, such that
\[ b_j = (-b_{j0}E, -b_{j1}E, \dots, D - b_{jj}E, \dots, -b_{j(n-1)}E, 0^n) \in L_A, \]
and
\[ c_j = (0^n, -c_{j0}D, -c_{j1}D, \dots, E - c_{jj}D, \dots, -c_{j(n-1)}D) \in L_A, \]
for all $j \in \{0, \dots, n-1 \}$. Then
\[ S_1 = \{b_j \mid j \in \{0, \dots, n-1 \} \}, \]
and
\[ S_2 = \{ c_j \mid j \in \{0, \dots, n-1 \} \}. \] 
Then,
\[ S_3 =  \{(0 \dots, t_jE^m, \dots, 0, 0, \dots, -s_iD^m, \dots, 0) \mid m \in \mathbb{Z}_{\geq 0}, 0 \leq i,j \leq n-1 \}, \]
where $t_jE^m$ is in the $(i+1)$st place, and $-s_iD^m$ is in the $(n+j+1)$st place, of the $2n$-tuple.\\
\end{jprop}

The calculation of these generators depends on understanding the action of $D,E$ on certain augmentation ideals.

\begin{jdef} \h
The affine group scheme $\mathbb{U}_{1,3}'$ is the group scheme of matrices in $\mathbb{U}_3'$ with $(2,3)$-entry zero.\\
The affine group scheme $\mathbb{U}_{2,3}'$ is the group scheme of matrices in $\mathbb{U}_3'$ with $(1,3)$-entry zero.\\
\end{jdef}

\begin{jlem} \h \label{explicit augmentation lemma}
Let $m \in \mathbb{Z}$. The subring of $kU_3'(F)$ corresponding to the Iwasawa algebra of 
\[\mathbb{U}_{1,3}'(p^m\mathcal{O}_F) = \begin{pmatrix} 1& 0& p^m\mathcal{O}_F \\ 0& 1 & 0 \\ 0 &0& 1 \end{pmatrix} \]
is the subring $k[[s_0^{p^m}, \dots, s_{n-1}^{p^m}]]$. The corresponding augmentation ideal
$ \epsilon_{\mathbb{U}_{1,3}'(p^m\mathcal{O}_F)}( \mathbb{U}_{1,3}'(p^m\mathcal{O}_F) ) $
is generated by $\{s_0^{p^m}, \dots, s_{n-1}^{p^m} \}$.\\
An identical result holds for $\mathbb{U}_{2,3}'(p^m\mathcal{O}_F)$.\\
\end{jlem}

This lemma is clear from the description of $kF$ in subsection \ref{Proof of Theorem kG_3 not coherent thm: change of rings}. It follows, by the description of augmentation ideals in Lemma \ref{augmentation generation lemma}, that if $\mathbb{U}_{1,3}'(p^m\mathcal{O}_F)$ is a subgroup of any $p$-adic Lie group $G$, then $\epsilon_G(\mathbb{U}_{1,3}'(p^m\mathcal{O}_F))$ is the left ideal generated by $\{s_0^{p^m}, \dots, s_{n-1}^{p^m} \}$.\\

We can now prove Proposition \ref{L_A generating set prop}, which we do in two parts.

\begin{jlem} \h \label{L_A finitely many generators lemma} 
There exist elements $b_{ji}, c_{ji} \in A$ for $i,j \in \{0, \dots, n-1 \}$, such that
\[ b_j = (-b_{j0}E, -b_{j1}E, \dots, D - b_{jj}E, \dots, -b_{j(n-1)}E, 0^n) \in L_A, \]
and
\[ c_j = (0^n, -c_{j0}D, -c_{j1}D, \dots, E - c_{jj}D, \dots, -c_{j(n-1)}D) \in L_A, \]
for all $j \in \{0, \dots, n-1 \}$.\\
\end{jlem}
\begin{jpf} \h
Note that $L_A$ has the following description:
\[ L_A = \left\{ (\lambda_0, \dots, \lambda_{n-1}, \mu_0, \dots, \mu_{n-1}) \in A[D,E]^{2n} : \sum_{j \geq 0}^{n-1} \lambda_js_j + \mu_jt_j \in A[D,E](D-E) \right\}. \]
Let $j \in \{0, \dots, n-1 \}$. Then
\[ Ds_j = \sigma_D(s_j)D = \sigma_D(s_j)E + \sigma_D(s_j)(D-E), \]
so
\[ Ds_j - \sigma_D(s_j)E \in A[D,E](D-E). \]
Now, notice that $s_j \in \epsilon_{\mathbb{U}_{3}'(\mathcal{O}_F)}( \mathbb{U}_{1,3}'(\mathcal{O}_F) ) $ and that $\sigma_D$ is a ring endomorphism of $A=k\mathbb{U}_{3}'(\mathcal{O}_F)$ corresponding to a group endomorphism of $\mathbb{U}_{3}'(\mathcal{O}_F)$. Moreover, the image of $\mathbb{U}_{1,3}'(\mathcal{O}_F)$ under this group homomorphism is $\mathbb{U}_{1,3}'(p^{n_u}\mathcal{O}_F)$. By Lemma \ref{augmentation endomorphism lemma}, it follows that \[ \sigma_D(s_j) \in \epsilon_{\mathbb{U}_{3}'(\mathcal{O}_F)}( \sigma_D(\mathbb{U}_{1,3}'(\mathcal{O}_F)) ) = \epsilon_{\mathbb{U}_3'(\mathcal{O}_F)}( \mathbb{U}_{1,3}'(p^{n_u}\mathcal{O}_F) ). \]
So by Lemma \ref{explicit augmentation lemma},
\[ \sigma_D(s_j) \in As_0^{p^{n_u}}+ \dots+ As_{n-1}^{p^{n_u}} \subseteq As_0+ \dots+ As_{n-1}. \vspace{3pt} \] 
So, there exist $b_{j0}, b_{j1}, \dots, b_{j(n-1)} \in A$ such that $\sigma_D(s_j) = \sum_{0 \leq i \leq n-1} b_{ji}s_i$, and then
\begin{align*}Ds_j - \sigma_D(s_j)E &= Ds_j - \sum_{0 \leq i \leq n-1} b_{ji}s_iE\\ &= Ds_j - \sum_{0 \leq i \leq n-1} b_{ji}Es_i\\ &= (D - b_{jj}E)s_j + \sum_{0 \leq i \leq n-1, i \neq j} -b_{ji}E s_i \in A[D,E](D-E) .\end{align*}
So $b_j = (-b_{j0}E, -b_{j1}E, \dots, D - b_{jj}E, \dots, -b_{j(n-1)}E, 0^n) \in L_A$, for all $j$.\\
Identical reasoning with the elements $E$, $t_j$ similarly gives the tuples $c_j \in L_A$. $\square$\\
\end{jpf}

\begin{jpfof}[Proposition \ref{L_A generating set prop}] \h
We can put a grading by total $(D,E)$-degree on the ring $A[D,E]$, and in this grading, each of the elements $s_0, \dots, s_{n-1}, t_0, \dots, t_{n-1}$ is homogeneous of degree zero, and $D-E$ is homogeneous of degree 1. Therefore every element of $L_A$ will be a sum of elements $(x_0, \dots, x_{n-1}, y_0, \dots, y_{n-1})$ where each $x_j, y_j$ is homogeneous of the same degree. \\
Let $L_A'$ be the $A[D,E]$-submodule generated by $S_1 \cup S_2$, as given by Lemma \ref{L_A finitely many generators lemma}.\\
We have that for any $j$,
\[ b_j = (0, \dots, D, \dots, 0, 0^n ) - (b_{j0}E, b_{j1}E, \dots, b_{j(n-1)}E, 0^n) \in L_A', \]
so for any $j$,
\[ (0, \dots, D, \dots, 0, 0^n ) \in L_A' + (A[E]^n, 0^n). \]
Now, because $DA[E] = \sigma_D(A)[E]D \subseteq A[E]D$, left multiplying by $D$ this equation gives that for any $j$, 
\[ (0, \dots, D^2, \dots, 0, 0^n ) \in L_A' + (A[E]^nD, 0^n) = L_A' + (A[E]^n, 0^n), \]
since $L_A'$ is a left $A[D,E]$-module. Continuing by induction we see that 
\[ (0, \dots, D^m, \dots, 0, 0^n ) \in L_A' + (A[E]^n, 0^n) \]
for any $m \geq 0$. Identical reasoning applies to show that
\[ ( 0^n, 0, \dots, E^m, \dots, 0, ) \in L_A' + ( 0^n, A[D]^n). \vspace{3pt} \]
It follows that if $x,y \in A[D,E]^n$ and $(x,y) \in L_A$, then $(x,y) + L_A' = (x', y') + L_A'$ for some $x' \in A[E]^n$ and $y' \in A[D]^n$. (Informally, the relations above give a way to replace any $D$s appearing in the first $n$ places with polynomials in $E$, and vice versa for the last $n$ places.)\\
Then $(x', y') \in L_A$ and is the sum of its homogeneous parts, which also lie in $L_A$, as discussed above. Therefore, $L_A/L_A'$ is generated by (the image of) elements $(x'', y'') \in L_A$, homogeneous of degree $m$, with $x'' \in (AE^m)^n, y'' \in (AD^m)^n$. \vspace{2pt}\\
Now, for $x_j', y_j' \in A$,
\begin{align*}  (x_0'E^m, \dots, x_{n-1}'E^m, y'_0D^m, \dots, y_{n-1}'D^m)  \in L_A &\Leftrightarrow \sum_{j=0}^{n-1} x_j'E^ms_j + y_j'D^mt_j \in A[D,E](D-E)\\ &\Leftrightarrow \sum_{j=0}^{n-1} x_j's_jE^m + y_j't_jD^m \in A[D,E](D-E) \\ &\Leftrightarrow \sum_{j=0}^{n-1} x_j's_j + y_j't_j = 0,
\end{align*}
since $D^m - E^m \in A[D,E](D-E)$.\\
Then, a generating set for $\Big\{   (x_0', \dots, x_{n-1}', y'_0, \dots, y_{n-1}') \in A^{2n} \mid \sum_{j=0}^{n-1} x_j's_j + y_j't_j = 0 \Big\}$ is
\[ S_3' = \{(0 \dots, t_j, \dots, 0, 0, \dots, -s_i, \dots, 0) \mid 0 \leq i,j \leq n-1 \}, \]
and thus a generating set for $\faktor{L_A}{L_A'}$ is (the image of)
\[ S_3 =  \{(0 \dots, t_jE^m, \dots, 0, 0, \dots, -s_iD^m, \dots, 0) \mid m \geq 0, 0 \leq i,j \leq n-1 \}. \]
It follows that $L_A$ is generated by $S = S_1 \cup S_2 \cup S_3$, as required. $\square$\\
\end{jpfof}

\subsection{Proof of Theorem \ref{kG_3 not coherent thm}: conclusion} \label{Proof of Theorem kG_3 not coherent thm: conclusion}

We can now finish proving Theorem \ref{kG_3 not coherent thm}. We have a set $S$ of generators for $L_A$, and hence for $L_{D,E} \cong K_{D,E}$. Using the set $S$ we show that a certain quotient of $L_{D,E}$ is not finitely-generated, implying $K_{D,E}$ cannot be finitely-generated, and the theorem follows.\\

Consider the left module $V =\faktor{kG_{3,D,E}}{kG_{3,D,E}(D-E)}$. Note that $V$ is a cyclic $kG_{3,D,E}$-module, and is a free left $kU_3'$-module, with basis $\{z_m = E^m + kG_{3,D,E}(D-E) \mid m \in \mathbb{Z} \}$.\vspace{2pt}\\
Let $q: L_{D,E} \rightarrow V$ be the left $kG_{3,D,E}$-module homomorphism given by
\[ q(\lambda_0, \dots, \lambda_{n-1}, \mu_0, \dots, \mu_{n-1}) = \sum_{j=0}^{n-1}\lambda_js_j + kG_{3,D,E}(D-E).  \]
Let $M_{D,E} = \textrm{Im }q$, so $M_{D,E}$ is a submodule of $V$. We will prove the following.

\begin{jprop} \h \label{M_{D,E} not fg prop}
The left $kG_{3,D,E}$-module $M_{D,E}$ is not finitely-generated.\\
\end{jprop}

From this we deduce the theorem.\\

\begin{jpfof}[Theorem \ref{kG_3 not coherent thm}] \h
Suppose $kG_{3,t'}$ is left coherent. Then, since $G_{3,g} \leq G_{3,t'}$ is a closed subgroup, $kG_{3,g}$ is left coherent by Proposition \ref{closed subgroup coherence prop}. Now, $I_g$ is a finitely-generated left ideal of $kG_{3,g}$, and so it is also finitely-presented. Then by Lemma \ref{subring presentation lemma}, it follows that the left ideal $I_{D,E}$ of $kG_{3,D,E}$ is finitely-presented, and therefore its relation module $K_{D,E}$ is finitely-generated. Then, $K_{D,E} \cong L_{D,E}$, and $M_{D,E}$ is a quotient of $L_{D,E}$, hence is finitely-generated. But this is a contradiction, by Proposition \ref{M_{D,E} not fg prop}. So $kG_{3,t'}$ is not left coherent, and so also not right coherent. $\square$ \\
\end{jpfof}

Now we prove the proposition.

\begin{jlem} \h \label{M_{D,E} specific generating set lemma}
The set $T = \{s_it_jz_m \mid i,j \in \{0, \dots, n-1 \}, m \geq 0 \}$ generates $M_{D,E}$ as a left $kG_{3,D,E}$-module.\\
\end{jlem}
\begin{jpf} \h
We have a generating set $S=S_1 \cup S_2 \cup S_3$ for $L_{D,E}$, and a surjection $q:L_{D,E} \rightarrow M_{D,E}$, so $M_{D,E}$ is generated by $q(S)$. It is easy to compute that $q(S_1) = q(S_2) = \{0\}$ and $q(S_3) = T$. $\square$\\
\end{jpf}

\begin{jpfof}[Proposition \ref{M_{D,E} not fg prop}]  \h
Suppose, for a contradiction, that $M_{D,E}$ is finitely-generated. In Lemma \ref{M_{D,E} specific generating set lemma}, we have an infinite generating set $T$ for $M_{D,E}$. By Lemma \ref{finite subset generation lemma}, a finite subset of $T$ generates $M_{D,E}$, in particular there exists a non-negative integer $N$ such that
\[ T_N = \{ s_it_jz_m \mid i,j \in \{0, \dots, n-1 \}, 0 \leq m \leq N \} \]
generates $M_{D,E}$.\vspace{2pt}\\
Then $s_{0}t_{0}z_{N+1} \in M_{D,E} = \langle T_N \rangle$. Now, $kG_{3,D,E} = kU_3'[D,D^{-1}, E, E^{-1}]$, and so
 there exist $N' \in \mathbb{N}$, $c_{abijm} \in kU_3'$, such that
\[ s_{0}t_{0}z_{N+1} = \sum_{m=0}^N \sum_{-N' \leq a,b \leq N'} \sum_{i,j \in \{0, \dots, n-1\}} c_{abijm}D^aE^bs_it_jz_m. \]
Now, $D^aE^bz_m = z_{m+a+b}$ for any $m,a,b$, so we have that
\begin{align*} s_{0}t_{0}z_{N+1} &= \sum_{m=0}^N \sum_{-N' \leq a,b \leq N'} \sum_{i,j \in \{0, \dots, n-1\}} c_{abijm}\sigma_D^a\sigma_E^b(s_it_j)D^aE^bz_m \\ &= \sum_{m=0}^N \sum_{-N' \leq a,b \leq N'} \sum_{i,j \in \{0, \dots, n-1\}} c_{abijm}\sigma_D^a\sigma_E^b(s_it_j)z_{m+a+b}. \end{align*}
Because $\{z_m \mid m \in \mathbb{Z} \}$ is a $kU_3'$-basis of $V$, we can compare the coefficients of $z_{N+1}$ to obtain
\[ s_{0}t_{0} = \sum_{\substack{-N' \leq a,b \leq N' \\ 1 \leq a+b \leq N+1}}  \sum_{i,j \in \{0, \dots, n-1\}} c_{abij(N+1-a-b)}\sigma_D^a\sigma_E^b(s_it_j). \]
Therefore $s_0t_0$ lies in the ideal of $kU_3'$ generated by 
\[ \{ \sigma_D^a\sigma_E^b(s_it_j) \mid a+b \geq 1, \; 0 \leq i,j \leq n-1 \}. \]
Then, $\sigma_D^a\sigma_E^b(s_it_j) = \sigma_D^a(s_i)\sigma_E^b(t_j)$ since $D,t_j$ and $E, s_i$ commute. Now,
\[ s_i \in \epsilon_{U_3'}(\mathbb{U}_{1,3}'(\mathcal{O}_F)), \quad t_j \in \epsilon_{U_3'}(\mathbb{U}_{2,3}'(\mathcal{O}_F)), \]
as seen in Lemma \ref{explicit augmentation lemma}, and
\[ \sigma_D(\mathbb{U}_{1,3}'(\mathcal{O}_F) )= \mathbb{U}_{1,3}'(p^{n_u}\mathcal{O}_F), \quad \sigma_E(\mathbb{U}_{2,3}'(\mathcal{O}_F) )= \mathbb{U}_{2,3}'(p^{n_v}\mathcal{O}_F). \]
So, again by Lemma \ref{explicit augmentation lemma}, the element $\sigma_D^a(s_i)\sigma_E^b(t_j)$ lies in the ideal of $kU_3'$,
\[ J_{ab} = (s_0^{p^{an_u}}, \dots, s_{n-1}^{p^{an_u}})(t_0^{p^{bn_v}}, \dots, t_{n-1}^{p^{bn_v}}). \]
It follows that $s_0t_0$ lies in the ideal
\[ \sum_{\substack{a,b \in \mathbb{Z},\\a+b \geq 1}} J_{ab} = \sum_{ \substack{a,b \in \mathbb{Z}, \\ a \geq 1, b \geq 1-a}} J_{ab} +  \sum_{ \substack{a,b \in \mathbb{Z}, \\ b \geq 1, a \geq 1-b}}J_{ab}. \]
Since $J_{ab} \subseteq J_{a'b'} \Leftrightarrow a \geq a', b \geq b'$, it follows that $s_0t_0$ lies in a finite sum of ideals,
\[ s_0t_0 \in \sum_{ a \geq 1} J_{ab'} +  \sum_{b \geq 1}J_{a'b} = J_{1b'} + J_{a'1} \]
for some integers $a', b' \in \mathbb{Z}$. So
\[ s_0t_0 \in  (s_0^{p^{n_u}}, \dots, s_{n-1}^{p^{n_u}})(t_0^{p^{b'n_v}}, \dots,  t_{n-1}^{p^{b'n_v}}) +  (s_0^{p^{a'n_u}}, \dots,  s_{n-1}^{p^{a'n_u}})(t_0^{p^{n_v}}, \dots, t_{n-1}^{p^{n_v}}). \]
But, since $n_u, n_v > 0$, this cannot occur, thus we have a contradiction. So $M_{D,E}$ is not finitely-generated. $\square$ \vspace{4pt}\\
\end{jpfof}

\subsection{Proof of Theorem \ref{kH_3 not coherent thm}} \label{Proof of Theorem kH_3 not coherent thm}

Recall the objects in the statement of Theorem \ref{kH_3 not coherent thm}.\\
Let $\mathbb{U}_3$ be the affine group scheme of upper unitriangular matrices in $\mathbb{GL}_3$, and $U_3 = \mathbb{U}_3(F)$, so
\[ U_3 =  \begin{pmatrix} 1& F& F \\ 0& 1 & F \\ 0 &0& 1 \end{pmatrix}. \]
As in the above subsections, let $u, v \in F^\times$ be such that $v_F(u) > 0, v_F(v) < 0$. Then let
\[ t = \begin{pmatrix} u &&\\ & 1 &\\ && v^{-1} \end{pmatrix} \in GL_3(F). \]
We define $H_{3,t}  =\langle t, U_3 \rangle \leq GL_3(F)$.\vspace{2pt}\\

Note that
\[Z(U_3) = \begin{pmatrix}1 & 0& F \\0 & 1 & 0\\ 0&0& 1 \end{pmatrix} \]
is a normal subgroup of $H_{3,t}$. Moreover, the quotient group is
\[ \faktor{H_{3,t}}{Z(U_3)} \cong G_{3,t'}, \]
as seen in Theorem \ref{kG_3 not coherent thm}, where $t' = \textrm{diag}(u,v,1)$. Then, by Proposition \ref{surjection of groups prop}, we have an isomorphism of augmented Iwasawa algebras,
\[ \phi: \faktor{kH_{3,t}}{J} \xrightarrow{\cong} kG_{3,t'}, \]
where $J = \epsilon_{H_{3,t}}(Z(U_3))$ is the augmentation ideal of $Z(U_3)$, which is two-sided.\\
Let us define notation for the elements of the augmented Iwasawa algebras. We have that
\[ U_3 = \varinjlim_{r \geq 0}  \begin{pmatrix} 1&p^{-r} \mathcal{O}_F&p^{-2r}O_ F \\ 0& 1 &p^{-r} \mathcal{O}_F \\ 0 &0& 1 \end{pmatrix}, \]
and so we write the augmented Iwasawa algebra as
\begin{align*} kU_3 &= \varinjlim_{r \geq 0} k[[ s_0^{\frac{1}{p^r}}, \dots, s_{n-1}^{\frac{1}{p^r}}, t_0^{\frac{1}{p^r}}, \dots, t_{n-1}^{\frac{1}{p^r}}, w_0^{\frac{1}{p^{2r}}}, \dots, w_{n-1}^{\frac{1}{p^{2r}}} ]]\\ &= k[[ s_0, \dots, s_{n-1}, t_0, \dots, t_{n-1}, w_0, \dots, w_{n-1} ]]^{\frac{1}{p^\infty}} , 
\end{align*}
where
\[ 1+s_i = \begin{pmatrix} 1& x_i& 0 \\ 0& 1 & 0 \\ 0 &0& 1 \end{pmatrix}, \quad 1+t_i = \begin{pmatrix} 1& 0& 0 \\ 0& 1 & x_i \\ 0 &0& 1 \end{pmatrix}, \quad 1+w_i = \begin{pmatrix} 1& 0& x_i \\ 0& 1 & 0 \\ 0 &0& 1 \end{pmatrix}. \ \]
Note that $kU_3$ is non-commutative -- it is a direct limit of non-commutative power series rings -- but the $w_j$ commute with any $s_i, t_i$.\\
The augmentation ideal of $Z(U_3)$ is generated by $\{ w_0^{\frac{1}{p^r}}, \dots, w_{n-1}^{\frac{1}{p^r}}  \mid r \geq 0 \}$. \vspace{4pt}\\
Then $kH_{3,t} = kU_3[t, t^{-1}]$ is a skew-Laurent ring over $kU_3$, and we have $\phi(t) = t'$, and ${\phi(s_i) = s_i, \phi(t_j) = t_j}$ with the labelling described above. Let $y = t^{n'}$, where $n'$ is the positive integer as defined at the beginning of subsection \ref{Proof of Theorem kG_3 not coherent thm: change of rings}, so that $\phi(y) = g = (t')^{n'} \in G_{3,t'}$.\\

Consider the left ideal of $kH_{3,t}$,
\[ I = kH_{3,t}s_0 + \dots + kH_{3,t}s_{n-1} + kH_{3,t}t_0 + \dots + kH_{3,t}t_{n-1} + kH_{3,t}(1-y). \]
Let $I'$ be the left ideal of $kG_{3,t'}$ given by
\[ I' = kG_{3,t'} \otimes_{kG_{3,g}} I_g, \]
that is explicitly,
\[ I' = kG_{3,t'}s_0 + \dots + kG_{3,t'}s_{n-1} + kG_{3,t'}t_0 + \dots + kG_{3,t'}t_{n-1} + kG_{3,t'}(1-g). \]

\begin{jlem} \h \label{I contains J lemma}
There is containment of ideals $J \subseteq I$, and $I/J$ corresponds to the left ideal $I'$ under the isomorphism $kH_{3,t}/J \cong kG_{3,t'}$.\\
\end{jlem}
\begin{jpf} \h
For all $x \in I$, $x - x(1-y)  = xy \in I$, so $y^{-1}xy \in I$. Similarly, $y^{-1}(1-y) = y^{-1}-1 \in I$, and so $y(x + x(y^{-1}-1)) = yxy^{-1} \in I$. Write $\sigma_y$ for the ring automorphism of $kU_3$ given by $\sigma_y(x) = yxy^{-1}$.\\
Then, by almost identical calculations as for $D,E$ and $g=DE^{-1} \in kG_{3,t'}$, we can compute how $\sigma_y$ acts on the elements, and ideals, of $kU_3$. We have that
 \[ \sigma_y(kU_3t_0+ \dots + kU_3t_{n-1} ) = kU_3t_0^{p^{-n_v}} + \dots + kU_3t_{n-1}^{p^{-n_v}}. \]
Therefore $t_0^{p^{-r}} \in I$ for all $r \geq 0$. Then, a matrix calculation gives that
\[ (1+s_j)(1+t_0^{p^{-r}} )(1+s_j)^{-1}(1+t_0^{p^{-r}} )^{-1} = 1+ w_j^{p^{-r}}. \]
(This uses that the element $x_0$ of the $\mathbb{Z}_p$-basis for $\mathcal{O}_F$ is chosen to be $x_0=1$.) We can rearrange to find that
\[ s_jt_0^{p^{-r}} - t_0^{p^{-r}}s_j = (1+s_j)(1+t_0^{p^{-r}})w_j^{p^{-r}}. \]
Since $ 1+s_j, 1+t_0^{p^{-r}}$ are units in $kU_3$, it follows that $w_j^{p^{-r}} \in I$, for any $j$, $r$.\\
Thus, $I$ contains $J$, which is generated by such elements.\\
It is clear that $g$ corresponds to $y$ and the $s_j$, $t_j$ correspond under the isomorphism ${ kH_{3,t}/J \cong kG_{3,t'} }$, therefore $I/J = I'$. $\square$\\
\end{jpf}

\begin{jlem} \h
The product of ideals $JI = J$.\\
\end{jlem}
\begin{jpf} \h
It is clear that $J^p = J$, and hence $J = J^2$. Since $J \subseteq I$, we have $ J = J^2 \subseteq JI \subseteq J \cap I = J$, hence $J = JI$. $\square$\\
\end{jpf}

We can then prove Theorem \ref{kH_3 not coherent thm}.\\

\begin{jpfof}[Theorem \ref{kH_3 not coherent thm}] \h
Consider the finitely-generated left ideal $I \subseteq kH_{3,t}$. We show that $I$ is not finitely-presented.\\
For a contradiction, suppose that $I$ is finitely-presented. Then we have a short exact sequence
\[ 0 \rightarrow K \rightarrow (kH_{3,t})^{2n+1} \rightarrow I \rightarrow 0 \]
with $K$ finitely-generated. We then apply the right exact functor $(kH_{3,t}/J) \otimes_{kH_{3,t}}$. Note that $(kH_{3,t}/J) \otimes_{kH_{3,t}} I = I/JI = I/J \cong I'$. Thus we have an exact sequence
\[  K/JK \rightarrow (kG_{3,t'})^{2n+1} \rightarrow I' \rightarrow 0. \]
Now $K$ is a finitely-generated $kH_{3,t}$-module, so $K/JK$ is a finitely-generated $kG_{3,t'}$-module, and hence $I'$ is a finitely-presented $kG_{3,t'}$-module.\\
But, $I' = kG_{3,t'} \otimes_{kG_{3,g}} I_g$, and $I_g$ is not finitely-presented (this follows from the proof of Theorem \ref{kG_3 not coherent thm}). Since $kG_{3,t'}$ is faithfully flat as a right ${kG_{3,g}}$-module, $I'$ cannot be finitely-presented, a contradiction.\\
So, the finitely-generated left ideal $I$ is not finitely-presented, and hence $kH_{3,t}$ is not left coherent. $\square$\\
\end{jpfof}

\subsection{The category of finitely-presented smooth mod $p$ representations} \label{The category of finitely-presented smooth mod $p$ representations}

Schraen, Hu, Vigneras and Wu have previously studied finitely-presented representations of $p$-adic Lie groups, particularly those of $GL_2(F)$. See \jcite{Schr15}, \jcite{Hu12}, \jcite{Vig11}, and \jcite{Wu21} respectively. Recently, Shotton has proved a link between coherence and finitely-presented smooth representations -- refer to Definition 1.1 and the proof of Theorem 4.5 in \jcite{Sho20}.

\begin{jthm} \h \label{coherence implies abelian thm}
Let $G$ be a $p$-adic Lie group. If $kG$ is a coherent ring, then the category of finitely-presented smooth mod $p$ representations of $G$ is an abelian category.\\
\end{jthm}

In Theorems \ref{kG_3 not coherent thm} and \ref{kH_3 not coherent thm}, we have proved that not all $p$-adic Lie groups have a coherent augmented Iwasawa algebra. It is natural to ask whether all $p$-adic Lie groups must have the category of smooth mod $p$ representations being abelian. We show in this subsection that this is also not the case -- again the counterexamples are the groups in Theorems \ref{kG_3 not coherent thm} and \ref{kH_3 not coherent thm}.\\
Recall, since $k$ is a (discrete) field of characteristic $p$, a smooth representation of $G$ is simply a representation where each vector is fixed by some compact open subgroup of $G$. We denote this category $\mathcal{C}_{k}(G)$, and the full subcategory of finitely-presented smooth representations $\mathcal{C}_{k}^{\mathrm{fp}}(G)$.\\
We give a class of examples of finitely-presented smooth mod $p$ representations, namely, compact induction of the trivial representation from a compact open subgroup.

\begin{jprop} \h \label{induction smooth prop}
Let $G$ be a $p$-adic Lie group, and $H \leq G$ be an open compact subgroup. The left $kG$-module $k\Big(\faktor{G}{H}\Big) \cong \faktor{kG}{\epsilon_G(H)}$ is a finitely-presented smooth (mod $p$) representation of $G$.\\
\end{jprop}
\begin{jpf} \h
The module $kG/\epsilon_G(H)$ is a module over $kG$ and hence a representation of $G$ where $G$ acts via its natural inclusion into $kG$. We also have that
\[ \faktor{kG}{\epsilon_G(H)} = kG \otimes_{kH} \faktor{kH}{\epsilon_H(H)} = kG \otimes_{kH} k = \bigoplus_{g \in G/H} g \otimes k, \]
where $k$ is the trivial module for $kH$. Now, for any $v \in k$ and $g \in G$, we have that for all $h \in H$,
\[ ghg^{-1} \cdot (g \otimes v) = g \otimes (h \cdot v) = g \otimes v, \]
so the compact open subgroup $gHg^{-1}$ fixes $g \otimes v$. Thus, for an arbitrary element
\[ x = \sum_{j=1}^n g_j \otimes v_j \in \faktor{kG}{\epsilon_G(H)}, \]
the compact open subgroup $K = \bigcap_{j=1}^n g_jHg_j^{-1}$ fixes $x$. Hence $kG/\epsilon_G(H)$ is a smooth representation of $G$.\\
Moreover, the trivial module $k$ is a finitely-presented left $kH$-module (since $kH$ is Noetherian), so has a presentation by finite-rank free $kH$-modules. Since $kG$ is a flat right $kH$-module, by Proposition \ref{locally profinite flatness prop}, it follows that $kG/\epsilon_G(H) = kG \otimes_{kH} k$ has a presentation by finite-rank free $kG$-modules, and so is finitely-presented as a $kG$-module.\\
By Proposition 3.8 of \jcite{Sho20}, it follows that $kG/\epsilon_G(H)$ is a finitely-presented smooth (mod $p$) representation of $G$. $\square$\\
\end{jpf}

\begin{jprop} \h \label{exact inclusion functor prop}
Let $G$ be a $p$-adic Lie group. The inclusion of $\mathcal{C}_{k}^{\mathrm{fp}}(G)$ into $\mathcal{C}_{k}(G)$ is exact. \\
\end{jprop}
\begin{jpf} \h
We will denote the inclusion functor $\mathcal{C}_{k}^{\mathrm{fp}}(G) \rightarrow \mathcal{C}_{k}(G)$ by $M \mapsto \underline{M}$.\\
Let $\phi: A \rightarrow B$ be a morphism in $\mathcal{C}_{k}^{\mathrm{fp}}(G)$, and suppose it has a kernel $j: K \rightarrow A$. We show that $\underline{j} : \underline{K} \rightarrow \underline{A}$ is the kernel of the morphism $\underline{\phi}$ in $\mathcal{C}_k(G)$.\\
For each $H \leq G$ a compact open subgroup, consider the functor 
\[ \mathrm{Hom}_{kG}(kG/\epsilon_G(H), \quad) : \mathcal{C}_k^{\mathrm{fp}}(G) \rightarrow kH\textrm{-Mod}. \]
These functors are left exact because $kG/\epsilon_G(H) \in \mathcal{C}_k^{\mathrm{fp}}(G)$ by Proposition \ref{induction smooth prop}. Moreover, for any $kG$-module $M$, there is a natural isomorphism
\[ \mathrm{Hom}_{kG}(kG/\epsilon_G(H), M) \cong M^H = \{x \in M \mid hx=x\textrm{ for all }h \in H\}. \]
Thus, the homomorphism $\underline{j}|_{\underline{K}^H} : \underline{K}^H \rightarrow \underline{A}^H$ is the kernel of $\underline{\phi}|_{\underline{A}^H} : \underline{A}^H \rightarrow \underline{B}^H$.\\
We now show that $\underline{j}$ satisfies the universal property of the kernel of $\underline{\phi}$. Let $K' \in \mathcal{C}_k(G)$ be a smooth representation and $j': K' \rightarrow \underline{A}$ be a $kG$-module homomorphism satisfying $\underline{\phi} \circ j' = 0$. Then for any compact open subgroup $H \leq G$, clearly $\underline{\phi}|_{\underline{A}^H} \circ j'|_{(K')^H} = 0$. So by the universal property for the kernel of $\underline{\phi}|_{\underline{A}^H}$, there is a unique morphism $u_H: (K')^H \rightarrow \underline{K}^H$ such that $\underline{j}|_{\underline{K}^H} \circ u_H = j'|_{(K')^H}$. If $H \leq H'$, then $u_H|_{(K')^{H'}} = u_{H'}$, by uniqueness. Moreover, $K'$ is smooth, meaning that
\[ K' = \bigcup_{\substack{\textrm{compact open}\\H \leq G}} (K')^H. \]
Therefore we can define the morphism of smooth representations,
\[ u: K' \rightarrow \underline{K}, \quad u(x) = u_H(x)\textrm{ if } x \in (K')^H, \]
and it follows that $\underline{j} \circ u = j'$ because this holds on each $(K')^H$ by definition of $u_H$. If $v$ is another morphism with this property, then $v|_{(K')^H} = u_H$ by uniqueness of $u_H$, so $v=u$. Thus $u$ is the unique such morphism. Therefore $j: \underline{K} \rightarrow \underline{A}$ satisfies the universal property of the kernel. Thus the inclusion functor $\mathcal{C}_{k}^{\mathrm{fp}}(G) \rightarrow \mathcal{C}_{k}(G)$ preserves kernels, so is left exact.\\
Furthermore, the cokernel of $\phi: A \rightarrow B$ in $\mathcal{C}_{k}^{\mathrm{fp}}(G)$ is $B/\phi(A)$ -- notice $B$ is finitely-presented and $\phi(A)$ is finitely-generated, so $B/\phi(A)$ is finitely-presented by Theorem 2.1.2 of \jcite{Glaz89}. Thus the inclusion functor preserves cokernels, so is right exact. $\square$ \\
\end{jpf}

Applying Proposition \ref{exact inclusion functor prop} and Lemma 1.6.2 of \jcite{Wei94} to this inclusion of categories, we have the following.

\begin{jcor} \h \label{closed under kernels cokernels cor}
The category $\mathcal{C}_k^{\mathrm{fp}}(G)$ of finitely-presented smooth representations is an abelian category if and only if it is closed under kernels and cokernels taken in the category $\mathcal{C}_k(G)$ of smooth representations, that is, $kG$-module kernels and cokernels.\\
\end{jcor} 

Notice that $\mathcal{C}_k^{\mathrm{fp}}(G)$ is always closed under cokernels taken in $\mathcal{C}_k(G)$, by Proposition 3.8 of \jcite{Sho20} and Theorem 2.1.2 of \jcite{Glaz89}. We now improve the result of Theorem \ref{coherence implies abelian thm} to a necessary and sufficient module-theoretic condition.

\begin{jprop} \h \label{abelian equivalent to coherent module prop}
Let $G$ be a $p$-adic Lie group. The category of finitely-presented smooth representations of $G$ is abelian if and only if ${kG/\epsilon_G(H)}$ is a coherent $kG$-module, for all compact open subgroups $H \leq G$.\\
\end{jprop}
\begin{jpf} \h
Suppose the category of finitely-presented smooth representations of $G$ is abelian. Now, for any open compact subgroup $H$, the $kG$-module ${kG/\epsilon_G(H)}$ is a smooth finitely-presented representation of $G$, by Proposition \ref{induction smooth prop}. Because ${kG/\epsilon_G(H)}$ is smooth, any $kG$-submodule or quotient is also a smooth representation. If $N \leq {kG/\epsilon_G(H)}$ is a finitely-generated $kG$-submodule, we have a short exact sequence of smooth representations
\[ 0 \rightarrow N \rightarrow {kG/\epsilon_G(H)} \rightarrow ({kG/\epsilon_G(H)})/N \rightarrow 0, \]
where the third and fourth terms are finitely-presented, by Proposition \ref{induction smooth prop} and Theorem 2.1.2 of \jcite{Glaz89}. Hence, $N$ must also be finitely-presented by Corollary \ref{closed under kernels cokernels cor}. Thus ${kG/\epsilon_G(H)}$ is a coherent $kG$-module. \vspace{2pt}\\
Conversely, suppose ${kG/\epsilon_G(H)}$ is a coherent $kG$-module for all open compact $H \leq G$. Let $M$ be a finitely-presented smooth representation of $G$. Then, $M$ is finitely-generated, so there is a surjection
\[ \phi: kG^n \rightarrow M \]
for some $n \in \mathbb{N}$. Because $M$ is smooth, for each $m \in M$, there exists an open compact subgroup $H \leq G$ such that $\epsilon_G(H) \cdot m =0$. By considering the vectors $m_i = \phi(e_i)$, where the $e_i$ are the standard basis of $kG^n$, we find that there is a surjection
\[ \tilde{\phi}: \faktor{kG}{\epsilon_G(H_1)} \oplus \dots \oplus \faktor{kG}{\epsilon_G(H_n)} \rightarrow M \]
for some open compact subgroups $H_1, \dots, H_n \leq G$. By assumption, ${kG/\epsilon_G(H_1) \oplus \dots \oplus kG/\epsilon_G(H_n)}$ is a coherent $kG$-module. Because $M$ is a finitely-presented quotient, it follows from Theorems 2.1.2 and 2.2.1 of \jcite{Glaz89} that $M$ is a coherent $kG$-module.\\
Therefore, the finitely-presented smooth representations of $G$ are exactly the smooth representations which are coherent $kG$-modules. The kernel of any homomorphism between coherent $kG$-modules is itself a coherent $kG$-module by Theorem 2.2.1 of \jcite{Glaz89}. Thus the kernel of any map between finitely-presented smooth representations must be a coherent, hence finitely-presented, representation. Therefore the category of finitely-presented smooth representations of $G$ is an abelian category, by Corollary \ref{closed under kernels cokernels cor}. $\square$\\
\end{jpf}

\begin{jlem} \h \label{coherent module extension lemma}
Let $G$ be a $p$-adic Lie group, $J \leq G$ be an open pro-$p$ subgroup. If ${kG/\epsilon_G(J)}$ is a coherent $kG$-module, then ${kG/\epsilon_G(H)}$ is a coherent $kG$-module for any compact open subgroup $H \leq G$.\\
\end{jlem}
\begin{jpf} \h
Let $H_1 \leq H_2$ be open compact subgroups of $G$. Suppose that ${kG/\epsilon_G(H_1)}$ is a coherent $kG$-module. We have a short exact sequence
\[ 0 \rightarrow \faktor{\epsilon_G(H_2)}{\epsilon_G(H_1)} \rightarrow \faktor{kG}{\epsilon_G(H_1)} \rightarrow \faktor{kG}{\epsilon_G(H_2)} \rightarrow 0, \]
and ${\epsilon_G(H_2)/\epsilon_G(H_1)}$ is finitely-generated because $H_2$ is compact. Thus ${\epsilon_G(H_2)/\epsilon_G(H_1)}$ is coherent, hence so is ${kG/\epsilon_G(H_2)}$.\vspace{2pt}\\
Now suppose ${kG/\epsilon_G(H_2)}$ is coherent, and that $H_1$ is normal and of $p$-power index in $H_2$. Denoting the trivial module for the group $\faktor{H_2}{H_1}$ by $k$, we have that
\[ \faktor{kG}{\epsilon_G(H_2)} = kG \otimes_{kH_2} k, \]
and
\[ \faktor{kG}{\epsilon_G(H_1)} = kG \otimes_{kH_2} k[H_2/H_1] .\]
Now, because $H_2/H_1$ is a finite $p$-group and $k$ is of characteristic $p$, the only irreducible $k[H_2/H_1]$-module is the trivial module $k$. Therefore the Jordan-Hölder series of $k[H_2/H_1]$ consists only of copies of $k$. That is, there exist $k[H_2/H_1]$-modules $M_1, \dots, M_n$ and short exact sequences
\begin{align*}
0 \rightarrow M_1 \rightarrow& k[H_2/H_1] \rightarrow k \rightarrow 0, \\
0 \rightarrow M_2 \rightarrow& M_1 \rightarrow k \rightarrow 0, \\
&\dots \\
0 \rightarrow M_n \rightarrow& M_{n-1} \rightarrow k \rightarrow 0,\\
0 \rightarrow k \rightarrow& M_n \rightarrow k \rightarrow 0.
\end{align*}
Now, $kG \otimes_{kH_2}$ is an exact functor by Theorem \ref{non-compact faithfully flat thm}. Applying it to the final sequence, we obtain
\[ 0 \rightarrow kG \otimes_{kH_2}k \rightarrow kG \otimes_{kH_2}M_n \rightarrow kG \otimes_{kH_2} k \rightarrow 0. \]
But $ kG \otimes_{kH_2}k$ is a coherent $kG$-module, and any extension of two coherent modules is coherent, by Theorem 2.2.1 of \jcite{Glaz89}. Thus $kG \otimes_{kH_2}M_n$ is coherent. Repeating this argument, we find that $kG \otimes_{kH_2}M_j$ is coherent for all $j$, and so \[ kG \otimes_{kH_2} k[H_2/H_1] = {kG/\epsilon_G(H_1)} \] is coherent. \vspace{2pt}\\
We now apply the above reasoning, supposing that ${kG/\epsilon_G(J)}$ is a coherent $kG$-module. Now, $J \cap H$ is an open subgroup of $J$, and so there exists an open normal subgroup $N \trianglelefteq J$ contained in $J \cap H$. Because $J$ is pro-$p$, $N$ has $p$-power index in $J$. Thus, because ${kG/\epsilon_G(J)}$ is coherent, ${kG/\epsilon_G(N)}$ is coherent. Then, $N \leq H$ and ${kG/\epsilon_G(N)}$ is coherent, therefore ${kG/\epsilon_G(H)}$ is coherent. $\square$\\
\end{jpf}

By combining Proposition \ref{abelian equivalent to coherent module prop} and Lemma \ref{coherent module extension lemma}, we can give a strengthening of Theorem \ref{coherence implies abelian thm}.

\begin{jcor} \h \label{fp smooth equivalences cor}
Let $G$ be a $p$-adic Lie group. The following are equivalent:
\begin{itemize}
\item The category of finitely-presented smooth representations of $G$ is abelian.
\item ${kG/\epsilon_G(H)}$ is a coherent $kG$-module, for all compact open subgroups ${H \leq G}$.
\item ${kG/\epsilon_G(J)}$ is a coherent $kG$-module, for some open pro-$p$ subgroup $J \leq G$.\\
\end{itemize}
\end{jcor}

Notice that if $kG$ is a coherent ring, ${kG/\epsilon_G(H)}$ is a finitely-presented quotient by Proposition \ref{induction smooth prop}, hence a coherent $kG$-module by Theorems 2.1.2 and 2.2.1 of \jcite{Glaz89}. Moreover, we can translate conditions on the modules $kG/\epsilon_G(H)$ into ideal-theoretic language.
\begin{jprop} \h
Let $G$ be a $p$-adic Lie group, $H \leq G$ be an open compact subgroup. The $kG$-module $kG/\epsilon_G(H)$ is coherent if and only if every finitely-generated ideal of $kG$ containing $\epsilon_G(H)$ is finitely-presented.\\
\end{jprop}
\begin{jpf} \h
The finitely-generated submodules of $kG/\epsilon_G(H)$ are exactly $I/\epsilon_G(H)$ for $I \subseteq kG$ a finitely-generated left ideal containing $\epsilon_G(H)$. By considering the short exact sequence
\[0 \rightarrow \epsilon_G(H) \rightarrow I \rightarrow I/\epsilon_G(H) \rightarrow 0,\]
by Corollary 6.25 of \jcite{Rot09} and Theorem 2.1.2 of \jcite{Glaz89}, $I$ is finitely-presented if and only if $I/\epsilon_G(H)$ is, since $\epsilon_G(H)$ is finitely-presented. The result follows. $\square$ \\
\end{jpf}

By Proposition \ref{abelian equivalent to coherent module prop}, it follows that if $kG$ has a finitely-generated ideal that contains the augmentation ideal of an open compact subgroup, but $I$ is not finitely-presented, then the category of finitely-presented smooth representations of $G$ is not abelian. This situation occurs with the two groups we have principally studied in this section.

\begin{jex} \hsmall
Recall the group $G_{3,g}$ and left ideal $I_g \subseteq kG_{3,g}$ from Definitions \ref{G_{3,g} def} and \ref{I_g def}. By Lemma \ref{explicit augmentation lemma}, $I_g$ contains the augmentation ideal of the compact open subgroup $\mathbb{U}_3'(\mathcal{O}_F)$. The proof of Theorem \ref{kG_3 not coherent thm} shows $I_g$ is finitely-generated but not finitely-presented. Thus the category of finitely-presented smooth representations of $G_{3,g}$ is not an abelian category.\\
\end{jex}

\begin{jex} \hsmall
Recall $H_{3,t}$ and $I$, as defined in subsection \ref{Proof of Theorem kH_3 not coherent thm}. By Lemma \ref{I contains J lemma}, $I$ contains the augmentation ideal of $\mathbb{U}_3(\mathcal{O}_F)$, and $I$ is finitely-generated but not finitely-presented. So the category of finitely-presented smooth representations of $H_{3,t}$ is not abelian.\\
\end{jex}

If $G$ is a closed subgroup of $G'$, the functor $kG' \otimes_{kG}$ takes finitely-presented $kG$-modules to finitely-presented $kG'$-modules. However it does not take smooth representations of $G$ to smooth representations of $G'$. Thus we cannot immediately deduce anything about the category of finitely-presented smooth representations of $G'$ from information about the corresponding category for $G$. For example, although $G_{3,g}$ is a closed subgroup of $GL_3(F)$, we cannot easily deduce any properties of smooth representations of $GL_3(F)$ from those of $G_{3,g}$.\\

If we make the stronger requirement that $G \leq G'$ be open, then any compact open subgroup of $G$ is compact open in $G'$, and ${kG' \otimes_{kG} (kG/\epsilon_G(H)) \cong kG'/\epsilon_{G'}(H)}$. So by Corollary \ref{fp smooth equivalences cor}, if the category of finitely-presented smooth mod $p$ representations of $G$ is not abelian, neither is this category for $G'$. For example, we can apply this to $G = G_{3,g}$ and $G' = G_{3,t'}$ or $G_{3,D,E}$.\\

\section{Characterisations of coherence} \label{Characterisations of coherence}

In this section we prove Theorem \ref{Main result general coherence thm}.

\subsection{The class of solvable groups} \label{The class of solvable groups}

We reiterate the notation and setup of Section \ref{Main results}.\\

Let $n \in \mathbb{N}$. Let $\mathbb{U}_n$ be the affine group scheme of upper unitriangular matrices in $\mathbb{GL}_n$ over $F$. Let $\mathbb{U}$ be a closed affine group subscheme of $\mathbb{U}_n$, defined and split over $F$. Let $U = \mathbb{U}(F)$.\\
Let $T$ be a closed subgroup -- in the $p$-adic topology -- of the diagonal elements $\mathbb{D}_n(F) \leq GL_n(F)$, such that $T$ normalises $U$. Let $G = \langle T, U \rangle \cong T \ltimes U$, so $G$ is a $p$-adic Lie group.\\ We define a set of roots of $U$ with respect to $T$, similarly to section 8.17 of \jcite{Bo91}.\vspace{2pt}\\
Let $\mathfrak{u} = \textrm{Lie}(U)$. Since $T$ acts on $U$ via conjugation, $T$ acts on $\mathfrak{u}$. Let
\[X(T) = \{ \beta: T \rightarrow F^{\times} \mid \beta \textrm{ a continuous group homomorphism} \} \] 
be the character group of $T$. For each $\beta \in X(T)$, define the subspace, called a weight space,
\[ \mathfrak{u}_\beta = \{ v \in \mathfrak{u} \mid t \cdot v = \beta(t)v \textrm{ }\forall t \in T \} \leq \mathfrak{u}. \]
The $\mathfrak{u}_\beta$ are finite-dimensional subspaces of $\mathfrak{u}$, and may have dimension larger than $1$.
We define the set of weights
\[ \Phi = \{ \beta \in X(T) \mid \mathfrak{u}_\beta \neq 0 \}. \]
Note that we allow the trivial character to be an element of $\Phi$.\\
Let $\mathbb{T}_{\mathbb{U}}$ be the normaliser of $\mathbb{U}$ in $\mathbb{D}_n$, so $T$ is a subgroup of $\mathbb{T}_{\mathbb{U}}(F)$. By Proposition 8.4 of \jcite{Bo91}, $\mathbb{T}_{\mathbb{U}}(F)$ acts diagonalisably on $\mathfrak{u}$, and hence so does $T$. \\
Thus $\mathfrak{u}$ is the direct sum of its weight spaces:
\[ \mathfrak{u} = \bigoplus_{\beta \in \Phi} \mathfrak{u}_\beta .\]
The $\mathfrak{u}_\beta$ correspond to root subgroups $U_\beta \leq U$.\\
For $\beta_1, \beta_2 \in X(T)$, $v_1 \in \mathfrak{u}_{\beta_1}, v_2 \in \mathfrak{u}_{\beta_2}$, we have that $t \cdot [v_1, v_2] = [t \cdot v_1, t \cdot v_2] = \beta_1(t)\beta_2(t) [v_1, v_2]$ for all $t \in T$. We use additive notation $\beta_1 + \beta_2$ to denote the character $t \mapsto \beta_1(t)\beta_2(t)$. Thus, we have that
\[ [\mathfrak{u}_{\beta_1}, \mathfrak{u}_{\beta_2} ] \leq \mathfrak{u}_{\beta_1+ \beta_2}. \]
We define the group homomorphism $f$ by
\[ f: T \rightarrow \mathbb{Z}^\Phi, \quad f(t) = \big(v_F( \beta(t) )\big)_{\beta \in \Phi}, \]
where $v_F : F^\times \rightarrow \mathbb{Z}$ is the discrete valuation on $F$. For each $\beta$, we will write $n_\beta(t) = v_F(\beta(t))$, so that $f(t) = (n_\beta(t))_{\beta \in \Phi}$. \vspace{2pt}
Note that if $\beta_1, \beta_2 \in \Phi$ are such that $\beta_1 + \beta_2 \in \Phi$, then for all $t \in T$,
\[n_{\beta_1 + \beta_2}(t) = v_F(\beta_1(t)\beta_2(t)) = v_F(\beta_1(t)) + v_F(\beta_2(t)) = n_{\beta_1}(t) + n_{\beta_2}(t). \]
In this section, we prove Theorem \ref{Main result general coherence thm}, which states that the augmented Iwasawa algebra $kG$ is coherent if and only if the image $f(T)$ is cyclic and generated by an element of $(\mathbb{Z}_{\geq 0})^\Phi$.

\subsection{Non-coherence: The case of cyclic torus}

In this subsection, assume that $T$ is cyclic and choose a fixed generator $t_0 \in T$. For convenience, we will write $n_\beta = n_\beta(t_0)$ in this subsection.\\

The next two lemmas are basic Lie algebra facts we need for examining the structure of $\mathfrak{u} = \textrm{Lie}(U)$.

\begin{jlem} \h \label{central series lemma}
Let $\mathfrak{g}$ be a Lie algebra over a field $K$, and
\[ \mathfrak{g} = \mathfrak{g}^1 \geq \mathfrak{g}^2 \geq \mathfrak{g}^3 \geq \dots \]
be its lower central series. Let $S$ be a generating set for $\mathfrak{g}$. Then, for any $m \in \mathbb{N}$, the quotient Lie algebra $\mathfrak{g}^m/\mathfrak{g}^{m+1}$ is generated by the (image of the) $m$-fold commutators
\[ S^m = \{ [s_1, [s_2, \dots, s_m]] \mid s_1, s_2, \dots, s_m \in S \} .\]
\end{jlem}
\begin{jpf} \h
Let $\mathfrak{h}^m$ be the Lie subalgebra of $\mathfrak{g}$ generated by $S^m$. We proceed by induction.\\
The case $m=1$ is clear. If $m \geq 1$, then
\[ \mathfrak{g}^{m+1} = [\mathfrak{g}, \mathfrak{g}^m] = [\mathfrak{g}, \mathfrak{h}^m + \mathfrak{g}^{m+1}] = [\mathfrak{g}, \mathfrak{h}^m] + \mathfrak{g}^{m+2}, \]
by induction. Then $[\mathfrak{g}, \mathfrak{h}^m]$ is the Lie subalgebra generated by $\cup_{r \geq m+1} S^k$. If $r > m+1$, then $S^r \subseteq \mathfrak{g}^{m+2}$, and so
\[ [\mathfrak{g}, \mathfrak{h}^m] + \mathfrak{g}^{m+2} = \mathfrak{h}^{m+1} +  \mathfrak{g}^{m+2}, \]
therefore $\mathfrak{g}^{m+1} / \mathfrak{g}^{m+2}$ is generated by $\mathfrak{h}^m$, hence by $S^m$, as required. $\square$\\
\end{jpf}

\begin{jlem} \h \label{nilpotent 3 element lemma}
Let $\mathfrak{g}$ be a non-abelian nilpotent Lie algebra over a field $K$, generated by two distinct 1-dimensional subspaces $\mathfrak{h}_1, \mathfrak{h}_2$. Let $\mathfrak{z} \leq \mathfrak{g}^2$ be a 1-dimensional subspace. Let $\mathfrak{g}'$ be the subalgebra generated by $\mathfrak{h}_1, \mathfrak{z}$. Then $\mathfrak{g}'$ is a proper subalgebra of $\mathfrak{g}$. \\
\end{jlem}

\begin{jpf} \h
Let $x \in \mathfrak{h}_1, y \in \mathfrak{h}_2, z  \in \mathfrak{z}$ be non-zero. Suppose that $\mathfrak{g} = \mathfrak{g}'$.\\
Then $\mathfrak{g} =  \mathfrak{g}'$ is generated by $x,z$, so $\mathfrak{g}^2 / \mathfrak{g}^3$ is generated by $[x,z]$, by Lemma \ref{central series lemma}. But $z \in \mathfrak{g}^2$, so $[x,z] \in \mathfrak{g}^3$, hence $\mathfrak{g}^2 = \mathfrak{g}^3$. Because $\mathfrak{g}$ is nilpotent, $\mathfrak{g}^2=0$, hence $z=0$, a contradiction. So $\mathfrak{g}' \neq \mathfrak{g}$. $\square$\\
\end{jpf}

Recall the $p$-adic Lie groups $G_{3,t'}, H_{3,t}$ from Theorems \ref{kG_3 not coherent thm} and \ref{kH_3 not coherent thm}. The following proposition characterises when these groups appear as closed subgroups of $G$. 

\begin{jprop} \h \label{1 variable non-coherence prop}
Suppose there exists $\alpha, \beta \in \Phi$ such that $n_{\alpha} > 0, n_{\beta} < 0$. Then $G$ contains a closed subgroup of the form $G_{3,t'}$ or $H_{3,t}$.\\
\end{jprop}
\begin{jpf} \h
Let $\mathfrak{u}_{\alpha}' \leq\mathfrak{u}_{\alpha}$ and $\mathfrak{u}_{\beta}' \leq \mathfrak{u}_{\beta}$ be 1-dimensional subspaces. The subalgebra $\mathfrak{v}$ generated by these subspaces corresponds to a closed subgroup $U' \leq U$ that is normalised by $T = \langle t_0 \rangle$, since $T$ fixes $\mathfrak{u}_{\alpha}$ and $\mathfrak{u}_{\beta}$. So, it is enough to prove that $U'$ contains a closed subgroup of the form $G_{3,t'}$ or $H_{3,t}$. That is, we may without loss of generality assume that $\mathfrak{u}$ is generated by $\mathfrak{u}_{\alpha}, \mathfrak{u}_{\beta}$, and that these spaces are 1-dimensional.\\
We proceed by induction on $\dim U = \dim \mathfrak{u}$. Clearly $\dim U \geq 2$, and if $\dim U = 2$, then $\mathfrak{u}_{\alpha}, \mathfrak{u}_{\beta}$ must commute, implying $G$ is of the form $G_{3,t'}$. \\
If $\dim U > 2$, notice that if we can find $\alpha_0, \beta_0 \in \Phi$ with $n_{\alpha_0} > 0, n_{\beta_0} < 0$, and subspaces $\mathfrak{u}_{\alpha_0}' \leq \mathfrak{u}_{\alpha_0}$, $\mathfrak{u}_{\beta_0}' \leq \mathfrak{u}_{\beta_0}$ such that $\mathfrak{u}_{\alpha_0}' , \mathfrak{u}_{\beta_0}' $ generate a proper subalgebra $\mathfrak{v} \leq \mathfrak{u}$, we are done by induction. This follows because $\mathfrak{v}$ corresponds to a proper closed subgroup $V \leq U$ with $V$ normalised by $T$, and satisfying the hypothesis of the proposition. \vspace{2pt}\\
Let $\dim U>2$. Then $\mathfrak{u}_{\alpha}, \mathfrak{u}_{\beta}$ do not commute. Let $\delta = \alpha+\beta$, and
\[\mathfrak{u}_{\delta}' = [\mathfrak{u}_{\alpha}, \mathfrak{u}_{\beta}]  \leq \mathfrak{u}_{\delta}. \]
If $\mathfrak{u}_{\delta}'$ commutes with both $\mathfrak{u}_{\alpha}$ and $\mathfrak{u}_{\beta}$, then $G$ must be of the form $H_{3,t}$.\vspace{2pt}\\
If not, suppose $n_{\alpha} + n_{\beta} = n_{\delta} > 0$. Let $\mathfrak{v}$ be the Lie subalgebra of $\mathfrak{u}$ generated by $\mathfrak{u}_{\delta}', \mathfrak{u}_{\beta}$. As $\mathfrak{u}_{\delta}' \leq \mathfrak{u}^2$, by Lemma \ref{nilpotent 3 element lemma}, $\mathfrak{v}$ is a proper subalgebra and $n_{\delta}  > 0, n_{\beta} < 0$, so we are done by induction. If $n_{\alpha} + n_{\beta} = n_{\delta} < 0$, the argument is very similar. \vspace{2pt}\\
Suppose that $n_{\alpha} + n_{\beta} = n_{\delta}=0$, and suppose that $[\mathfrak{u}_{\alpha}, \mathfrak{u}_{\delta}'] \neq 0$. Let $\gamma = \alpha + \delta$, and $[\mathfrak{u}_{\alpha}, \mathfrak{u}_{\delta}'] = \mathfrak{u}_{\gamma}' \leq\mathfrak{u}_{\gamma}$. Then $n_{\gamma} = n_{\alpha} + n_{\delta} = n_{\alpha} > 0$. Let $\mathfrak{v}$ be the subalgebra generated by $\mathfrak{u}_{\beta}, \mathfrak{u}_{\gamma}'$. Then $\mathfrak{u}_{\gamma}' \leq \mathfrak{u}^3 \leq \mathfrak{u}^2$, so by Lemma \ref{nilpotent 3 element lemma}, the subalgebra $\mathfrak{v}$ is proper. Also $n_{\gamma} > 0, n_{\beta} < 0$, so we are done by induction. If $[\mathfrak{u}_{\alpha}, \mathfrak{u}_{\delta}'] = 0$, then $[\mathfrak{u}_{\beta}, \mathfrak{u}_{\delta}'] \neq 0$ and the argument is very similar. $\square$\\
\end{jpf}

Since $kG_{3,t'}, kH_{3,t}$ are not coherent, by Theorems \ref{kG_3 not coherent thm} and \ref{kH_3 not coherent thm}, if the hypothesis of the above proposition is satisfied, then $kG$ is not coherent, by Proposition \ref{closed subgroup coherence prop}.

\begin{jcor} \h \label{cyclic torus non-coherence cor}
Suppose $f(t_0) \not \in (\mathbb{Z}_{\geq 0})^\Phi \cup (\mathbb{Z}_{\leq 0})^\Phi$, and let $G = \langle T, U \rangle = \langle t_0, U \rangle$. Then the augmented Iwasawa algebra $kG$ is not coherent.\\
\end{jcor}

\subsection{Non-coherence for solvable groups}

Let $G$ be as in subsection \ref{The class of solvable groups}. We show that the class of groups $G$ for which $kG$ is coherent is quite small.

\begin{jthm} \h \label{general T non-coherence thm}
The augmented Iwasawa algebra $kG$ is coherent only if $f(T) \leq \mathbb{Z}^{\Phi}$ is cyclic and generated by an element in ${(\mathbb{Z}_{\geq 0})}^{\Phi}$.\\
\end{jthm}

The proof of this theorem is a deduction from Corollary \ref{cyclic torus non-coherence cor}, using some information about subgroups of $\mathbb{Z}^N$.

\begin{jlem} \h \label{Z^N lemma}
Let $N \in \mathbb{N}$. Let $x, y \in {(\mathbb{Z}_{\geq 0})}^N$. Suppose that the submodule
\[ A = \mathbb{Z}x + \mathbb{Z}y \leq \mathbb{Z}^N \]
is contained in ${(\mathbb{Z}_{\geq 0})}^N \cup {(\mathbb{Z}_{\leq 0})}^N$. Then $A$ is generated by one element.\\
\end{jlem}
\begin{jpf} \h
Note that if $x=0$ or $y=0$ we are done so assume $x,y \neq 0$. If there exist $i, j \in \{1,2, \dots, N\}$ such that $x_i > y_i$ but $x_j < y_j$, then $x-y \not \in {(\mathbb{Z}_{\geq 0})}^N \cup {(\mathbb{Z}_{\leq 0})}^N$, a contradiction. So, without loss of generality, assume that $x_i \geq y_i$ for all $i \in \{1,2, \dots, N\}$.\\
If there exists $i \in \{1,2, \dots, N\}$ such that $x_i=y_i=0$, then we may consider $A$ as a submodule of $\mathbb{Z}^{N-1}$ by removing the $i$th component. Moreover, if there exists $i \in \{1,2, \dots, N\}$ such that $y_i=0$ and $x_i > 0$, let $j \in \{1,2, \dots, N\}$ be such that $y_j > 0$. Then, there exists $m \in \mathbb{N}$ large enough such that $x_j -  my_j < 0$. But $x_i - my_i = x_i >0$, thus $x-my \not \in {(\mathbb{Z}_{\geq 0})}^N \cup {(\mathbb{Z}_{\leq 0})}^N$, a contradiction. \vspace{4pt}\\
Therefore, we may assume that $x_i \geq y_i >0$ for all $i \in \{1,2, \dots N \}$.\\
For each $i$, let $q_i = \frac{x_i}{y_i} \in \mathbb{Q}$. Let $q = \min\{q_i \mid i \in \{1,2, \dots N \} \}$, and by reordering, assume that $q=q_1$ without loss of generality. Let $q = \frac{a}{b}$ for coprime $a,b \in \mathbb{N}$, and $z = bx-ay \in A$. Suppose $z \neq 0$. Then $z_1=0$ -- but $y_1 > 0$, so we have a contradiction as seen above. Thus, $z=0$. \vspace{2pt}\\
So, $bx=ay$. By Euclid's algorithm, $\exists r, s \in \mathbb{Z}$ such that $ra-sb=1$, since $a,b$ are coprime. Then
\begin{align*} x &= (ra-sb)x = rax-sbx =rax - say = a(rx-sy),\\ y &= (ra-sb)y = ray - sby = rbx - sby = b(rx-sy), \end{align*}
and so $A$ is generated by $rx-sy$. $\square$\\
\end{jpf}

\begin{jpfof}[Theorem \ref{general T non-coherence thm}] \h
Suppose $kG$ is coherent. Note that for any $t \in T$, $\langle t, U \rangle$ is a closed subgroup of $G$, so by Proposition \ref{closed subgroup coherence prop}, $k\langle t, U \rangle$ is coherent. Thus, by Corollary \ref{cyclic torus non-coherence cor}, $f(t) \in {(\mathbb{Z}_{\geq 0})}^{\Phi} \cup {(\mathbb{Z}_{\leq 0})}^{\Phi}$ for all $t \in T$. That is,
\[ f(T) \subseteq {(\mathbb{Z}_{\geq 0})}^{\Phi} \cup {(\mathbb{Z}_{\leq 0})}^{\Phi} .\]
Now, $f(T)$ is finitely-generated, thus by Lemma \ref{Z^N lemma} and induction, $f(T)$ is a cyclic subgroup of $\mathbb{Z}^\Phi$, and must be generated by an element of ${(\mathbb{Z}_{\geq 0})}^{\Phi}$.   $\square$\\
\end{jpfof}

\subsection{Coherence for solvable groups}

In this subsection we show that the converse to Theorem \ref{general T non-coherence thm} also holds.\\
First, we examine the structure of the subgroups of $(F^\times)^n$.

\begin{jlem} \h \label{closed subgroup torus lemma}
Let $T \leq (F^\times)^n$ be a closed subgroup. Then,
\[ T \cong \mathbb{Z}^d \times \large( T\cap (\mathcal{O}_F^\times)^n \large) \]
as $p$-adic Lie groups, and $T \cap (\mathcal{O}_F^\times)^n$ is compact.\\
\end{jlem}
\begin{jpf} \h
Note that $ (F^\times)^n = (\pi^{\mathbb{Z}} \times \mathcal{O}_F^\times)^n$, where $\pi \in F$ is a uniformiser. So let $q : T \rightarrow  (\pi^{\mathbb{Z}})^n$ be the natural projection restricted to $T$. Then $\textrm{Im }q$ is (isomorphic to) a subgroup of $\mathbb{Z}^n$, hence is isomorphic to $\mathbb{Z}^d$ for some $d \leq n$. Then, $\Ker q = T \cap  (\mathcal{O}_F^\times)^n $. Thus we have a short exact sequence of abelian groups
\[ 0 \rightarrow T \cap  (\mathcal{O}_F^\times)^n \rightarrow T \rightarrow \mathbb{Z}^d \rightarrow 0, \]
which must split because $ \mathbb{Z}^d$ is a free $\mathbb{Z}$-module. So $T \cong \mathbb{Z}^d \times \large( T\cap (\mathcal{O}_F^\times)^n \large)$, and $T\cap (\mathcal{O}_F^\times)^n$ is a closed subgroup of 
the compact group $ (\mathcal{O}_F^\times)^n$, hence $T\cap (\mathcal{O}_F^\times)^n$ is compact. $\square$\\
\end{jpf}

\begin{jthm} \h \label{general T coherence thm}
If $f(T) \leq \mathbb{Z}^\Phi$ is cyclic and generated by an element of $(\mathbb{Z}_{\geq 0})^\Phi$, then the augmented Iwasawa algebra $kG$ is coherent. \\
\end{jthm}
\begin{jpf} \h
By Lemma \ref{closed subgroup torus lemma}, $T = T_0 \times T_1$ where $T_0 \cong \mathbb{Z}^d$ and $T_1$ is a compact $p$-adic Lie group, with $T_1$ a subgroup of the group of diagonal matrices in $GL_n(\mathcal{O}_F)$. Because $U$ is a subgroup of the upper unitriangular matrices in $GL_n(F)$, we can show easily that $T_1 \leq \Ker f$. Thus $f(T_0)=f(T)$ is cyclic. If $f(T_0)$ is non-trivial, let $\{t_1, \dots, t_{d-1}\} \subseteq T_0$ generate $\mathrm{Ker}(f|_{T_0})$ and $t_0 \in T_0$ be such that $f(t_0)$ generates $f(T_0)$. If $f(T_0)$ is trivial, let $\{ t_0, t_1, \dots, t_{d-1} \}$ be any set of generators for $T_0$. Then $T_0 = \langle t_0 \rangle \times \langle t_1, \dots, t_{d-1} \rangle \cong \mathbb{Z} \times \mathbb{Z}^{d-1}$, and $\langle t_1, \dots, t_{d-1} \rangle \leq \Ker f$. \vspace{2pt}\\
Next, choose $\mathcal{O}_F$-Lie subalgebras of $\mathfrak{u}$,
\[ \mathfrak{v}_0 \leq \mathfrak{v}_1 \leq \mathfrak{v}_2 \leq \dots \]
of the form \[\mathfrak{v}_m = \bigoplus_{\beta \in \Phi} \mathfrak{v}_{\beta, m} \]
where $\mathfrak{v}_{\beta, m} \leq \mathfrak{u}_{\beta}$ is an $\mathcal{O}_F$-lattice, and such that $\mathfrak{u} = \varinjlim_{m \geq 0} \mathfrak{v}_m$.\\
These correspond to compact open subgroups $V_0 \leq V_1 \leq V_2 \leq \dots$ of $U$, and we have that $V_{m-1}$ is an open (hence closed) subgroup of $V_m$ for all $m \geq 1$.\\ 
Then,
\[ t_0 \cdot \mathfrak{v}_{m} = \bigoplus_{\beta \in \Phi} t_0 \cdot \mathfrak{v}_{\beta, m} = \bigoplus_{\beta \in \Phi} {\pi}^{n_{\beta}(t_0)}\mathfrak{v}_{\beta, m}, \] 
where $\pi \in F$ is the uniformiser of $F$. Thus $t_0 \cdot \mathfrak{v}_{m}$ has finite index in $\mathfrak{v}_{m}$, and so $t_0V_mt_0^{-1}$ is a finite index subgroup of $V_m$. Thus $t_0V_mt_0^{-1} \leq V_m$ is an open subgroup (since it is also closed). Moreover, $t_0 \cdot \mathfrak{v}_{m}$ is an $\mathcal{O}_F$-Lie ideal of $\mathfrak{v}_{m}$, thus $t_0V_mt_0^{-1}$ is an open normal subgroup of $V_m$. \vspace{2pt}\\
Let $V_m' = \langle T_1, V_m \rangle = T_1 \ltimes V_m$ for each $m \geq 0$. Let $\mathfrak{v}_m'$ be the corresponding Lie algebra. Because $t_0 \cdot \mathfrak{v}_m$ is a finite index Lie ideal of $\mathfrak{v}_m$, it follows that $t_0 \cdot \mathfrak{v}_m'$ is a finite index Lie ideal of $\mathfrak{v}_m'$, for all $m$. Thus $t_0V_m't_0^{-1}$ is a finite index, closed, normal subgroup of $V_m'$, and hence is also an open normal subgroup. By Corollary 19.4 of \jcite{Schn11}, it follows that the Iwasawa algebra $kV_m'$ is a free right $k(t_0V_m't_0^{-1})$-module.\vspace{2pt}\\
Let $V_m'' = \langle V_m', t_1, \dots, t_{d-1} \rangle$ for each $m \geq 0$. Now, each of the $t_1, \dots t_{d-1}$ commute with $T_1$, and act, via conjugation, as an automorphism on $V_m$. Thus the $t_1, \dots, t_{d-1}$ act as automorphisms of $V_m'$, and therefore as ring automorphisms of the Iwasawa algebra $kV_m'$. So, the augmented Iwasawa algebra of $V_m''$ is a skew Laurent polynomial ring, $kV_m'' = kV_m'[t_1, t_1^{-1}, \dots, t_{d-1}, t_{d-1}^{-1}]$. \vspace{2pt}\\
By the non-commutative Hilbert basis theorem, $kV_m''$ is a Noetherian ring for each $m$. Let
\[ \sigma_{t_0}: kV_m'' \rightarrow kV_m'', \quad \sigma_{t_0}(x) = t_0xt_0^{-1} \]
be the ring endomorphism given by conjugation by $t_0$. Then $\sigma_{t_0}$ is injective. Note $\sigma_{t_0}$ fixes the $t_1, \dots, t_{d-1}$. The image of $\sigma_{t_0}$ is
\[ \sigma_{t_0}(kV_m'') =  k(t_0V_m''t_0^{-1}) =  k(t_0V_m't_0^{-1})[t_1, t_1^{-1}, \dots, t_{d-1}, t_{d-1}^{-1}]. \] Thus $kV_m''$ is a free, hence flat, right $\sigma_{t_0}(kV_m'')$-module. \\
By Emerton's result on the coherence of skew polynomial rings, Theorem \ref{coherent 1 variable thm}, it follows that the skew polynomial ring $kV_m''[t_0] = kV_m''[t_0; \sigma_{t_0}]$ is a left coherent ring, for all $m \geq 0$. Now, since $V_m$ is a closed subgroup of $V_{m+1}$, it follows that $V_m''$ is a closed subgroup of $V_{m+1}''$, therefore $kV_{m+1}''$ is a flat right $kV_m''$-module, by Proposition \ref{locally profinite flatness prop}. Hence $kV_{m+1}''[t_0]$ is a flat right $kV_m''[t_0]$-module, for all $m \geq 0$. Thus, by Lemma \ref{direct limit coherence lemma}, the direct limit
\[ \varinjlim_{m \geq 0} kV_m''[t_0] = \left(\varinjlim_{m \geq 0} kV_m'' \right) [t_0] = kG'  [t_0] \]
is a left coherent ring, where $G' = \langle  t_1, \dots t_{d-1},T_1, U \rangle$. Then, the set $\{1, t_0, t_0^2, \dots \}$ is a left denominator set in $kG'  [t_0]$, and $kG = kG'[t_0, t_0^{-1}]= \{1, t_0, t_0^2, \dots \}^{-1} kG'[t_0]$. 
Therefore, by Corollary \ref{localisations are coherent rings cor}, $kG$ is a left coherent ring. $\square$\\
\end{jpf}

We can conclude from Theorem \ref{general T non-coherence thm} and Theorem \ref{general T coherence thm} a necessary and sufficient criterion for $kG$ to be coherent -- this is Theorem \ref{Main result general coherence thm}.

\begin{jcorr}[Theorem \ref{Main result general coherence thm}] \h
The augmented Iwasawa algebra $kG$ is coherent if and only if $f(T) \leq \mathbb{Z}^\Phi$ is a cyclic subgroup generated by an element of $(\mathbb{Z}_{\geq 0})^\Phi$.\\
\end{jcorr}

\section{Coherence for algebraic groups} \label{Coherence for algebraic groups}

In this section we prove Corollary \ref{Main result solvable algebraic group non-coherence cor}, deducing Theorem \ref{Main result semisimple algebraic group non-coherence thm} and Corollary \ref{Main result general linear group coherence cor}.

\subsection{Coherence for solvable and split-semisimple groups}

\begin{jprop} \h \label{solvable algebraic group non-coherence prop}
Let $\mathbb{U}$ be an affine group subscheme of the upper unitriangular matrices $\mathbb{U}_n$, defined and split over $F$. Let $\mathbb{T}$ be a split affine group subscheme of the diagonal group $\mathbb{D}_n \leq \mathbb{GL}_n$ such that $\mathbb{T}$ normalises $\mathbb{U}$, and let $\mathbb{G} = \mathbb{T} \ltimes \mathbb{U} \leq \mathbb{GL}_n$. Let $G = \mathbb{G}(F) = T \ltimes U$, where $T = \mathbb{T}(F)$, $U = \mathbb{U}(F)$.\\
If the rank of the root system of $U$ with respect to $T$ is at least two, then $kG$ is not a coherent ring.\\
\end{jprop}
\begin{jpf} \h
The character group of $\mathbb{T}$ is $X(\mathbb{T}) = \{\textrm{morphisms of group schemes }\mathbb{T} \rightarrow \mathbb{G}_m \}$. Let $\mathbb{T} \cong \mathbb{G}_m^d$. Then $X(\mathbb{T}) \cong \mathbb{Z}^d$, and $T \cong (F^{\times})^d$, so
\[X(T) = \{ \alpha_{(j_1, \dots, j_d)} \mid (j_1, \dots, j_d) \in \mathbb{Z}^d \} \cong \mathbb{Z}^d,\]
where
\[\alpha_{(j_1, \dots, j_d)}: T \rightarrow {F}^\times, \quad \alpha_{(j_1, \dots, j_d)} (t_1, \dots, t_d) = t_1^{j_1} \dots t_d^{j_d}. \]
We have that $\Phi  = \{ \beta \in X(T) \mid u_\beta \neq 0 \}$. For each $\beta \in \Phi$, let $(j_1^\beta, \dots , j_d^\beta) \in \mathbb{Z}^d$ such that $\beta =  \alpha_{(j_1^\beta, \dots, j_d^\beta)}$.\vspace{2pt}\\
Let $t = (\lambda_1, \dots, \lambda_d) \in T$, and let $n_i = v_F(\lambda_i)$ for each $i \in \{1,2, \dots, d\}$. Then
\[ v_F(\beta(t)) = v_F(\lambda_1^{j_1^\beta}\dots \lambda_d^{j_d^\beta}) = j_1^\beta n_1 + \dots + j_d^\beta n_d. \]
Let $M$ be the $(|\Phi| \times d)$-matrix with coefficients $M_{\beta i} = j_i^\beta$. Recall from subsection \ref{The class of solvable groups} the group homomorphism $f: T \rightarrow \mathbb{Z}^\Phi$. Then, $f$ is given by
\[ f(t) = M \cdot \begin{pmatrix} n_1 \\ \dots \\ n_d \end{pmatrix}, \]
where $n_i = v_F(\lambda_i)$. So, the image $f(T) \leq \mathbb{Z}^\Phi$ is given by the $\mathbb{Z}$-span of the columns of $M$.\\
If $kG$ is coherent, then by Theorem \ref{Main result general coherence thm}, $f(T)$ is cyclic, so $\mathbb{Q} \otimes_\mathbb{Z} f(T) \leq \mathbb{Q}^{\Phi}$ is either trivial or a 1-dimensional vector space. Hence the $\mathbb{Q}$-rank of $M$ is at most 1, and so the vectors $\{(j_1^\beta, \dots, j_d^\beta) \mid \beta \in \Phi\}$ are contained in a 1-dimensional vector subspace of $\mathbb{Q}^d$. So, the rank of the root system of $U$ with respect to $T$ is at most 1.\vspace{2pt}\\
Hence, if the rank of the root system of $U$ with respect to $T$ is 2 or greater, $kG$ is not coherent. $\square$ \\
\end{jpf}

We then deduce Corollary \ref{Main result solvable algebraic group non-coherence cor}.

\notinsubfile{\solvablealgebraicgroups*}
\onlyinsubfile{
\begin{jcor} \h 
Let $\mathbb{G}$ be a finite-dimensional split-solvable affine group scheme defined over $F$, and $G=\mathbb{G}(F)$. If the root system of $\mathbb{G}$ has rank 2 or greater, then $kG$ is not coherent.\\
\end{jcor}
}

This follows directly from Proposition \ref{solvable algebraic group non-coherence prop} because the class of groups considered is the same.

We can then show that nearly all split-semisimple groups over $F$ have a non-coherent augmented Iwasawa algebra.

\begin{jcor} \h \label{semisimple algebraic group non-coherence cor}
Let $\mathbb{G}$ be a split-semisimple affine group scheme defined over $F$, with root system of rank at least 2. Let $G = \mathbb{G}(F)$. Then the augmented Iwasawa algebra $kG$ is not coherent.\\
\end{jcor}
\begin{jpf} \h
This follows by considering the Borel subgroup $\mathbb{B} \leq \mathbb{G}$, with $B = \mathbb{B}(F) \leq G$. We have that $\mathbb{B}$ is a split-solvable affine group scheme with a root system of rank at least 2, thus $\mathbb{B}$ is as in the statement of Corollary \ref{Main result solvable algebraic group non-coherence cor}, and hence $kB$ is not coherent. But $B$ is a closed subgroup of $G$, so by Proposition \ref{closed subgroup coherence prop}, $kG$ is not coherent. $\square$\\ 
\end{jpf}

\notinsubfile{\semisimplealgebraicgroups*}
\onlyinsubfile{
\begin{jcor} \h
Let $\mathbb{G}$ be a split-semisimple affine group scheme defined over $F$. Let $G = \mathbb{G}(F)$. Then the augmented Iwasawa algebra $kG$ is coherent if and only if the rank of the root system of $\mathbb{G}$ is $1$.\\
\end{jcor}
}
\begin{jpf} \h
The ``only if'' part of the statement follows from Corollary \ref{semisimple algebraic group non-coherence cor}.\\
If the root system of $\mathbb{G}$ has rank 1, then $\mathbb{G} = \mathbb{SL}_2$ or $\mathbb{PGL}_2$, by the classification of split-semisimple linear algebraic groups, see for example Table 9.2 of \jcite{MaTe11}.\\
If $\mathbb{G} = \mathbb{SL}_2$, then $G = SL_2(F)$, and Corollary 4.4 of \jcite{Sho20} shows $kG = kSL_2(F)$ is coherent.\\
If $\mathbb{G} = \mathbb{PGL}_2$, then $G = PGL_2(F)$ is isomorphic to a quotient of $SL_2(F)$, namely ${ G \cong \faktor{SL_2(F)}{H} }$, where $H = \{I, -I\}$ is the centre of $SL_2(F)$. Then, by Proposition \ref{surjection of groups prop}, the augmented Iwasawa algebra $kG \cong kSL_2(F)/\epsilon(H)$. The two-sided ideal $\epsilon(H)$ is certainly finitely-generated as a left ideal, because $H$ is finite. Thus, $kG$ is coherent by Lemma \ref{finite ideal quotient coherence lemma}, since $kSL_2(F)$ is coherent.  $\square$\\
\end{jpf}

\subsection{Coherence for the general linear group}

Notice that applying Proposition \ref{closed subgroup coherence prop} and Theorem \ref{Main result semisimple algebraic group non-coherence thm} to the inclusion of the special linear group $SL_n(F)$ into the general linear group $GL_n(F)$ gives the following.

\begin{jthm} \label{general linear group non-coherence thm} \h
If $n \geq 3$, then the augmented Iwasawa algebra $kGL_n(F)$ is not coherent.\\
\end{jthm}

In the rest of this subsection, we improve Theorem \ref{general linear group non-coherence thm} to a full characterisation. Firstly notice that, as in the second Example of subsection \ref{Definitions}, the augmented Iwasawa algebra
\[ kGL_1(F) \cong kF^\times \cong k\mathcal{O}_F^\times[X, X^{-1}] \]
is a Noetherian ring, and hence coherent. We now prove that $GL_2(F)$ also has a coherent augmented Iwasawa algebra.\\
We reuse much of the notation from Section 5 of \jcite{Sho20}. In particular, let $G=GL_2(F)$ throughout the remainder of the subsection.

\begin{jdef} \h
Let
\[ K_1 = GL_2(\mathcal{O}_F), \quad \alpha = \begin{pmatrix} 1 & 0 \\ 0 & \pi \end{pmatrix} \in G, \quad K_2 = \alpha K_1 \alpha^{-1}, \quad  z =  \begin{pmatrix} \pi & 0 \\ 0 & \pi \end{pmatrix} \in G, \]
and
\[ G^0 = \mathrm{Ker}\big(v_F \circ \mathrm{det}: GL_2(F) \rightarrow \mathbb{Z} \big), \]
where $\mathrm{det}: GL_2(F) \rightarrow F^\times$ is the matrix determinant.\\
\end{jdef}

Because $G/G^0 \cong \mathbb{Z}$ is a discrete group, $G^0$ is an open subgroup of $G$. Therefore any compact open subgroup $H \leq G^0$ is also a compact open subgroup of $G$, and so the augmented Iwasawa algebra of $G$ is given by
\[ kG = \bigoplus_{g \in H \backslash G} kH \otimes g = \bigoplus_{\substack{g \in H \backslash G^0\\ g' \in G^0 \backslash G}} kH \otimes gg' = \bigoplus_{g' \in G^0 \backslash G} kG^0 g. \]
In fact, a set of representatives for the cosets $G^0 \backslash G$ is
\[  \{ z^n \mid n \in \mathbb{Z} \} \cup \{ z^n\alpha \mid n \in \mathbb{Z} \}, \]
by considering determinants. Therefore
\[ kG = \bigoplus_{n \in \mathbb{Z}} kG^0z^n \oplus \bigoplus_{n \in \mathbb{Z}} kG^0z^n\alpha. \]

\begin{jdef} \h
Let $G' = \langle G^0, z \rangle \leq G$ be the subgroup of $G$ generated by $G^0$ and $z$. \\
\end{jdef}

Then $G'$ is also an open subgroup of $G$. In fact, $G' = (v_F \circ \mathrm{det})^{-1}(2\mathbb{Z})$, so $G'$ is of index 2 in $G$, and 
\[ kG = kG' \oplus kG' \alpha, \]
thus $kG$ is a free left $kG'$-module of rank 2.

\begin{jthm} \h \label{kG' coherent thm}
The augmented Iwasawa algebra $kG'$ is coherent.\\
\end{jthm}

We prove Theorem \ref{kG' coherent thm} below by considering amalgmated products of rings. The coherence of $kGL_2(F)$ then follows straightforwardly.

\begin{jthm} \h \label{GL2 coherence thm}
The augmented Iwasawa algebra of the general linear group $GL_2(F)$ is a coherent ring.\\
\end{jthm}
\begin{jpf} \h
Clearly $kG = kGL_2(F)$ is a finitely-presented left $kG'$-module, since it is free of rank 2. By Theorem \ref{kG' coherent thm}, $kG'$ is a left coherent ring, thus by Lemma \ref{fp extension coherence lemma}, $kG$ is a left coherent ring. $\square$\\
\end{jpf}

This theorem implies an important consequence for the category of finitely-presented smooth representations of $GL_2(F)$.

\begin{jcor} \h \label{GL2 fp smooth reps cor}
The category of finitely-presented smooth representations of $GL_2(F)$ over $k$ is an abelian category.\\
\end{jcor}
\begin{jpf} \h
Let $G=GL_2(F)$ and $J \leq G$ be an open pro-$p$ subgroup. Then $kG$ is a coherent left $kG$-module, and $kG/\epsilon_G(J)$ is a finitely-presented quotient. So by Theorems 2.1.2 and 2.2.1 of \jcite{Glaz89}, $kG/\epsilon_G(J)$ is also a left coherent $kG$-module. Hence, by Corollary \ref{fp smooth equivalences cor}, the category of finitely-presented smooth representations of $GL_2(F)$ is abelian. $\square$\\
\end{jpf}

Corollary \ref{GL2 fp smooth reps cor} improves Corollary 5.2 of \jcite{Sho20} by removing the $Z$-finiteness assumption on the representations of $GL_2(F)$.\\

Notice that the proof of Theorem \ref{GL2 coherence thm} does not require $k$ to be a field of characteristic $p$ -- in particular Theorem \ref{GL2 coherence thm} and Corollary \ref{GL2 fp smooth reps cor} hold when $k$ is the ring of integers of a finite extension of $\mathbb{Q}_p$. This proves the $GL_2(F)$ case of Conjecture 6.1.4 of \jcite{EmGeHe23}.\\

It now remains to prove Theorem \ref{kG' coherent thm}, which again requires an intermediate result.

\begin{jprop} \h \label{Laurent amalgamation prop}
Let $D = B *_{A} C$ be an amalgamated product of rings. The amalgamated product of Laurent polynomial rings
\[ B[X, X^{-1}] *_{A[X, X^{-1}]} C[X, X^{-1}] \cong D[X, X^{-1}]. \]
\end{jprop}
\begin{jpf} \h
We show that $D[X, X^{-1}]$ satisfies the universal property of being a pushout.\\
Let $D$ be the amalgamated product of $B,C$ over $A$ via the ring homomorphisms $f: A \rightarrow B$, $g: A \rightarrow C$, giving the following commutative diagram.
\[
\begin{tikzcd}
D                   & B \arrow[l, "j_B"]               \\
C \arrow[u, "j_C"'] & A \arrow[l, "g"'] \arrow[u, "f"]
\end{tikzcd}
\]
Let $F, G, J_B, J_C$ be the appropriate ring homomorphisms between $A[X,X^{-1}]$, $B[X,X^{-1}]$, $C[X,X^{-1}]$, $D[X,X^{-1}]$ which each map $X \mapsto X$,
\[
\begin{tikzcd}
{D[X,X^{-1}]}                    & {B[X,X^{-1}]} \arrow[l, "J_B"]                \\
{C[X, X^{-1}]} \arrow[u, "J_C"'] & {A[X, X^{-1}]} \arrow[l, "G"'] \arrow[u, "F"]
\end{tikzcd}
\]
and let $i_A: A \rightarrow A[X,X^{-1}]$ be the natural inclusion, similarly defining $i_B, i_C, i_D$.\\
Now, let $Q$ be a ring, with ring homomorphisms $\theta: B[X, X^{-1}] \rightarrow Q$, $\phi: C[X, X^{-1}] \rightarrow Q$ such that $\theta \circ F = \phi \circ G$, as shown in the following diagram.
\[
\begin{tikzcd}
Q &                                                      &                                                       \\
  & {D[X,X^{-1}]} \arrow[lu, "u" description, dashed]    & {B[X,X^{-1}]} \arrow[l, "J_B"] \arrow[llu, "\theta"'] \\
  & {C[X, X^{-1}]} \arrow[u, "J_C"'] \arrow[luu, "\phi"] & {A[X, X^{-1}]} \arrow[l, "G"'] \arrow[u, "F"]        
\end{tikzcd}
\]
We show that there exists a unique ring homomorphism $u:D[X, X^{-1}] \rightarrow Q$ such that the diagram commutes.\\
Notice that $F \circ i_A = i_B \circ f$, $J_B \circ i_B = i_D \circ j_B$, similarly for $G, J_C$, and that $\theta \circ i_B = \theta|_B$, $\phi \circ i_C = \phi|_C$. Thus we have a commutative diagram,
\[
\begin{tikzcd}
Q &                                            &                                              \\
  & D \arrow[lu, "v" description]              & B \arrow[l, "j_B"] \arrow[llu, "\theta|_B"'] \\
  & C \arrow[u, "j_C"'] \arrow[luu, "\phi|_C"] & A \arrow[l, "g"'] \arrow[u, "f"]            
\end{tikzcd}
\]
and there exists a unique $v: D \rightarrow Q$ such that $v \circ j_B = \theta|_B$, $v \circ j_C = \phi|_C$, by the universal property of pushout for $D$.
Now, $\theta(X) = \phi(X)$ since $\theta \circ F = \phi \circ G$, so $\theta(X)$ commutes with all elements of $v \circ j_B(B)$ and $v \circ j_C(C)$. Since $D$ is generated by $j_B(B), j_C(C)$, it follows that $\theta(X)$ commutes with all elements of $v(D)$. Moreover, $\theta(X)$ is invertible with inverse $\theta(X^{-1})$.\\
Hence, there is a unique ring homomorphism $u: D[X, X^{-1}] \rightarrow Q$ satisfying $u \circ i_D = v$ and $u(X) = \theta(X)$. Then, 
\[ u \circ J_B \circ i_B = u \circ i_D \circ j_B = v \circ j_B = \theta|_B, \quad u(X) = \theta(X), \]
so $u \circ J_B = \theta$, and
\[ u \circ J_C \circ i_C = u \circ i_D \circ j_C = v \circ j_C = \phi|_C, \quad u(X) = \phi(X), \]
so $u \circ J_C = \phi$.\\
Let $u': D[X, X^{-1}] \rightarrow Q$ be any ring homomorphism such that $u' \circ J_B = \theta$ and $u' \circ J_C = \phi$. Then $u'(X) = \theta(X) = \phi(X) = u(X)$. Also
\[ (u' \circ i_D) \circ j_B = u' \circ (i_D \circ j_B) = u' \circ J_B \circ i_B = \theta|_B = v \circ j_B, \]
and similarly $(u' \circ i_D) \circ j_C = v \circ j_C$. Since $D$ is generated by $j_B(B), j_C(C)$, it follows $u' \circ i_D = v$. Therefore $u' = u$.\\
Therefore, there is a unique ring homomorphism $u: D[X, X^{-1}] \rightarrow Q$ making the diagram commute. So $D[X, X^{-1}]$ satisfies the universal property of pushout, as required. $\square$\\
\end{jpf}

\begin{jpfof}[Theorem \ref{kG' coherent thm}] \h
Let $H_1, H_2, H'$ be the closed subgroups of $G'$ generated by $z$ and $K_1, K_2, K_1 \cap K_2$, respectively:
\[ H_1 = \langle K_1, z \rangle, \quad H_2 = \langle K_2, z \rangle, \quad H' = \langle K_1 \cap K_2, z \rangle. \]
Then $K_1$ is a compact open subgroup of $H_1$, $z$ commutes with all elements of $K_1$, and $z^n \in K_1$ only if $n=0$. It follows that the augmented Iwasawa algebra $kH_1 \cong kK_1[X, X^{-1}]$. Also $kK_1$ is Noetherian, hence $kH_1$ is Noetherian by Hilbert's basis theorem. Similarly, the augmented Iwasawa algebras $kH_2 \cong kK_2[X, X^{-1}]$ and $kH' \cong k(K_1 \cap K_2)[X, X^{-1}]$ are Noetherian. Moreover the natural inclusion maps from $H'$ to $H_1, H_2$ send $X$ to $X$.\\
By Proposition \ref{Laurent amalgamation prop}, it follows that
\[ kH_1 *_{kH'} kH_2 \cong kK_1[X, X^{-1}] *_{k(K_1 \cap K_2)[X, X^{-1}]} kK_2[X, X^{-1}] \cong \left(kK_1 *_{k(K_1 \cap K_2)} kK_2\right) [X, X^{-1}] . \]
Now, by Proposition 4.2 of \jcite{Sho20} and Theorem II.3 of \jcite{Ser80}, (as shown during the proof of Corollary 5.2 of \jcite{Sho20}),
\[kG^0 \cong kK_1 *_{k(K_1 \cap K_2)} kK_2. \]
Moreover, as $z$ is a central element, and no power lies in $G^0$, we have a natural isomorphism $kG' \cong kG^0[X, X^{-1}]$. It follows that
\[ kG' \cong kG^0[X, X^{-1}] \cong kH_1 *_{kH'} kH_2,\]
with the natural maps from $kH_1, kH_2, kH'$ corresponding to inclusions of $H_1, H_2, H'$ into $G'$. In particular, by Theorem \ref{non-compact faithfully flat thm}, $kG'$ is a flat module over each of the Noetherian rings $kH_1, kH_2, kH'$. Thus, by Theorem 12 of \jcite{Abe82}, $kG'$ is a left coherent ring. $\square$\\
\end{jpfof}

Thus we have shown that the augmented Iwasawa algebra $kGL_2(F)$ is a coherent ring. Consequently, by Theorem \ref{general linear group non-coherence thm}, we know precisely when the augmented Iwasawa algebra of each general linear group is coherent. 

\notinsubfile{
\generallineargroups*
}
\onlyinsubfile{
\begin{jcor} \h \label{general linear group coherence cor}
The augmented Iwasawa algebra of $GL_n(F)$ is coherent if and only if $n \leq 2$.\\
\end{jcor}
}

\bibliographystyle{coherencejt.bst}
\bibliography{0_Whole_paper.bib}

\begin{thebibliography}{EGH23}

\bibitem[{\AA}be82]{Abe82}
Hans {\AA}berg, \emph{Coherence of Amalgamations}, Journal of Algebra,
  volume~78, no.~2:pp. 372--385 (1982).

\bibitem[Ard21]{Ard21}
Konstantin Ardakov, \emph{{Equivariant $\mathcal{D}$-modules on rigid analytic
  spaces}}, Astérisque, volume 423 (2021).

\bibitem[Bor91]{Bo91}
Armand Borel, \emph{{Linear Algebraic Groups}}, number 126 in Graduate Texts in
  Mathematics, Springer-Verlag, second edition (1991).

\bibitem[Bru66]{Bru66}
Armand Brumer, \emph{{Pseudocompact algebras, profinite groups and class
  formations}}, Bulletin of the American Mathematical Society, volume~72,
  no.~2:pp. 321--324 (1966).

\bibitem[CD03]{ChaDo03}
M.J. Chasco, X.~Dominguez, \emph{{Topologies on the direct sum of topological
  Abelian groups}}, Topology and its Applications, volume 133:pp. 209--223
  (2003).

\bibitem[EGH23]{EmGeHe23}
Matthew Emerton, Toby Gee, Eugen Hellmann, \emph{An introduction to the
  categorical $p$-adic {L}anglands program}, arXiv preprint
  \href{http://arXiv.org/abs/2210.01404}{arXiv:2210.01404} (2023).

\bibitem[EH76]{EdHa76}
David~A. Edwards, Harold~M. Hastings, \emph{Čech and Steenrod homotopy
  theories with applications to geometric topology}, number 542 in Lecture
  Notes in Mathematics, Springer-Verlag (1976).

\bibitem[Eme08]{Eme08}
Matthew Emerton, \emph{{On a class of coherent rings, with applications to the
  smooth representation theory of $GL_2(\mathbb{Q}_p)$ in characteristic $p$}},
  preprint available at
  \url{http://math.uchicago.edu/~emerton/pdffiles/frob.pdf} (2008).

\bibitem[Eme10]{Eme10}
Matthew Emerton, \emph{{Ordinary parts of admissible representations of
  $p$-adic reductive groups {I. Definition and first properties}}},
  Ast{\'e}risque, volume 331:pp. 355--402 (2010).

\bibitem[Gla89]{Glaz89}
Sarah Glaz, \emph{Commutative Coherent Rings}, number 1371 in Lecture Notes in
  Mathematics, Springer-Verlag (1989).

\bibitem[GW04]{GoWa04}
K.R. Goodearl, R.B. Warfield, \emph{{An Introduction to Noncommutative
  Noetherian Rings}}, number~61 in London Mathematical Society Student Texts,
  Cambridge University Press (2004).

\bibitem[Har66]{Har66}
Morton~E. Harris, \emph{{Some Results on Coherent Rings}}, Proceedings of the
  American Mathematical Society, volume~17:pp. 474--479 (1966).

\bibitem[Har67]{Har67}
Morton~E. Harris, \emph{{Some Results on Coherent Rings II}}, Glasgow
  Mathematical Journal, volume~8:pp. 123--126 (1967).

\bibitem[Hu12]{Hu12}
Yongquan Hu, \emph{{Diagrammes canoniques et représentations modulo p de
  $GL_2(F)$}}, Journal of the Institute of Mathematics of Jussieu (2012).

\bibitem[Koh17]{Koh17}
Jan Kohlhaase, \emph{{Smooth duality in natural characteristic}}, Advances in
  Mathematics, volume 317:pp. 1--49 (2017).

\bibitem[Lan02]{Lang02}
Serge Lang, \emph{Algebra}, number 211 in Graduate Texts in Mathematics,
  Springer-Verlag, third edition (2002).

\bibitem[MT11]{MaTe11}
Gunter Malle, Donna Testerman, \emph{{Linear Algebraic Groups and Finite Groups
  of Lie Type}}, number 133 in Cambridge Studies in Advanced Mathematics,
  Cambridge University Press (2011).

\bibitem[Rot09]{Rot09}
Joseph~J Rotman, \emph{{An Introduction to Homological Algebra}}, Universitext,
  Springer, second edition (2009).

\bibitem[Sch02]{Schn02}
Peter Schneider, \emph{Nonarchimedean Functional Analysis}, Springer Monographs
  in Mathematics, Springer (2002).

\bibitem[Sch11]{Schn11}
Peter Schneider, \emph{{$p$-Adic Lie Groups}}, number 344 in Grundlehren der
  mathematischen Wissenschaften, Springer (2011).

\bibitem[Sch15]{Schr15}
Benjamin Schraen, \emph{{Sur la présentation des représentations
  supersingulières de $GL_2(F)$}}, Journal für die Reine und Angewandte
  Mathematik, volume 704:pp. 187--208 (2015).

\bibitem[Ser77]{Ser77}
Jean-Pierre Serre, \emph{Linear Representations of Finite Groups}, number~42 in
  Graduate Texts in Mathematics, Springer-Verlag, translated from the second
  French edition by Leonard L. Scott (1977).

\bibitem[Ser79]{Ser79}
Jean-Pierre Serre, \emph{Local Fields}, number~67 in Graduate Texts in
  Mathematics, Springer-Verlag, translated from the French by Marvin Jay
  Greenberg (1979).

\bibitem[Ser80]{Ser80}
Jean-Pierre Serre, \emph{Trees}, Springer-Verlag, translated from the French by
  John Stillwell (1980).

\bibitem[Sho20]{Sho20}
Jack Shotton, \emph{{The category of finitely presented smooth mod $p$
  representations of $GL_2(F)$}}, Documenta Mathematica, volume~25:pp. 143--157
  (2020).

\bibitem[Vig11]{Vig11}
Marie-France Vignéras, \emph{{Le foncteur de Colmez pour $GL(2,F)$}}, in
  \emph{Arithmetic Geometry and Automorphic Forms}, number~19 in Advanced
  Lectures in Mathematics, International Press, Somerville, MA, pp. 531--557
  (2011).

\bibitem[Wei94]{Wei94}
Charles~A. Weibel, \emph{{An introduction to homological algebra}}, number~38
  in Camridge Studies in Advanced Mathematics, Cambridge University Press
  (1994).

\bibitem[Wil98]{Wil98}
John~S. Wilson, \emph{{Profinite Groups}}, number~19 in London Mathematical
  Society Monographs, New Series, Clarendon Press (1998).

\bibitem[Wu21]{Wu21}
Zhixiang Wu, \emph{A note on presentations of supersingular representations of
  ${GL_2(F)}$}, Manuscripta Mathematica, volume 165, no.~3:pp. 583--596 (2021).

\end{thebibliography}

\end{document}